\documentclass[12pt,reqno]{amsart}

\usepackage{diagrams,amssymb,amsmath}

\newcommand{\mynote}[1]{}

\newcommand{\newsect}[1]{\vskip7mm\section{#1}}
\newcommand{\newsub}[1]{\vskip4mm\stepcounter{subsection}%
    \noindent\textbf{\underbar{\arabic{section}.\arabic{subsection}\quad 
    #1}}\quad}
\diagramstyle[height=2em,width=2em,midshaft,labelstyle=\scriptstyle]
\newarrow{To}----{->}
\newarrow{Epi}----{->>}
\newarrow{Mono}>---{->}
\newarrow{Iso}>---{->>}
\newarrow{Mapsto}|---{->}
\newarrow{Igual}=====
\newarrow{Dashto}{}{dash}{}{dash}{->}

\newcommand{\Z}{{\mathbb Z}}

\newcommand{\F}{{\mathbb F}}

\newcommand{\pcom}{{}_{p}^{\wedge}}

\newcommand{\zploc}{\Z_{(p)}}

\DeclareMathAlphabet\EuR{U}{eur}{m}{n}
\SetMathAlphabet\EuR{bold}{U}{eur}{b}{n}

\newcommand{\Res}{\operatorname{Res}\nolimits}

\newcommand{\Ext}{\operatorname{Ext}\nolimits}
\newcommand{\Inj}{\operatorname{Inj}\nolimits}

\newcommand{\defeq}{\overset{\text{\textup{def}}}{=}}
\newcommand{\gen}[1]{{\langle}#1{\rangle}}

\renewcommand{\:}{\colon}

\newcommand{\calb}{\mathcal{B}}
\newcommand{\calc}{\mathcal{C}}

\newcommand{\cale}{\mathcal{E}}
\newcommand{\calf}{\mathcal{F}}

\newcommand{\calk}{\mathcal{K}}
\newcommand{\call}{\mathcal{L}}

\newcommand{\calp}{\mathcal{P}}
\newcommand{\calq}{\mathcal{Q}}

\newcommand{\calz}{\mathcal{Z}}

\newcommand{\calfq}{\calf^q}
\newcommand{\calfc}{\calf^c}
\newcommand{\callq}{\call^q}

\newcommand{\orb}{\mathcal{O}}

\newcommand{\nv}[1]{|#1|} 

\newcommand{\curs}{\EuR}
\newcommand{\Ab}{\curs{Ab}}

\renewcommand{\mod}{\mbox{-}\curs{mod}}

\newcommand{\widebar}[1]{\overset{\mskip3mu\hrulefill\mskip3mu}{#1}
                \vphantom{#1}}
\newcommand{\sminus}{\smallsetminus}

\newcommand{\nsg}{\vartriangleleft}
\newcommand{\Id}{\textup{Id}}
\newcommand{\incl}{\operatorname{incl}\nolimits}
\newcommand{\proj}{\operatorname{proj}\nolimits}
\newcommand{\Inn}{\textup{Inn}}

\newcommand{\op}{^{\textup{op}}}
\newcommand{\isotyp}{_{\textup{typ}}}      

\newcommand{\fus}{_{\textup{fus}}}

\let\oldcirc=\circ
\renewcommand{\circ}{\mathchoice
    {\mathbin{\scriptstyle\oldcirc}}{\mathbin{\scriptstyle\oldcirc}}
    {\mathbin{\scriptscriptstyle\oldcirc}}
    {\mathbin{\scriptscriptstyle\oldcirc}}}

\newcommand{\hclim}[1]{\setbox1=\hbox{\rm hocolim}
    \setbox2=\hbox to \wd1{\rightarrowfill} \ht2=0pt \dp2=-1pt
    \mathop{\vtop{\baselineskip=5pt\box1\box2}}
    _{#1}}

\newcommand{\invlim}[1]{\higherlim{#1}{}}

\newcommand{\map}{\operatorname{Map}\nolimits}
\renewcommand{\hom}{\operatorname{Hom}\nolimits}
\newcommand{\Hom}{\operatorname{Hom}\nolimits}
\newcommand{\rep}{\operatorname{Rep}\nolimits}
\newcommand{\Rep}{\operatorname{Rep}\nolimits}
\newcommand{\Iso}{\operatorname{Iso}\nolimits}
\newcommand{\Aut}{\operatorname{Aut}\nolimits}
\newcommand{\Out}{\operatorname{Out}\nolimits}

\newcommand{\Ob}{\operatorname{Ob}\nolimits}
\newcommand{\Mor}{\operatorname{Mor}\nolimits}
\newcommand{\Ker}{\operatorname{Ker}\nolimits}
\renewcommand{\ker}{\operatorname{Ker}\nolimits}

\renewcommand{\Im}{\operatorname{Im}\nolimits}

\newcommand{\q}[1]{\widehat{#1}}

\newcommand{\longleft}[1]{\;{\leftarrow%
\count255=0 \loop \mathrel{\mkern-6mu}%
    \relbar\advance\count255 by1\ifnum\count255<#1\repeat}\;}
\newcommand{\longright}[1]{\;{\count255=0 \loop \relbar\mathrel{\mkern-6mu}%
    \advance\count255 by1\ifnum\count255<#1\repeat\rightarrow}\;}
\newcommand{\Right}[2]{\overset{#2}{\longright#1}}
\newcommand{\RIGHT}[3]{\mathrel{\mathop{\kern0pt\longright#1}
        \limits^{#2}_{#3}}}

\newcommand{\LEFT}[3]{\mathrel{\mathop{\kern0pt\longleft#1}\limits^{#2}_{#3}}
}
\newcommand{\dRIGHT}[3]{\mathrel{%
   \mathop{\vcenter{\baselineskip=0pt\hbox{$\kern0pt\longright#1$}%
   \hbox{$\kern0pt\longright#1$}}}\limits^{#2}_{#3}}}
\newcommand{\LRIGHT}[3]{\mathrel{%
   \mathop{\vcenter{\baselineskip=0pt\hbox{$\kern0pt\longleft#1$}%
   \hbox{$\kern0pt\longright#1$}}}\limits^{#2}_{#3}}}
\newcommand{\RLEFT}[3]{\mathrel{%
   \mathop{\vcenter{\baselineskip=0pt\hbox{$\kern0pt\longright#1$}%
   \hbox{$\kern0pt\longleft#1$}}}\limits^{#2}_{#3}}}
\newcommand{\onto}[1]{\;{\count255=0 \loop \relbar\joinrel
    \advance\count255 by1
    \ifnum\count255<#1 \repeat \twoheadrightarrow}\;}
\newcommand{\Onto}[2]{\overset{#2}{\onto#1}}

\newtheorem{Thm}{Theorem}[section]
\newtheorem{Prop}[Thm]{Proposition}
\newtheorem{Cor}[Thm]{Corollary}
\newtheorem{Lem}[Thm]{Lemma}

\newtheorem{Defi}[Thm]{Definition}
\newtheorem{Th}{Theorem}

\newtheorem{Ex}[Thm]{Example}

\newcommand{\longline}{\bigskip\centerline{\hbox to 5cm{\hrulefill}}\bigskip}

\newcommand{\hfocal}[2]{O^p_{#1}(#2)}  

\newcommand{\cj}[3][]{\def\test{#1}\def\tst{x}\ifx\test\tst{#2^{-1}#3#2}
	\else{#2#3#2^{-1}}\fi}
\newcommand{\cjup}[3][]{\def\test{#1}\def\tst{x}\def\tset{-}
	\ifx\test\tst{{}^{#2^{-1}}\!#3} 
	\else{\ifx\test\tset{{}^{#2}#3} \else{{}^{#2}\!#3} \fi} \fi}

\renewcommand{\b}[1]{\check{#1}}

\newcommand{\higherlim}[2]{\displaystyle\setbox1=\hbox{\rm lim}
	\setbox2=\hbox to \wd1{\leftarrowfill} \ht2=0pt \dp2=-1pt
	\setbox3=\hbox{$\scriptstyle{#1}$}
	\def\test{#1}\ifx\test\empty
	\mathop{\mathop{\vtop{\baselineskip=5pt\box1\box2}}}\nolimits^{#2}
	\else
	\ifdim\wd1<\wd3
	\mathop{\hphantom{^{#2}}\vtop{\baselineskip=5pt\box1\box2}^{#2}}_{#1}
	\else
	\mathop{\mathop{\vtop{\baselineskip=5pt\box1\box2}}_{#1}}%
	\nolimits^{#2}
	\fi\fi}

\newcommand{\calh}{\mathcal{H}}

\newcommand{\SFL}[1][]{(S#1,\calf#1,\call#1)}

\newcommand{\sset}{\mathfrak{Sub}}

\renewcommand{\labelenumi}{\textup{(\alph{enumi})}}%

\newcommand{\homf}{\Hom_{\calf}}
\newcommand{\repf}{\rep_{\calf}}
\newcommand{\autf}{\Aut_{\calf}}
\newcommand{\outf}{\Out_{\calf}}
\newcommand{\isof}{\Iso_{\calf}}

\newcommand{\sylp}[1]{\textup{Syl}_p(#1)}

\renewenvironment{enumerate}{\begin{list}%
{\labelenumi}
{\usecounter{enumi}%
\setlength{\itemindent}{0pt}%
\settowidth{\labelwidth}{\labelenumi}%
\addtolength{\labelwidth}{\labelsep}%
\setlength{\leftmargin}{\labelsep}%
\addtolength{\leftmargin}{\labelwidth}%
\setlength{\listparindent}{0pt}%
\setlength{\itemsep}{6pt}%
\setlength{\parsep}{0pt}%
\setlength{\topsep}{6pt}%
}}{\end{list}}

\newenvironment{enumtwo}{\begin{list}%
{\labelenumii}
{\usecounter{enumii}%
\setlength{\itemindent}{0pt}%
\settowidth{\labelwidth}{\labelenumii}%
\addtolength{\labelwidth}{\labelsep}%
\setlength{\leftmargin}{\labelsep}%
\addtolength{\leftmargin}{\labelwidth}%
\setlength{\listparindent}{0pt}%
\setlength{\itemsep}{6pt}%
\setlength{\parsep}{0pt}%
\setlength{\topsep}{6pt}%
}}{\end{list}}

\newcommand{\hilim}[4][]{\def\test{#1}\def\tst{-}
\ifx\test\tst{\setbox1=\hbox{\rm lim}
\setbox2=\hbox to \wd1{\leftarrowfill} \ht2=0pt \dp2=-1pt
\mathop{\vtop{\baselineskip=5pt\box1\box2}}\nolimits^{#3}(#4)}
\else{\higherlim{#2}{#3}(#4)}\fi}
\let\2=\partial
\def\Syl{\textup{Syl}}

\newcommand{\dirlim}[1]{\setbox1=\hbox{\rm colim}
	\setbox2=\hbox to \wd1{\rightarrowfill} \ht2=0pt \dp2=-1pt
	\mathop{\vtop{\baselineskip=5pt\box1\box2}}
	_{#1}}

\renewcommand{\onto}[1]{\;{\count255=0 \loop \relbar\mathrel{\mkern-6mu}%
	\advance\count255 by1 \ifnum\count255<#1\repeat\twoheadrightarrow}\;}

\def\beq#1\eeq{\begin{equation*}#1\end{equation*}}

\newcommand{\oppf}{O^{p'}_*(\calf)}
\newcommand{\oppfc}{O^{p'}_*(\calf)^c}

\newcommand{\opf}{O^{p}_*(\calf)}
\newcommand{\opfc}{O^{p}_*(\calf)^c}

\newcommand{\Op}[1]{O^p(#1)}
\newcommand{\Opp}[1]{O^{p'}(#1)}
\newcommand{\Oppf}{\Opp{\calf}}
\newcommand{\Opf}{\Op{\calf}}
\newcommand{\Oppl}{\Opp{\call}}
\newcommand{\Opl}{\Op{\call}}
\newcommand{\gpp}[1]{\Gamma_{p'}(#1)} 
\newcommand{\gp}[1]{\Gamma_{p}(#1)} 
\newcommand{\gpf}{\gp{\calf}}
\newcommand{\gppf}{\gpp{\calf}}
\newcommand{\lf}[2]{\def\test{#2}\def\tst{}\ifx\test\tst\calf_{1}%
        \else\calf_{#2}\fi}      
\newcommand{\lfx}[2]{\lf{#1}{#2}^{\bullet}}
\renewcommand{\ll}[2]{\def\test{#2}\def\tst{}\ifx\test\tst\call_{1}%
        \else\call_{#2}\fi}      
\newcommand{\llx}[2]{\ll{#1}{#2}^{\bullet}}

\title{Extensions of $p$-local finite groups}

\author{C. Broto}
\address{Departament de Matem\`atiques, Universitat Aut\`onoma de 
Barcelona, E--08193 Bellaterra, Spain}
\email{broto@mat.uab.es}
\thanks{C. Broto is partially supported by MCYT grant BFM2001--2035}

\author{N. Castellana}
\address{Departament de Matem\`atiques, Universitat Aut\`onoma de 
Barcelona, E--08193 Bellaterra, Spain}
\email{natalia@mat.uab.es}
\thanks{N. Castellana is partially supported by MCYT grant BFM2001--2035}

\author{J. Grodal}
\address{Department of Mathematics, University of Chicago, Chicago, IL 
60637, USA}
\email{jg@math.uchicago.edu}
\thanks{J. Grodal is partially supported by NSF grants DMS-0104318 and 
DMS-0354633}

\author{R. Levi}
\address{Department of Mathematical Sciences, University of Aberdeen, 
Meston Building 339, Aberdeen AB24 3UE, U.K.}
\email{ran@maths.abdn.ac.uk}
\thanks{R. Levi is partially supported by EPSRC grant GR/M7831.}

\author{B. Oliver}
\address{LAGA, Institut Galil\'ee, Av. J-B Cl\'ement, 93430 
Villetaneuse, France}
\email{bob@math.univ-paris13.fr}
\thanks{B. Oliver is partially supported by UMR 7539 of the CNRS}



\begin{document}

\begin{abstract}  
A $p$-local finite group consists of a finite $p$-group $S$, together with 
a pair of categories which encode ``conjugacy'' relations among subgroups 
of $S$, and which are modelled on the fusion in a Sylow $p$-subgroup of a 
finite group.  It contains enough information to define a classifying 
space which has many of the same properties as $p$-completed  classifying 
spaces of finite groups.  In this paper, we study and classify extensions 
of $p$-local finite groups, and also compute the fundamental group of the 
classifying space of a $p$-local finite group. 
\end{abstract}

\maketitle

\markboth {\hfill \uppercase{C. Broto, N. Castellana, J. Grodal, R. Levi, 
\& B. Oliver}\hfill} {\hfill\uppercase{Extensions of $p$-local finite 
groups}\hfill} 


A $p$-local finite group consists of a finite $p$-group S, together with a 
pair of categories $(\calf,\call)$, of which $\calf$ is modeled on the 
conjugacy (or fusion) in a Sylow subgroup of a finite group. The category 
$\call$ is essentially an extension of $\calf$ and contains just enough 
extra information so that its $p$-completed nerve has many of the same 
properties as $p$-completed classifying spaces of finite groups.  We 
recall the precise definitions of these objects in Section 1, and refer to 
\cite{BLO2} and \cite{bcglo1} for motivation for their study.

In this paper, we study extensions of saturated fusion systems and of 
$p$-local finite groups.  This is in continuation of our more general 
program of trying to understand to what extent properties of finite groups 
can be extended to properties of $p$-local finite groups, and to shed light 
on the question of how many (exotic) $p$-local finite groups there are. 
While we do not get a completely general theory of extensions of one 
$p$-local finite group by another, we do show how certain types of 
extensions can be described in manner very similar to the situation for 
finite groups.

From the point of view of group theory, developing an extension theory for 
$p$-local finite groups is related to the question of to what extent the 
extension problem for groups is a {\em local} problem, i.e., a problem 
purely described in terms of a Sylow $p$-subgroup and conjugacy
relations inside it. In 
complete generality this is {\em not} the case.  For example, strongly 
closed subgroups of a Sylow $p$-subgroup $S$ of $G$ need not correspond to 
normal subgroups of $G$. However, special cases where this does happen 
include the case of existence of $p$-group quotients (the focal subgroup 
theorems, see \cite[\S\S 7.3--7.4]{Gorenstein}) and central subgroups 
(described via the $Z^*$-theorem of Glauberman \cite{Glauberman}).

From the point of view of homotopy theory, one of the problems which comes 
up when looking for a general theory of extensions of $p$-local finite 
groups is that while an extension of groups $1\to{}K\to\Gamma\to{}G\to1$ 
always induces a (homotopy) fibration sequence of classifying spaces, it 
does not in general induce a fibration sequence of $p$-completed 
classifying spaces.  Two cases where this does happen are those where $G$ 
is a $p$-group, and where the extension is central.  In both of these 
cases, $BH\pcom$ is the homotopy fiber of the map 
$B\Gamma\pcom\Right2{}BG\pcom$.  Thus, also from the point of view of 
homotopy theory, it is natural to study extensions of $p$-local finite 
groups with $p$-group quotient, and to study central extensions of 
$p$-local finite groups.  The third case we study is that of extensions 
with quotient of order prime to $p$; and the case of $p$ and $p'$-group 
quotients to some extent unify to give a theory of extensions with 
$p$-solvable quotient. Recall in this connection that by a previous result 
of ours \cite[Proposition~C]{bcglo1}, solvable $p$-local finite groups all 
come from $p$-solvable groups. In all three of these situations, we 
develop a theory of extensions which parallels the situation for finite 
groups.  

We now describe the contents of this paper in more detail, stating 
simplified versions of our main results on extensions.  Stronger and more 
precise versions of some of these theorems will be stated and proven later.

In Section 3, we construct a very general theory of fusion subsystems 
(Proposition \ref{F:p-solv.quot.}) and linking subsystems (Theorem 
\ref{L:p-p'-quot.}) with quotient a $p$-group or a group of order prime to 
$p$. As a result we get the following theorem (Corollary 
\ref{L:p-solv.quot.}), which for a $p$-local finite group $\SFL$, 
describes a correspondence between covering spaces of the geometric 
realization $|\call|$ and certain $p$-local finite subgroups of $\SFL$.

\begin{Th} \label{ThA}
Suppose that $\SFL$ is a $p$-local finite group.  Then there is a normal 
subgroup $H\nsg\pi_1(|\call|)$ which is minimal among all those whose 
quotient is finite and $p$-solvable.  Any covering space of the geometric 
realization $|\call|$ whose fundamental group contains $H$ is homotopy 
equivalent to $|\call'|$ for some $p$-local finite group $\SFL[']$, where 
$S'\le{}S$ and $\calf'\subseteq{}\calf$.
\end{Th}

Moreover, the $p$-local finite group $\SFL[']$ of Theorem \ref{ThA} can be 
explicitly described in terms of $\call$, as we will explain in Section 3.

In order to use this theorem, it is useful to have ways of finding the 
finite $p$-solvable quotients of $\pi_1(|\call|)$. This can be done by 
iteration, using the next two theorems.  In them, the maximal $p$-group 
quotient of $\pi_1(|\call|)$, and the maximal quotient of order prime to 
$p$, are described solely in terms of the fusion system $\calf$.  

When $G$ is an infinite group, we define $O^p(G)$ and $O^{p'}(G)$ to be 
the intersection of all normal subgroups of $G$ of $p$-power index, or 
index prime to $p$, respectively.  These clearly generalize the usual 
definitions for finite $G$ (but are not the only possible 
generalizations). 

\begin{Th}[Hyperfocal subgroup theorem for $p$-local finite groups] 
\label{ThB}
For a $p$-local finite group $\SFL$, the natural homomorphism
	$$S \Right4{} \pi_1(|\call|)/O^p(\pi_1(|\call|)) \cong 
	\pi_1(|\call|\pcom) $$
is surjective, with kernel equal to 
	$$ O^p_\calf(S) \defeq \bigl\langle g^{-1}\alpha(g) \in S
	\,\big|\, g\in{}P\le{}S,\ \alpha\in O^p(\Aut_\calf(P))  
	\bigr\rangle.$$
\end{Th}

For any saturated fusion system $\calf$ over a $p$-group $S$, we let 
$\oppf\subseteq\calf$ be the smallest fusion subsystem of $\calf$ (in the 
sense of Definition \ref{D:fus.syst.}) which contains all automorphism 
groups $O^{p'}(\autf(P))$ for $P\le{}S$.  Equivalently, $\oppf$ is the 
smallest subcategory of $\calf$ with the same objects, and which contains 
all restrictions of all automorphisms in $\calf$ of $p$-power order.  This 
subcategory is needed in the statement of the next theorem.

\begin{Th} \label{ThC}
For a $p$-local finite group $\SFL$, the natural map 
	$$\Out_\calf(S) \Right4{} \pi_1(|\call|)/O^{p'}(\pi_1(|\call|))$$ 
is surjective, with kernel equal to
	$$ \outf^0(S) \defeq \bigl\langle \alpha\in\outf(S) \,\big|\,
	\alpha|_P\in\Mor_{\oppf}(P,S), \textup{ some $\calf$-centric
	$P\le{}S$} \bigr\rangle. $$
\end{Th}

Theorem \ref{ThB} is proved as Theorem~\ref{pi1(hfg)}, and Theorem 
\ref{ThC} is proved as Theorem~\ref{L:p'-quot.}.

In fact, we give a purely algebraic description of these subsystems of 
``$p$-power index'' or of ``index prime to $p$'' (Definition 
\ref{def-p'-index}), and then show in Sections 4.1 and 5.1 that they in 
fact all arise as finite covering spaces of $|\call|$ (see Theorems 
\ref{L:p-quot.} and \ref{L:p'-quot.}). Similar results were also 
discovered independently by Puig \cite{Puig-ppt}.  Then, in Sections 4.2 
and 5.2, we establish converses to these concerning the extensions of 
a $p$-local finite group, which include the following theorem.

\begin{Th} \label{ThD}
Let $\SFL$ be a $p$-local finite group.  Suppose we are given a fibration 
sequence $|\call|\pcom \rTo E \rTo BG$, where $G$ is a finite $p$-group or 
has order prime to $p$.  Then there exists a $p$-local finite group 
$\SFL[']$ such that $|\call'|\pcom\simeq{}E\pcom$. 
\end{Th} 

This is shown as Theorems \ref{|L|->X->Bpi} and \ref{|L|->X->Bpi'}.  
Moreover, when $G$ is a $p$-group, we give in Theorem \ref{|L|->X->Bpi} an 
explicit algebraic construction of the $p$-local finite group $\SFL[']$.  

Finally, in Section 6, we develop the theory of central extensions of 
$p$-local finite groups. Our main results there (Theorems 
\ref{P:L/A-centext} and \ref{centthm'}) give a more elaborate 
version of the following theorem. Here, the  center of a $p$-local finite 
group $\SFL$ is defined to be the subgroup of elements $x \in Z(S)$ such 
that $\alpha(x)=x$ for all $\alpha\in\Mor(\calf^c)$.

\begin{Th} \label{ThE}
Suppose that $A$ is a central subgroup of a $p$-local finite group $\SFL$. 
Then there exists a canonical quotient $p$-local finite group $\SFL[/A]$, 
and the canonical projection of $|\call|$ onto $|\call/A|$ is a principal 
fibration with fiber $BA$.  

Conversely, for any principal fibration $E \rTo |\call|$ with fiber $BA$, 
where $A$ is a finite abelian $p$-group, there exists a $p$-local finite 
group $(\widetilde{S},\widetilde{\calf},\widetilde{\call})$ such that 
$|\widetilde\call|\simeq{}E$. Furthermore, this correspondence sets up a 
$1-1$ correspondence between central extensions of $\call$ by $A$ and 
elements in $H^2(|\call|;A)$. 
\end{Th} 

One motivation for this study was the question of whether extensions of 
$p$-local finite groups coming from finite groups can produce exotic 
$p$-local finite groups. In the case of central extensions, we are able to 
show that $\SFL$ comes from a finite group if and only if $\SFL[/A]$ comes 
from a group (Corollary \ref{F.exo=>F/A.exo}).  
For the other types of extensions studied in this paper, this is still an 
open question. We have so far failed to produce exotic examples in this 
way, and yet we have also been unable to show that exotic examples cannot 
occur. This question seems to be related to some rather subtle and 
interesting group theoretic issues relating local to global structure; see 
Corollary~\ref{exotic-ext} for one partial result in this direction.

This paper builds on the earlier paper \cite{bcglo1} by the same authors, 
and many of the results in that paper were originally motivated by this 
work on extensions.

The authors would like to thank the University of Aberdeen, Universitat 
Aut\`onoma de Barcelona, Universit\'e Paris 13, and Aarhus Universitet for 
their hospitality.  In particular, the origin of this project goes back to 
a three week period in the spring of 2001, when four of the authors met in 
Aberdeen.


\newsect{A quick review of $p$-local finite groups} 
\label{review}

We first recall the definitions of a fusion system, and a saturated fusion 
system, in the form given in \cite{BLO2}.  For any group $G$, and any pair 
of subgroups $H,K\le{}G$, we set
	$$ N_G(H,K) = \{x\in{}G\,|\, xHx^{-1}\le K \}, $$
let $c_x$ denote conjugation by $x$ ($c_x(g)=xgx^{-1}$), and set 
	$$ \Hom_G(H,K) = \bigl\{c_x\in\Hom(H,K) 
	\,\big|\, x\in{}N_G(H,K) \bigr\} \cong N_G(H,K)/C_G(H). $$
By analogy, we also write
	$$ \Aut_G(H)=\Hom_G(H,H)=\bigl\{c_x\in\Aut(H) 
	\,\big|\, x\in{}N_G(H) \bigr\} \cong N_G(H)/C_G(H). $$

\begin{Defi}[{\cite{Puig} and \cite[Definition~1.1]{BLO2}}] 
\label{D:fus.syst.}
A fusion system over a finite $p$-group $S$ is a category $\calf$, where 
$\Ob(\calf)$ is the set of all subgroups of $S$, and which satisfies the 
following two properties for all $P,Q\le{}S$:
{\renewcommand{\labelenumi}{\hskip-4pt$\bullet$}
\begin{enumerate}
\item  $\Hom_S(P,Q) \subseteq \homf(P,Q) \subseteq \Inj(P,Q)$; and
\item  each $\varphi\in\homf(P,Q)$ is the composite of an isomorphism in
$\calf$ followed by an inclusion.
\end{enumerate}
}
\end{Defi}

The following additional definitions and conditions are needed in order for 
these systems to be very useful.  If $\calf$ is a fusion system over a 
finite $p$-subgroup $S$, then two subgroups $P,Q\le S$ are said to be 
\emph{$\calf$-conjugate} if they are isomorphic as objects of the category 
$\calf$. 

\begin{Defi}[{\cite{Puig}, see \cite[Def.~1.2]{BLO2}}] \label{sat.Frob.}
\renewcommand{\labelenumi}{\hskip-4pt$\bullet$}
Let $\calf$ be a fusion system over a $p$-group $S$.
\begin{enumerate}
\item A subgroup $P\le{}S$ is \emph{fully centralized in $\calf$} if
$|C_S(P)|\ge|C_S(P')|$ for all $P'\le{}S$ which is $\calf$-conjugate to
$P$.
\item A subgroup $P\le{}S$ is \emph{fully normalized in $\calf$} if
$|N_S(P)|\ge|N_S(P')|$ for all $P'\le{}S$ which is $\calf$-conjugate to
$P$.
\item $\calf$ is a \emph{saturated fusion system} if the following
two conditions hold:
\begin{enumtwo}  \renewcommand{\labelenumii}{\textup{(\Roman{enumii})}}
\item For all $P\le{}S$ which is fully normalized in $\calf$, $P$ is fully
centralized in $\calf$ and $\Aut_S(P)\in\sylp{\autf(P)}$.
\item If $P\le{}S$ and $\varphi\in\homf(P,S)$ are such that $\varphi{}P$ is
fully centralized, and if we set
    $$ N_\varphi = \{ g\in{}N_S(P) \,|\, \varphi c_g\varphi^{-1} \in
    \Aut_S(\varphi{}P) \}, $$
then there is $\widebar{\varphi}\in\homf(N_\varphi,S)$ such that
$\widebar{\varphi}|_P=\varphi$.
\end{enumtwo}
\end{enumerate}
\end{Defi}

If $G$ is a finite group and $S\in\sylp{G}$, then by 
\cite[Proposition~1.3]{BLO2}, the category $\calf_S(G)$ defined by 
letting $\Ob(\calf_S(G))$ be the set of all subgroups of $S$ and setting 
$\Mor_{\calf_S(G)}(P,Q)=\Hom_G(P,Q)$ is a saturated fusion system.

An alternative pair of axioms for a fusion system being saturated have been 
given by Radu Stancu \cite{Stancu}.  He showed that axioms (I) and (II) 
above are equivalent to the two axioms:
\begin{enumerate}\renewcommand{\labelenumi}{\textup{(\Roman{enumi}$'$)}}
\item $\Inn(S)\in\sylp{\autf(S)}$.
\item If $P\le{}S$ and $\varphi\in\homf(P,S)$ are such that $\varphi{}P$ is
fully normalized, and if we set
    $$ N_\varphi = \{ g\in{}N_S(P) \,|\, \varphi c_g\varphi^{-1} \in
    \Aut_S(\varphi{}P) \}, $$
then there is $\widebar{\varphi}\in\homf(N_\varphi,S)$ such that
$\widebar{\varphi}|_P=\varphi$.
\end{enumerate}

The following consequence of conditions (I) and (II) above will be needed 
several times throughout the paper. 

\begin{Lem} \label{N->N}
Let $\calf$ be a saturated fusion system over a $p$-group $S$.  Let 
$P,P'\le{}S$ be a pair of $\calf$-conjugate subgroups such that $P'$ is 
fully normalized in $\calf$.  Then there is a homomorphism 
$\alpha\in\homf(N_S(P),N_S(P'))$ such that $\alpha(P)=P'$.
\end{Lem}

\begin{proof}  This is shown in \cite[Proposition A.2(b)]{BLO2}.
\end{proof}

In this paper, it will sometimes be necessary to work with fusion systems 
which are not saturated.  This is why we have emphasized the difference 
between fusion systems, and saturated fusion systems, in the above 
definitions.

We next specify certain collections of subgroups relative to a given 
fusion system. 

\begin{Defi}\label{centric-radical-def}
Let $\calf$ be a fusion system over a finite
$p$-subgroup $S$. 
\renewcommand{\labelenumi}{\hskip-4pt$\bullet$}%
\begin{enumerate}
\item A subgroup $P \leq S$ is \emph{$\calf$-centric} if $C_S(P') =Z(P')$ 
for all $P' \leq S$ which is $\calf$-conjugate to $P$. 
\item A subgroup $P \leq S$ is \emph{$\calf$-radical} if $\Out_\calf(P)$ 
is $p$-reduced; i.e., if $O_p(\outf(P))=1$. 

\item For any $P\le{}S$ which is fully centralized in $\calf$, the 
\emph{centralizer fusion system} $C_\calf(P)$ is the fusion system over 
$C_S(P)$ defined by setting
	$$ \Hom_{C_\calf(P)}(Q,Q') = \bigl\{ \alpha|Q \,\big|\, 
	\alpha\in\homf(QP,Q'P),\ \alpha|_P=\Id_P,\ \alpha(Q)\le Q' \bigr\}. 
	$$
A subgroup $P\le S$ is \emph{$\calf$-quasicentric} if for all $P'\le S$ 
which is $\calf$-conjugate to $P$ and fully centralized in $\calf$, 
$C_\calf(P')$ is the fusion system of the $p$-group $C_S(P')$.  

\item $\calfc\subseteq\calfq\subseteq\calf$ denote the full subcategories 
of $\calf$ whose objects are the $\calf$-centric subgroups, and 
$\calf$-quasicentric subgroups, respectively, of $S$.
\end{enumerate}
\end{Defi}

If $\calf = \calf_S(G)$ for some finite group $G$, then $P\le S$ is 
$\calf$-centric if and only if $P$ is $p$-centric in $G$ (i.e., 
$Z(P)\in\Syl_p(C_G(P))$), and $P$ is $\calf$-radical if and only if 
$N_G(P)/(P{\cdot}C_G(P))$ is $p$-reduced.  Also, $P$ is $\calf$-quasicentric 
if and only if $C_G(P)$ contains a normal subgroup of order prime to $p$ 
and of $p$-power index.

In fact, when working with $p$-local finite groups, it suffices to have a 
fusion system $\calfc$ defined on the centric subgroups of $S$, and 
which satisfies axioms (I) and (II) above for those centric subgroups.  In 
other words, fusion systems defined only on the centric subgroups are 
equivalent to fusion systems defined on all subgroups, as described in the 
following theorem.

\begin{Thm} \label{Fc<=>F} \label{Alp.fusion} \label{centr->sat}
Fix a $p$-group $S$ and a fusion system $\calf$ over $S$.  
\begin{enumerate}  
\item Assume $\calf$ is saturated.  Then each morphism in $\calf$ is a 
composite of restrictions of morphisms between subgroups of $S$ which are 
$\calf$-centric, $\calf$-radical, and fully normalized in $\calf$.  More 
precisely, for each $P,P'\le{}S$ and each $\varphi\in\isof(P,P')$, there 
are subgroups $P=P_0,P_1,\ldots,P_k=P'$, subgroups 
$Q_i\ge\gen{P_{i-1},P_i}$ ($i=1,\ldots,k$) which are $\calf$-centric, 
$\calf$-radical, and fully normalized in $\calf$, and automorphisms 
$\varphi_i\in\autf(Q_i)$, such that $\varphi_i(P_{i-1})=P_i$ for all $i$ 
and $\varphi=\varphi_k\circ\cdots\circ\varphi_1|_P$.

\item Assume conditions (I) and (II) in Definition \ref{sat.Frob.} are 
satisfied for all $\calf$-centric subgroups $P\le{}S$.  Assume also that 
each morphism in $\calf$ is a composite of restrictions of morphisms 
between $\calf$-centric subgroups of $S$.  Then $\calf$ is saturated.
\end{enumerate}
\end{Thm}

\begin{proof}  Part (a) is Alperin's fusion theorem for saturated fusion 
systems, in the form shown in \cite[Theorem A.10]{BLO2}.  Part (b) is a 
special case of \cite[Theorem 2.2]{bcglo1}:  the case where $\mathcal{H}$ 
is the set of all $\calf$-centric subgroups of $S$. 
\end{proof}

Theorem \ref{Alp.fusion}(a) will be used repeatedly throughout this paper.  
The following lemma is a first easy application of the theorem, and 
provides a very useful criterion for a subgroup to be quasicentric or not. 

\begin{Lem} \label{L:qcentric}
Let $\calf$ be a saturated fusion system over a $p$-group $S$. Then the 
following hold for any $P\le{}S$.
\begin{enumerate} 
\item Assume that $P\le{}Q\le{}P{\cdot}C_S(P)$ and $\Id\ne\alpha\in\autf(Q)$ 
are such that $\alpha|_P=\Id_P$ and $\alpha$ has order prime to $p$.  Then 
$P$ is not $\calf$-quasicentric.
\item Assume that $P$ is fully centralized in $\calf$, and is not 
$\calf$-quasicentric.  Then there are $P\le{}Q\le{}P{\cdot}C_S(P)$ and 
$\Id\ne\alpha\in\autf(Q)$ such that $Q$ is $\calf$-centric, 
$\alpha|_P=\Id_P$, and $\alpha$ has order prime to $p$.  
\end{enumerate}
\end{Lem}

\begin{proof} \textbf{(a) } Fix any $P'$ which is $\calf$-conjugate to $P$ 
and fully centralized in $\calf$.  By axiom (II), there is 
$\varphi\in\homf(Q,S)$ such that $\varphi(P)=P'$; set $Q'=\varphi(Q)$.  
Thus $\varphi\alpha\varphi^{-1}|_{C_{Q'}(P')}$ is an automorphism in 
$C_\calf(P')$ whose order is not a power of $p$, so $C_\calf(P')$ is not 
the fusion system of $C_S(P')$, and $P$ is not $\calf$-quasicentric.

\noindent\textbf{(b) } Assume that $P$ is fully centralized in $\calf$ and 
not $\calf$-quasicentric.  Then $C_\calf(P)$ strictly contains the fusion 
system of $C_S(P)$ (since $C_\calf(P')$ is isomorphic as a category to 
$C_\calf(P)$ for all $P'$ which is $\calf$-conjugate to $P$ and fully 
centralized in $\calf$).  Since $C_\calf(P)$ is saturated by 
\cite[Proposition A.6]{BLO2}, Theorem \ref{Alp.fusion}(a) implies there is 
a subgroup $Q\le{}C_S(P)$ which is $C_\calf(P)$-centric and fully 
normalized in $C_\calf(P)$, and such that 
$\Aut_{C_\calf(P)}(Q)\gneqq\Aut_{C_S(P)}(Q)$.  Since $Q$ is fully 
normalized, $\Aut_{C_S(P)}(Q)$ is a Sylow $p$-subgroup of 
$\Aut_{C_\calf(P)}(Q)$, and hence this last group is not a $p$-group.  
Also, by \cite[Proposition 2.5(a)]{BLO2}, $PQ$ is $\calf$-centric since 
$Q$ is $C_\calf(P)$-centric.
\end{proof}

Since orbit categories --- both of fusion systems and of groups --- will 
play a role in certain proofs in the last three sections, we define them 
here.  

\begin{Defi} \label{D:orbit}
\begin{enumerate}  
\item If $\calf$ is a fusion system over a $p$-group $S$, then 
$\orb^c(\calf)$ (the \emph{centric orbit category} of $\calf$) is the 
category whose objects are the $\calf$-centric subgroups of $S$, and whose
morphism sets are given by
	$$ \Mor_{\orb^c(\calf)}(P,Q) = \repf(P,Q) \defeq 
	Q{\backslash}\homf(P,Q). $$
Let $\calz_\calf\:\orb^c(\calf)\rTo\zploc\mod$ be the functor which sends 
$P$ to $Z(P)$ and $[\varphi]$ (the class of $\varphi\in\homf(P,Q)$) to 
$Z(Q)\Right3{\varphi^{-1}}Z(P)$. 

\item If $G$ is a finite group and $S\in\sylp{G}$, then $\orb_S^c(G)$ 
(the \emph{centric orbit category} of $G$) is the category whose objects 
are the subgroups of $S$ which are $p$-centric in $G$, and where 
	$$ \Mor_{\orb^c_S(G)}(P,Q) = Q{\backslash}N_G(P,Q)
	\cong \map_G(G/P,G/Q). $$
Let $\calz_G\:\orb^c_S(G)\rTo\zploc\mod$ be the functor which sends 
$P$ to $Z(P)$ and $[g]$ (the class of $g\in N_G(P,Q)$) to conjugation by 
$g^{-1}$.
\end{enumerate}
\end{Defi}

We now turn to linking systems associated to abstract fusion systems.

\begin{Defi}[{\cite[Def.~1.7]{BLO2}}]  \label{L-cat}
Let $\calf$ be a fusion system over the $p$-group $S$.  A \emph{centric 
linking system associated to $\calf$} is a category $\call$ whose objects 
are the $\calf$-centric subgroups of $S$, together with a functor 
$\pi\:\call\Right2{}\calfc$, and ``distinguished'' monomorphisms 
$P\Right1{\delta_P}\Aut_{\call}(P)$ for each $\calf$-centric subgroup 
$P\le{}S$, which satisfy the following conditions.
\begin{enumerate}
\renewcommand{\labelenumi}{\textup{(\Alph{enumi})}}%
\item  $\pi$ is the identity on objects.  For each pair of objects
$P,Q\in\call$, $Z(P)$ acts freely on $\Mor_{\call}(P,Q)$ by composition
(upon identifying $Z(P)$ with $\delta_P(Z(P))\le\Aut_{\call}(P)$), and
$\pi$ induces a bijection
	$$ \Mor_{\call}(P,Q)/Z(P) \Right5{\cong} \homf(P,Q). $$

\item  For each $\calf$-centric subgroup $P\le{}S$ and each $x\in{}P$,
$\pi(\delta_P(x))=c_x\in\Aut_{\calf}(P)$.

\item  For each $f\in\Mor_{\call}(P,Q)$ and each $x\in{}P$, the following
square commutes in $\call$:
    \begin{diagram}
    P & \rTo^{f} & Q \\ \dTo>{\delta_P(x)} &&
    \dTo>{\delta_Q(\pi(f)(x))} \\
    P & \rTo^{f} & Q.
    \end{diagram}
\end{enumerate}
\end{Defi}

A \emph{$p$-local finite group} is defined to be a triple $\SFL$, where $S$
is a finite $p$-group, $\calf$ is a saturated fusion system over $S$, and
$\call$ is a centric linking system associated to $\calf$.  The
\emph{classifying space} of the triple $\SFL$ is the $p$-completed nerve
$|\call|\pcom$.

For any finite group $G$ with Sylow $p$-subgroup $S$, a category 
$\call^c_S(G)$ was defined in \cite{BLO1}, whose objects are the 
$p$-centric subgroups of $G$, and whose morphism sets are defined by
	$$ \Mor_{\call^c_S(G)}(P,Q) = N_G(P,Q)/O^p(C_G(P)). $$
Since $C_G(P)=Z(P)\times{}O^p(C_G(P))$ when $P$ is $p$-centric in $G$, 
$\call^c_S(G)$ is easily seen to satisfy conditions (A), (B), and (C) above, 
and hence is a centric linking system associated to $\calf_S(G)$.  Thus 
$(S,\calf_S(G),\call_S^c(G))$ is a $p$-local finite group, with 
classifying space $|\call^c_S(G)|\pcom\simeq{}BG\pcom$ (see 
\cite[Proposition~1.1]{BLO1}).

It will be of crucial importance in this paper that any centric linking 
system associated to $\calf$ can be extended to a quasicentric linking 
system; a linking system with similar properties, whose objects are the 
$\calf$-quasicentric subgroups of $S$.  We first make more precise what 
this means.

\begin{Defi}  \label{L^q}
Let $\calf$ be any saturated fusion system over a $p$-group $S$.  A 
\emph{quasicentric linking system} associated to $\calf$ consists of a 
category $\callq$ whose objects are the $\calf$-quasicentric subgroups of 
$S$, together with a functor $\pi\:\callq\Right2{}\calfq$, and 
distinguished monomorphisms 
	$$ P{\cdot}C_S(P)\Right3{\delta_P}\Aut_{\callq}(P), $$ 
which satisfy the following conditions. 
\begin{enumerate}\renewcommand{\labelenumi}{\textup{(\Alph{enumi})$_q$}}%
\item  $\pi$ is the identity on objects and surjective on morphisms. For 
each pair of objects $P,Q\in\callq$ such that $P$ is fully centralized, 
$C_S(P)$ acts freely on $\Mor_{\callq}(P,Q)$ by composition (upon 
identifying $C_S(P)$ with $\delta_P(C_S(P))\le\Aut_{\callq}(P)$), and 
$\pi$ induces a bijection
    $$ \Mor_{\callq}(P,Q)/C_S(P) \Right5{\cong} \homf(P,Q). $$

\item  For each $\calf$-quasicentric subgroup $P\le{}S$ and each $g\in{}P$,
$\pi$ sends $\delta_P(g)\in\Aut_{\callq}(P)$ to
$c_g\in\Aut_{\calf}(P)$.

\item  For each $f\in\Mor_{\callq}(P,Q)$ and each $x\in{}P$,
$f\circ\delta_P(x)=\delta_Q(\pi(f)(x))\circ{}f$ in $\Mor_{\callq}(P,Q)$.

\item  For each $\calf$-quasicentric subgroup $P\le{}S$, there is a 
morphism $\iota_P\in\Mor_{\callq}(P,S)$ such that 
$\pi(\iota_P)=\incl_P^S\in\Hom(P,S)$, and such that for each 
$g\in{}P{\cdot}C_S(P)$, $\delta_S(g)\circ\iota_P=\iota_P\circ\delta_P(g)$ 
in $\Mor_{\callq}(P,S)$.
\end{enumerate}
\end{Defi}

If $P$ and $P'$ are $\calf$-conjugate and $\calf$-quasicentric, then for any 
$Q\le{}S$, $\Mor_{\callq}(P,Q)\cong\Mor_{\callq}(P',Q)$ and 
$\homf(P,Q)\cong\homf(P',Q)$, while the centralizers $C_S(P)$ and $C_S(P')$ 
need not have the same order.  This is why condition (A)$_q$ makes sense 
only if we assume that $P$ is fully centralized; i.e., that $C_S(P)$ is as 
large as possible.  When $P$ is $\calf$-centric, then this condition is 
irrelevant, since every subgroup $P'$ which is $\calf$-conjugate to $P$ is 
fully centralized ($C_S(P')=Z(P')\cong{}Z(P)$).  

Note that (D)$_q$ is a special case of (C)$_q$ when $P$ is 
$\calf$-centric; this is why the axiom is not needed for centric linking 
systems.  We also note the following relation between these axioms:

\begin{Lem} \label{C=>B}
In the situation of Definition \ref{L^q}, axiom (C)$_q$ implies axiom 
(B)$_q$.
\end{Lem}

\begin{proof} Fix an $\calf$-quasicentric subgroup $P\le{}S$, and an 
element $g\in{}P$.  We apply (C)$_q$ with $f=\delta_P(g)$.  For each 
$x\in{}P$, if we set $y=\pi(\delta_P(g))(x)$, then 
$\delta_P(g)\circ\delta_P(x)=\delta_P(y)\circ\delta_P(g)$.  Since 
$\delta_P$ is an injective homomorphism, this implies that $gx=yg$, and 
thus that $y=c_g(x)$.  So $\pi(\delta_P(g))=c_g$.
\end{proof}

When $\SFL$ is a $p$-local finite group, and $\callq$ is a quasicentric 
linking system associated to $\calf$, then we say that $\callq$ 
\emph{extends} $\call$ if the full subcategory of $\callq$ with objects the 
$\calf$-centric subgroups of $S$ is isomorphic to $\call$ via a functor 
which commutes with the projection functors to $\calf$ and with the 
distinguished monomorphisms.  In \cite[Propositions 3.4 \& 3.12]{bcglo1}, 
we constructed an explicit quasicentric linking system $\callq$ associated 
to $\calf$ and extending $\call$, and showed that it is unique up to an 
isomorphism of categories which preserves all of these structures.  So 
from now on, we will simply refer to $\callq$ as \emph{the} quasicentric 
linking system associated to $\SFL$.  

Condition (D)$_q$ above helps to motivate the following definition of 
\emph{inclusion morphisms} in a quasicentric linking system.

\begin{Defi} \label{D:incl}
Fix a $p$-local finite group $\SFL$, with associated quasicentric linking 
system $\callq$.
\begin{enumerate}\renewcommand{\labelenumi}{\textup{(\alph{enumi})}}
\item A morphism $\iota_P\in\Mor_{\callq}(P,S)$ is 
an \emph{inclusion morphim} if it satisfies the hypotheses of axiom 
(D)$_q$:  if $\pi(\iota_P)=\incl_P^S$, and if 
$\delta_S(g)\circ\iota_P=\iota_P\circ\delta_P(g)$
in $\Mor_{\callq}(P,S)$ for all $g\in{}P{\cdot}C_S(P)$.  If $P=S$, then we 
also require that $\iota_S=\Id_S$.
\item A \emph{compatible set of inclusions} for $\callq$ is a choice of 
morphisms $\{\iota_P^Q\}$ for all pairs of $\calf$-quasicentric subgroups 
$P\leq Q$, such that $\iota_P^Q\in\Mor_{\callq}(P,Q)$, such that 
$\iota_P^R=\iota_Q^R \circ\iota_P^Q$ for all $P\leq Q\leq R$, and such 
that $\iota_P^S$ is an inclusion morphism for each $P$.
\end{enumerate}
\end{Defi}

The following properties of quasicentric linking systems were also proven 
in \cite{bcglo1}.

\begin{Prop}  \label{L-props}
The following hold for any $p$-local finite group $\SFL$, with associated
quasicentric linking system $\callq$.
\begin{enumerate}  
\item The inclusion $\call\subseteq\callq$ induces a homotopy equivalence 
$|\call|\simeq|\callq|$ between geometric realizations.  More generally, 
for any full subcategory $\call'\subseteq\callq$ which contains as objects 
all subgroups of $S$ which are $\calf$-centric and $\calf$-radical, the 
inclusion $\call'\subseteq\callq$ induces a homotopy equivalence 
$|\call'|\simeq|\callq|$.  

\item Let $\varphi\in\Mor_{\callq}(P,R)$ and $\psi\in\Mor_{\callq}(Q,R)$ 
be any pair of morphisms in $\callq$ with the same target group such that 
$\Im(\pi(\varphi))\le\Im(\pi(\psi))$.  Then there is a unique morphism 
$\chi\in\Mor_{\callq}(P,Q)$ such that $\varphi=\psi\circ\chi$.
\end{enumerate}
\end{Prop}

\begin{proof}  The homotopy equivalences 
$|\callq|\simeq|\call|\simeq|\call'|$ are shown in \cite[Theorem 
3.5]{bcglo1}.  Point (b) is shown in \cite[Lemma 3.6]{bcglo1}.  
\end{proof}

Point (b) above will be frequently used throughout the paper.  In 
particular, it makes it possible to embed the linking system of $S$ (or an 
appropriate full subcategory) in $\callq$, depending on the choice of an 
inclusion morphisms $\iota_P$ as defined above, for each object $P$.  
Such inclusion morphisms always exist by axiom (D)$_q$.  In the 
following proposition, $\call_S(S)|_{\Ob(\callq)}$ denotes the full 
subcategory of $\call_S(S)$ whose objects are the $\calf$-quasicentric 
subgroups of $S$.  

In general, for a functor $F\:\calc\rTo\calc'$, and objects 
$c,d\in\Ob(\calc)$, we let $F_{c,d}$ denote the map from 
$\Mor_{\calc}(c,d)$ to $\Mor_{\calc'}(F(c),F(d))$ induced by $F$.

\begin{Prop}  \label{deltaPQ}
Fix a $p$-local finite group $\SFL$. let $\callq$ be its associated 
quasicentric linking system, and let $\pi\:\callq\Right1{}\calfq$ be the 
projection.  Then any choice of inclusion morphisms 
$\iota_P=\iota_P^S\in\Mor_{\callq}(P,S)$, for all $\calf$-quasicentric 
subgroups $P\le{}S$, extends to a unique inclusion of categories
        $$ \delta\: \call_S(S)|_{\Ob(\callq)} \Right5{} \callq $$
such that $\delta_{P,S}(1)=\iota_P$ for all $P$; and such that
\begin{enumerate}
\item $\delta_{P,P}(g)=\delta_P(g)$ for all $g\in P\cdot C_S(P)$, and
\item $\pi(\delta_{P,Q}(g))=c_g\in\Hom(P,Q)$ for all $g\in N_S(P,Q)$.
\end{enumerate}
In addition, the following hold.
\begin{enumerate}\setcounter{enumi}{2}
\item If we set $\iota_P^Q=\delta_{P,Q}(1)$ for all $P\le Q$, then
$\{\iota_P^Q\}$ is a compatible set of inclusions for $\callq$.  

\item For any $P\nsg{}Q\le{}S$, where $P$ and $Q$ are both 
$\calf$-quasicentric and $P$ is fully centralized in $\calf$, and any 
morphism $\psi\in\Aut_{\callq}(P)$ which normalizes $\delta_{P,P}(Q)$, 
there is a unique $\widehat{\psi}\in\Aut_{\callq}(Q)$ such that 
$\widehat{\psi}\circ\iota_P^Q=\iota_P^Q\circ\psi$.  Furthermore, for any 
$g\in{}Q$, $\psi\delta_{P,P}(g)\psi^{-1}= 
\delta_{P,P}(\pi(\widehat{\psi})(g))$. 

\item Every morphism $\varphi\in\Mor_{\callq}(P,Q)$ in $\callq$ is a 
composite $\varphi=\iota_{P'}^Q\circ\varphi'$ for a unique morphism 
$\varphi'\in\Iso_{\callq}(P,P')$, where $P'=\Im(\pi(\varphi))$.
\end{enumerate}
\end{Prop}

\begin{proof}  For each $P$ and $Q$, and each $g\in{}N_S(P,Q)$, there is
by Proposition \ref{L-props}(b) a unique morphism $\delta_{P,Q}(g)$ such
that
        $$ \delta_S(g)\circ\iota_P = \iota_Q\circ\delta_{P,Q}(g). $$
This defines $\delta$ on morphism sets, and also allows us to define 
$\iota_P^Q=\delta_{P,Q}(1)$.  Then by the axioms in Definition \ref{L^q}, 
$\{\iota_P^Q\}$ is a compatible set of inclusions for $\callq$, and 
$\delta$ is a functor which satisfies (a), (b), and (c).  Point (e) is a 
special case of Proposition \ref{L-props}(b) (where $P'=\Im(\pi(\varphi))$).

If $\delta_{P,Q}(g)=\delta_{P,Q}(g')$ for $g,g'\in{}N_G(P,Q)$, then 
$\delta_S(g)\circ\iota_P=\delta_S(g')\circ\iota_P$, and hence $g=g'$ by 
\cite[Lemma 3.9]{bcglo1}.  Thus each $\delta_{P,Q}$ is injective.

It remains to prove (d).  Set $\varphi=\pi(\psi)\in\autf(P)$ for short.  
Since $\psi$ normalizes $\delta_{P,P}(Q)$, for all $g\in{}Q$ there is 
$h\in{}Q$ such that $\psi\delta_{P,P}(g)\psi^{-1}=\delta_{P,P}(h)$, and 
this implies the relation $\varphi{}c_g\varphi^{-1}=c_h$ in $\autf(P)$.  
Thus $Q$ is contained in $N_\varphi$.  So by 
axiom (II) (and since $P$ is fully centralized), 
$\varphi$ extends to $\widebar{\varphi}\in\isof(Q,Q')$ for some 
$P\nsg{}Q'\le{}S$.  Let $\widehat{\psi}_0\in\Iso_{\callq}(Q,Q')$ be any 
lifting of $\widebar{\varphi}$ to $\callq$.  

By axiom (A)$_q$ (and since $\psi$ is an isomorphism), there is 
$x\in{}C_S(P)$ such that $\iota_P^{Q'}\circ\delta_P(x)\circ\psi= 
\widehat{\psi}_0\circ\iota_P^Q$.  By (C)$_q$, 
$\widehat{\psi}_0\delta_{Q}(Q)\widehat{\psi}_0^{-1}=\delta_{Q'}(Q')$, and 
hence after restriction, $\delta_P(x)\circ\psi$ conjugates 
$\delta_{P,P}(Q)$ to $\delta_{P,P}(Q')$.  Since $\psi$ normalizes 
$\delta_{P,P}(Q)$, this shows that $\delta_P(x)$ conjugates 
$\delta_{P,P}(Q)$ to $\delta_{P,P}(Q')$, and hence (since $\delta_{P,P}$ is 
injective) that $xQx^{-1}=Q'$.  We thus have the following commutative 
diagram:
	\begin{diagram}[w=40pt]
	P & \rTo^{\psi} & P & \rTo^{\delta_P(x)} & P & 
	\rTo^{\delta_P(x)^{-1}} & P \\
	\dTo<{\iota_P^Q} &&&& \dTo<{\iota_P^{Q'}} && \dTo<{\iota_P^Q} \\
	Q && \rTo^{\widehat{\psi}_0} && Q' & \rTo^{\delta_{Q,Q'}(x)^{-1}} 
	& Q \rlap{\,.}
	\end{diagram}
So if we set $\widehat{\psi}=\delta_{Q,Q'}(x)^{-1}\circ\widehat{\psi}_0$, 
then $\widehat{\psi}\in\Aut_{\callq}(Q)$ and 
$\iota_P^Q\circ\psi=\widehat{\psi}\circ\iota_P^Q$.  

The uniqueness of $\widehat{\psi}$ follows from \cite[Lemma 3.9]{bcglo1}.  
Finally, for any $g\in{}Q$, $\widehat{\psi}\delta_Q(g)\widehat{\psi}^{-1}= 
\delta_Q(\pi(\widehat{\psi})(g))$ by (C)$_q$, and hence 
$\psi\delta_{P,P}(g)\psi^{-1}= \delta_{P,P}(\pi(\widehat{\psi})(g))$ since 
morphisms have unique restrictions (Proposition \ref{L-props}(b) again).
\end{proof}


Once we have fixed a compatible set of inclusions $\{\iota_P^Q\}$ in a 
linking system $\callq$, then for any $\varphi\in\Mor_{\callq}(P,Q)$, and 
any $P'\le{}P$ and $Q'\le{}Q$ such that $\pi(\varphi)(P')\le{}Q'$, there 
is a unique morphism $\widebar{\varphi}\in\Mor_{\callq}(P',Q')$ such that 
$\iota_{Q'}^Q\circ\widebar{\varphi}=\varphi\circ\iota_{P'}^P$.  We think 
of $\widebar{\varphi}$ as the \emph{restriction} of $\varphi$.

Note, however, that all of this depends on the choice of morphisms 
$\iota_P\in\Mor_{\callq}(P,S)$ which satisfy the hypotheses of axiom 
(D)$_q$, and that not just any  lifting of the inclusion 
$\incl\in\Hom_{\calf}(P,S)$ can be chosen. To see why, assume for 
simplicity that $P$ is also fully centralized. From the axioms in 
Definition~\ref{L^q} and Proposition~\ref{L-props}(b), we see that if 
$\iota_P, \iota'_P\in\Mor_{\callq}(P,S)$ are two liftings of 
$\incl\in\Hom_{\calf}(P,S)$, then $\iota'_P=\iota_P\circ\delta_P(g)$ for 
some unique $g\in C_S(P)$.  But if $\iota_P$ satisfies the conditions of 
(D)$_q$, then $\iota'_P$ also satisfies those conditions only if 
$g\in{}Z(C_S(P))$.

One situation where the choice of inclusion morphisms is useful is when 
describing the fundamental group of $|\call|$ or of its $p$-completion.  
For any group $\Gamma$, we let $\calb(\Gamma)$ denote the category with 
one object, and with morphism monoid the group $\Gamma$.  Recall that 
$|\call|\simeq|\callq|$ (Proposition \ref{L-props}(a)), so we can work 
with either of these categories; we will mostly state the results for 
$|\callq|$.  Let the vertex $S$ be the basepoint of $|\callq|$.  For each 
morphism $\varphi\in\Mor_{\callq}(P,Q)$, let 
$J(\varphi)\in\pi_1(|\callq|)$ denote the homotopy class of the loop 
$\iota_Q{\cdot}\varphi{\cdot}\iota_P{}^{-1}$ in $|\callq|$ (where paths 
are composed from right to left).  This defines a functor 
	$$ J\: \callq \Right5{} \calb(\pi_1(|\callq|)), $$
where all objects are sent to the unique object of 
$\calb(\pi_1(|\callq|))$, and where all inclusion morphisms are sent to the 
identity.  Let $j\:S\rTo\pi_1(|\callq|)$ denote 
the composite of $J$ with the distinguished monomorphism 
$\delta_S\:S\rTo\Aut_\call(S)$.  

The next proposition describes how $J$ is universal among functors of this 
type, and also includes some other technical results for later use about 
the structure of $\pi_1(|\call|)$.

\begin{Prop}\label{conj-pi1L} 
Let $\SFL$ be a p-local finite group, and let $\callq$ be the associated 
quasicentric linking system.  Assume a compatible set of inclusions 
$\{\iota_P^Q\}$ has been chosen for $\callq$.  Then the following hold.
\begin{enumerate}
\item For any group $\Gamma$, and any functor 
$\lambda\:\callq\Right3{}\calb(\Gamma)$ 
which sends inclusions to the identity, there is a unique homomorphism 
$\bar\lambda\:\pi_1(|\callq|)\Right3{}\Gamma$ such that 
$\lambda=\calb(\bar\lambda)\circ{}J$.  
\item For $g\in P\leq S$ with $P$ $\calf$-quasicentric, 
$J(\delta_P(g))=J(\delta_S(g))$.  In particular, $J(\delta_P(g))=1$ in 
$\pi_1(|\callq|)$ if and only if $\delta_P(g)$ is nulhomotopic as a loop 
based at the vertex $P$ of  $\nv{\callq}$.
\item If $\alpha\in\Mor_{\callq}(P,Q)$, and $\pi(\alpha)(x)= y$, then 
$j(y)= J(\alpha)j(x)J(\alpha)^{-1}$ in $\pi_1(|\callq|)$. 
\item If $x$ and $y$ are $\calf$-conjugate elements of $S$, then 
$j(x)$ and $j(y)$ are conjugate in $\pi_1(|\callq|)$. 
\end{enumerate}
\end{Prop}

\begin{proof}  Clearly, \emph{any} functor $\lambda\:\callq\rTo\calb(\Gamma)$ 
induces a homomorphism $\bar\lambda=\pi_1(|\lambda|)$ between the 
fundamental groups of their geometric realizations.  If $\lambda$ sends 
inclusion morphisms to the identity, then 
$\lambda=\calb(\bar\lambda)\circ{}J$ by definition of $J$.  

The other points follow easily, using condition (C)$_q$ for quasicentric 
linking systems.  Point (d) is shown by first reducing to a map between 
centric subgroups of $S$ which sends $x$ to $y$. 
\end{proof}

We finish this introductory section with two unrelated results which will 
be needed later in the paper.  The first is a standard, group theoretic 
lemma.

\begin{Lem} \label{gorenstein}
Let $Q \nsg P$ be $p$-groups.  If $\alpha$ is a $p'$-automorphism
of $P$ which acts as the identity on $Q$ and on $P/Q$, then
$\alpha=\Id_P$. Equivalently, the group of all automorphisms of
$P$ which restrict to the identity on $Q$ and on $P/Q$ is a
$p$-group.
\end{Lem}

\begin{proof}  See \cite[Corollary 5.3.3]{Gorenstein}.
\end{proof}

\newcommand{\repffc}{\repf^{\textup{fc}}}
\newcommand{\repffn}{\repf^{\textup{fn}}}
\newcommand{\calpfn}{\calp_{\textup{fn}}}

The following proposition will only be used in Section 4, but we include 
it here because it seems to be of wider interest.  Note, for any fusion 
system $\calf$ over $S$, any subgroup $P\le{}S$ fully normalized in 
$\calf$, and any $P'$ which is $S$-conjugate to $P$, that $P'$ is also 
fully normalized in $\calf$ since $N_S(P')$ is $S$-conjugate to $N_S(P)$. 

\begin{Prop} \label{|Rep(P,S)|}
Let $\calf$ be a saturated fusion system over a $p$-group $S$.  Then for 
any subgroup $P\le{}S$, the set of $S$-conjugacy classes of subgroups 
$\calf$-conjugate to $P$ and fully normalized in $\calf$ has order prime 
to $p$.
\end{Prop}

\begin{proof}  By \cite[Proposition 5.5]{BLO2}, there is an 
$(S,S)$-biset $\Omega$ which, when regarded as a set with 
$(S\times{}S)$-action, satisfies the following three conditions.
\begin{enumerate}  
\item The isotropy subgroup of each point in $\Omega$ is of the form
	$$ P_\varphi \defeq \{(x,\varphi(x))\,|\,x\in{}P\} $$
for some $P\le{}S$ and some $\varphi\in\homf(P,S)$.
\item For each $P\le{}S$ and each $\varphi\in\homf(P,S)$, the two 
structures of $(S\times{}P)$-set on $\Omega$ obtained by restriction and by 
$\Id\times\varphi$ are isomorphic.  
\item $|\Omega|/|S|\equiv1$ (mod $p$).
\end{enumerate}
Note that by (a), the actions of $S\times1$ and $1\times{}S$ on $\Omega$ are 
both free.

Now fix a subgroup $P\le{}S$.  Set $S_2=1\times{}S$ for short, and let 
$\Omega_0\subseteq\Omega$ be the subset such that 
$\Omega_0/S_2=(\Omega/S_2)^P$.  In other words, $\Omega_0$ is the set of 
all $x\in\Omega$ such that for each $g\in P$, there is some $h\in S$ 
satisfying $(g,h){\cdot}x=x$. Since the action of $S_2$ on $\Omega$ is 
free, this element $h\in S$ is uniquely determined for each $x\in 
\Omega_0$ and $g\in P$. Let $\theta(x)\:P\rTo S$ denote the function such 
that for each $g\in P$, $(g,\theta(x)(g)){\cdot}x = x$. The isotropy 
subgroup at $x$ of the $(P\times{}S)$-action is thus the subgroup 
$P_{\theta(x)} = \{(g,\theta(x)(g))\;|\; g\in P\}$; and by (a), 
$\theta(x)\in\homf(P,S)$. This defines a map 
	$$ \theta \: \Omega_0\rTo \homf(P,S). $$

By definition, for each $\varphi\in\homf(P,S)$, $\theta^{-1}(\varphi)$ is 
the set of elements of $\Omega$ fixed by $P_\varphi$.  
By condition (b) above, the action of $P\times{}P$ on $\Omega$ induced by 
the homomorphism $1\times\varphi\in\Hom(P\times{}P,S\times{}S)$ is 
isomorphic to the action defined by restriction, and thus 
$|\theta^{-1}(\varphi)|=|\theta^{-1}(\incl)|$.  This shows that the point 
inverses of $\theta$ all have the same fixed order $k$.  

Now let $\repf(P,S)=\homf(P,S)/\Inn(S)$:  the set of $S$-conjugacy classes 
of morphisms from $P$ to $S$.  Let $\calp$ be the set of $S$-conjugacy 
classes of subgroups $\calf$-conjugate to $P$, and let 
$\calpfn\subseteq\calp$ be the subset of classes of subgroups fully 
normalized in $\calf$.  If $x\in\Omega_0$ and $\theta(x)=\varphi$, then for 
all $s\in{}S$ and $g\in{}P$,
	$$ (1,s){\cdot}x = (g,s\varphi(g)){\cdot}x = 
	(g,c_s\circ\varphi(g)){\cdot}(1,s){\cdot}x, $$
and this shows that $\theta((1,s){\cdot}x)=c_s\circ\theta(x)$.  Thus
$\theta$ induces a map
	$$ \bar{\theta} \: (\Omega/S_2)^P = \Omega_0/S_2 
	\Right5{\theta/S_2} \repf(P,S) \Right5{\Im(-)} \calp, $$
where $\repf(P,S)=\homf(P,S)/\Inn(S)$.  
Furthermore, $|\Omega/S_2|\equiv1$ (mod $p$) by (c), and thus 
$|(\Omega/S_2)^P|\equiv1$ (mod $p$).

For each $P'$ which is $\calf$-conjugate to $P$, there are $|S|/|N_S(P')|$ 
distinct subgroups in the $S$-conjugacy class $[P']$.  Hence there are
	$$ |\autf(P)|{\cdot}|S|/|N_S(P')| $$
elements of $\homf(P,S)$ whose image lies in $[P']$.  Since each of these 
is the image of $k$ elements in $\Omega_0$, this shows that
	$$ |\bar{\theta}^{-1}([P'])| = k{\cdot}|\autf(P)|/|N_S(P')|. $$
Thus $|N_S(P')|\big|k{\cdot}|\autf(P)|$ for all $P'$ $\calf$-conjugate to 
$P$, and so $\bar{\theta}^{-1}([P'])$ has order a multiple of $p$ if $P'$ 
is not fully normalized in $\calf$ (if $[P']\in\calp{\sminus}\calpfn$).  
Hence
	$$ |\bar{\theta}^{-1}(\calpfn)| \equiv |(\Omega/S_2)^P|
	\equiv 1 \pmod{p} \,. $$
So if we set $m=|N_S(P')|$ for $[P']\in\calpfn$ (i.e., the maximal value of 
$|N_S(P')|$ for $P'$ $\calf$-conjugate to $P$), then
	$$ |\bar{\theta}^{-1}(\calpfn)| = 
	|\calpfn|{\cdot}\Bigl(\frac{k{\cdot}|\autf(P)|}{m}\Bigr), $$
and thus $|\calpfn|$ is prime to $p$.
\end{proof}

Using a similar argument, one can also show that the set 
$\repf^{\textup{fc}}(P,S)$ of elements of $\repf(P,S)$ whose image is fully 
centralized also has order prime to $p$.


\newsect{The fundamental group of $|\call|\pcom$} 

The purpose of this section is to give a simple description of the 
fundamental group of $|\call|\pcom$, for any $p$-local finite group 
$\SFL$, purely in terms of the fusion system $\calf$. The result is
analogous to the (hyper-) focal subgroup theorem for finite groups, as
we explain below.

In Section 1, we defined a functor $J\:\callq\rTo\calb(\pi_1(|\call|))$, 
for any $p$-local finite group $\SFL$, and a homomorphism 
$j=J\circ\delta_S$ from $S$ to $\pi_1(|\call|)$.  Let $\tau\: S\to 
\pi_1(|\call|\pcom)$ be the composite of $j$ with the natural homomorphism 
from $\pi_1(|\call|)$ to $\pi_1(|\call|\pcom)$. 

In \cite[Proposition 1.12]{BLO2}, we proved that 
$\tau\:S\rTo\pi_1(\nv{\call}\pcom)$ is a surjection. In this section, we 
will show that $\Ker(\tau)$ is the \emph{hyperfocal subgroup} of $\calf$, 
defined by Puig \cite{Puig-ppt} (see also \cite{puig00}). 


\begin{Defi}  \label{O^p(F)}
For any saturated fusion system $\calf$ over a $p$-group $S$,
the hyperfocal subgroup of $\calf$ is the normal subgroup of $S$ defined by 
    $$ \hfocal{\calf}{S} = \bigl\langle g^{-1}\alpha(g)
    \,\big|\, g\in{}P\le{}S,\ \alpha\in O^p(\autf(P))
    \bigr\rangle\,. $$
\end{Defi}

We will prove, for any $p$-local finite group $\SFL$, that 
$\pi_1(|\call|\pcom)\cong{}S/\hfocal{\calf}S$.  This is motivated by 
Puig's hyperfocal theorem, and we will also need that theorem in order to 
prove it.  Before stating Puig's theorem, we first recall the standard 
focal subgroup theorem.  If $G$ is a finite group and $S\in\sylp{G}$, then 
this theorem says that $S\cap[G,G]$ (the focal subgroup) is the subgroup 
generated by all elements of the form $x^{-1}y$ for $x,y\in{}S$ which are 
$G$-conjugate (cf. \cite[Theorem 7.3.4]{Gorenstein} or 
\cite[5.2.8]{Suz2}).  

The quotient group $S/(S\cap[G,G])$ is isomorphic to the $p$-power torsion 
subgroup of $G/[G,G]$, and can thus be identified as a quotient group of 
the maximal $p$-group quotient $G/O^p(G)$.  Since $G/O^p(G)$ is a 
$p$-group, $G=S{\cdot}O^p(G)$, and hence 
$G/O^p(G)\cong{}S/(S\cap{}O^p(G))$.  Hence $S\cap{}O^p(G)\le{}S\cap[G,G]$. 
This subgroup $S\cap{}O^p(G)$ is what Puig calls the hyperfocal subgroup, 
and is described by the hyperfocal subgroup theorem in terms of $S$ and 
fusion.

For $S\in\sylp{G}$ as above, let $\hfocal{G}S$ be the normal subgroup of 
$S$ defined by 
	\begin{align*} 
	\hfocal{G}{S} =\hfocal{\calf_S(G)}S &= \bigl\langle g^{-1}\alpha(g) 
	\,\big|\, g\in{}P\le{}S,\ \alpha\in O^p(\Aut_G(P)) \bigr\rangle \\
	&= \bigl\langle [g,x] \,\big|\, g\in{}P\le{}S,\ 
	x\in{}N_G(P) \textup{ of order prime to $p$} \bigr\rangle \,. 
	\end{align*}

\begin{Lem}[\cite{puig00}] \label{G:hyperfocal}
Fix a prime $p$, a finite group $G$, and a Sylow subgroup $S\in\sylp{G}$.  
Then $\hfocal{G}S=S\cap{}O^p(G)$.
\end{Lem}

\begin{proof}  This is stated in \cite[\S1.1]{puig00}, but the proof is 
only sketched there, and so we elaborate on it here.  Following 
standard notation, for any $P\le{}S$ and any $A\le\Aut(P)$, we write 
$[P,A]=\gen{x^{-1}\alpha(x)\,|\,x\in{}P,\ \alpha\in{}A}$.  Thus 
$\hfocal{G}S$ is generated by the subgroups $[P,O^p(\Aut_G(P))]$ for all 
$P\le{}S$. It is clear that $\hfocal{G}S\le{}S\cap{}O^p(G)$; the problem 
is to prove the opposite inclusion.  

Set $G_*=O^p(G)$ and $S_*=S\cap{}G_*$ for short.  Then $[G_*,G_*]$ has 
index prime to $p$ in $G_*$, so it contains $S_*$.  By the focal subgroup 
theorem (cf. \cite[Theorem 7.3.4]{Gorenstein}), applied to 
$S_*\in\sylp{G_*}$, $S_*$ is generated by all elements of the form 
$x^{-1}y$ for $x,y\in{}S_*$ which are $G_*$-conjugate.  Combined with 
Alperin's fusion theorem (in Alperin's original version \cite{Alp} or in 
the version of Theorem \ref{Alp.fusion}(a)), this implies that $S_*$ is 
generated by all subgroups $[P,N_{G_*}(P)]$ for $P\le{}S_*$ such that 
$N_{S_*}(P)\in\sylp{N_{G_*}(P)}$.  (This last condition is equivalent to 
$P$ being fully normalized in $\calf_{S_*}(G_*)$.)  Also, $N_{G_*}(P)$ is 
generated by $O^p(N_{G_*}(P))$ and the Sylow subgroup $N_{S_*}(P)$, so 
$[P,N_{G_*}(P)]$ is generated by 
$[P,O^p(N_{G_*}(P))]\le{}\hfocal{G}S$ and $[P,N_{S_*}(P)]\le[S_*,S_*]$. 
Thus $S_*=\gen{\hfocal{G}S,[S_*,S_*]}$.  Since $\hfocal{G}S$ is normal in 
$S$ (hence also normal in $S_*$), this shows that $S_*/\hfocal{G}S$ is 
equal to its commutator subgroup, which for a $p$-group is possible only 
if $S_*/\hfocal{G}S$ is trivial, and hence $S_*=\hfocal{G}S$.  
\end{proof}

By Proposition \ref{conj-pi1L}, the key to getting information about 
$\pi_1(|\call|)$ is to construct functors from $\call$ or $\callq$ to 
$\calb(\Gamma)$, for a group $\Gamma$, which send inclusions to the 
identity.  The next lemma is our main inductive tool for doing this.
Whenever $\alpha\in\Mor_{\callq}(P,P')$ and $\beta\in\Mor_{\callq}(Q,Q')$
are such that $P\le Q$, $P'\le Q'$, we write $\alpha=\beta|_P$ to mean
that $\alpha$ is the restriction of $\beta$ in the sense defined in 
Section 1; i.e., $\iota_{P'}^{Q'}\circ\alpha=\beta\circ\iota_P^Q$.

\begin{Lem} \label{ext.L}
Fix a $p$-local finite group $\SFL$, and let $\callq$ be its associated 
quasicentric linking system.  Assume a compatible set of inclusions 
$\{\iota_P^Q\}$ has been chosen for $\callq$.  Let $\calh_0$ be a set of 
$\calf$-quasicentric subgroups of $S$ which is closed under 
$\calf$-conjugacy and overgroups.  Let $\calp$ be an $\calf$-conjugacy 
class of $\calf$-quasicentric subgroups maximal among those not in 
$\calh_0$, set $\calh=\calh_0\cup\calp$, and let 
$\call^{\calh_0}\subseteq\call^{\calh}\subseteq\callq$ be the full 
subcategories with these objects.  Assume, for some group $\Gamma$, that 
	$$ \lambda_0 \: \call^{\calh_0} \Right4{} \calb(\Gamma) $$
is a functor which sends inclusions to the identity.  Fix some $P\in\calp$ 
which is fully normalized in $\calf$, and fix a homomorphism 
$\lambda_P\:\Aut_{\callq}(P)\rTo \Gamma$.  Assume that
\begin{enumerate} 
\item[\textup{($*$)}]  for all $P\lneqq{}Q\le{}N_S(P)$ such that 
$Q$ is fully normalized in $N_\calf(P)$, 
and for all $\alpha\in\Aut_{\callq}(P)$ and $\beta\in\Aut_{\callq}(Q)$ such 
that $\alpha=\beta|_P$, $\lambda_P(\alpha)=\lambda_0(\beta)$.
\end{enumerate}
Then there is a unique extension of $\lambda_0$ to a functor 
$\lambda\:\call^{\calh}\to\calb(\Gamma)$ which sends inclusions to 
the identity, and such that $\lambda(\alpha)=\lambda_P(\alpha)$ 
for all $\alpha\in\autf(P)$.  
\end{Lem}

\begin{proof}  The uniqueness of the extension is an immediate consequence 
of Theorem \ref{Alp.fusion}(a) (Alperin's fusion theorem).  

To prove the existence of the extension $\lambda$, we first show that 
($*$) implies the following (a prori stronger) statement:
\begin{enumerate} 
\item[\textup{($**$)}]  for all $Q,Q'\le{}S$ which strictly contain $P$, 
and for all $\beta\in\Mor_{\callq}(Q,Q')$ and $\alpha\in\Aut_{\callq}(P)$ 
such that $\alpha=\beta|_P$, $\lambda_P(\alpha)=\lambda_0(\beta)$.
\end{enumerate}
To see this, note first that it suffices to consider the case where $P$ is 
normal in $Q$ and $Q'$.  By assumption, $\pi(\beta)(P)=\pi(\alpha)(P)=P$, 
hence $\pi(\beta)(N_Q(P))\le{}N_{Q'}(P)$, and therefore $\beta$ restricts 
to a morphism $\widebar\beta\in\Mor_{\callq}(N_Q(P),N_{Q'}(P))$ by Proposition 
\ref{L-props}(b) (applied with $\varphi=\beta\circ\iota_{N_Q(P)}^Q$ and 
$\psi=\iota_{N_{Q'}(P)}^{Q'}$).  Since $N_Q(P),N_{Q'}(P)\gneqq{}P$, by the 
induction hypothesis, $\lambda_0(\beta)=\lambda_0(\widebar\beta)$, so we 
are reduced to proving that $\lambda_P(\alpha)=\lambda_0(\widebar\beta)$.  

Thus $\beta\in\Mor_{N_{\callq}(P)}(Q,Q')$.  We now apply Alperin's fusion 
theorem (Theorem \ref{Alp.fusion}(a)) to the morphism $\pi(\beta)$ in the 
fusion system $N_{\calf}(P)$ (which is saturated by \cite[Proposition 
A.6]{BLO2}). Thus $\pi(\beta)=\varphi_k\circ\cdots\circ\varphi_1$, where 
each $\varphi_i\in\Hom_{N_\calf(P)}(Q_{i-1},Q_i)$ is the restriction to 
$Q_{i-1}$ of an automorphism 
$\widehat{\varphi}_i\in\Aut_{N_{\calf}(P)}(R_i)$, where $R_i\geq 
Q_{i-1},Q_i$ is an $N_{\calf}(P)$-centric subgroup of $N_S(P)$ which is 
fully normalized in $N_{\calf}(P)$, and where $Q=Q_0$ and $Q'=Q_k$. Each 
$R_i$ contains $P$, and hence is $\calf$-centric by 
\cite[Lemma~6.2]{BLO2}.  For each $\widehat\varphi_i$, also regarded as an 
automorphism in $\calf$, we choose a lifting 
$\widehat\beta_i\in\Aut_{\call}(R_i)$, and let 
$\beta_i\in\Mor_{\callq}(Q_{i-1},Q_i)$ be its restriction.  By (A)$_q$, 
$\beta=\beta_k\circ\cdots\circ\beta_1\circ\delta_Q(g)$ for some 
$g\in{}C_S(Q)$, and hence 
	\begin{align*} 
	\lambda_0(\beta)&=\lambda_0(\beta_k)\cdots\lambda_0(\beta_1)
	\cdot \lambda_0(\delta_Q(g)) 
	= \lambda_0(\widehat\beta_k)\cdots\lambda_0(\widehat\beta_1)
	\cdot \lambda_0(\delta_{N_S(P)}(g)) \\
	&= \lambda_P(\widehat\beta_k|_P)\cdots\lambda_P(\widehat\beta_1|_P)
	\cdot \lambda_P(\Id_P) = \lambda_P(\beta|_P) = 
	\lambda_P(\alpha) \,,
	\end{align*}
where the third equality follows from ($*$).  This finishes the proof of 
($**$).

We can now extend $\lambda$ to be defined on all morphisms in 
$\call^{\calh}$ not in $\call^{\calh_0}$.  Fix such a morphism 
$\varphi\in\Mor_{\callq}(P_1,Q)$.  Set $P_2=\pi(\varphi)(P_1)\le{}Q$; then 
$P_1,P_2\in\calp$, and $\varphi=\iota_{P_2}^Q\circ\varphi'$ for some 
unique $\varphi'\in\Iso_{\callq}(P_1,P_2)$.  By Lemma \ref{N->N} (and then 
lifting to the linking category), there are isomorphisms 
$\widebar{\varphi}_i\in\Iso_{\callq}(N_S(P_i),N_i)$, for some 
$N_i\le{}N_S(P)$ containing $P$, which restrict to isomorphisms 
$\varphi_i\in\Iso_{\callq}(P_i,P)$.  Set 
$\psi=\varphi_2\circ\varphi'\circ\varphi_1^{-1}\in\Aut_{\callq}(P)$.  We have 
thus decomposed $\varphi'$ as the composite 
$\varphi_2^{-1}\circ\psi\circ\varphi_1$, and can now define
	$$ \lambda(\varphi) = \lambda(\varphi') = 
	\lambda_0(\widebar{\varphi}_2)^{-1} {\cdot}
	\lambda_P(\psi) {\cdot} \lambda_0(\widebar{\varphi}_1). $$

Now let $\varphi'=(\varphi'_2)^{-1}\circ\psi'\circ\varphi'_1$
be another such decomposition, where $\varphi'_i$ is the restriction of 
$\widebar{\varphi}'_i\in\Iso_{\callq}(N_S(P_i),N'_i)$.  We thus have a 
commutative diagram
	\begin{diagram}[w=30pt]
	P & \lTo^{\varphi_1} & P_1 & \rTo^{\varphi'_1} & P \\
	\dTo>{\psi} && \dTo>{\varphi'} && \dTo>{\psi'} \\
	P & \lTo^{\varphi_2} & P_2 & \rTo^{\varphi'_2} & P \rlap{\,,}
	\end{diagram}
where for each $i$, $\varphi_i$ and $\varphi'_i$ are restrictions of
isomorphisms $\widebar{\varphi}_i$ and $\widebar{\varphi}'_i$ defined on
$N_S(P_i)$.  To see that the two decompositions
give the same value of $\lambda(\varphi)$, it remains to show that
    $$ \lambda_P(\psi') {\cdot} 
    \lambda_0(\widebar{\varphi}'_1\circ(\widebar{\varphi}_1)^{-1}) =
    \lambda_0(\widebar{\varphi}'_2\circ(\widebar{\varphi}_2)^{-1}) {\cdot} 
    \lambda_P(\psi). $$
And this holds since 
$\lambda_0(\widebar{\varphi}'_i\circ(\widebar{\varphi}_i)^{-1})=
\lambda_P(\varphi'_i\circ(\varphi_i)^{-1})$ by ($**$).

We have now defined $\lambda$ on all morphisms in $\call^{\calh}$, and it 
sends inclusion morphisms to the identity by construction.  By 
construction, $\lambda$ sends composites to products, and thus the proof 
is complete.
\end{proof}

Lemma \ref{ext.L} provides the induction step when proving the following 
proposition, which is the main result needed to compute 
$\pi_1(|\call|\pcom)$. 
 
\begin{Prop} \label{ext-function}
Fix a $p$-local finite group $\SFL$, and let $\callq$ be its associated 
quasicentric linking system.  Assume a compatible set of inclusions 
$\{\iota_P^Q\}$ has been chosen for $\callq$.  Then there is a unique 
functor
	$$ \lambda\: \callq \Right5{} \calb(S/\hfocal{\calf}{S}) $$
which sends inclusions to the identity, and such that 
$\lambda(\delta_S(g))=g$ for all $g\in{}S$.  
\end{Prop}

\begin{proof} The functor $\lambda$ will be constructed inductively, using 
Lemma \ref{ext.L}.  Let $\calh_0\subseteq\Ob(\callq)$ be a subset 
(possibly empty) which is closed under $\calf$-conjugacy and overgroups.  
Let $\calp$ be an $\calf$-conjugacy class of $\calf$-quasicentric 
subgroups maximal among those not in $\calh_0$, set 
$\calh=\calh_0\cup\calp$, and let 
$\call^{\calh_0}\subseteq\call^{\calh}\subseteq\callq$ be the full 
subcategories with these objects.  Assume that
	$$ \lambda_0 \: \call^{\calh_0} \Right4{} \calb(S/\hfocal{\calf}S) $$
has already been constructed, such that $\lambda_0(\delta_S(g))=g$ for all 
$g\in{}S$ (if $S\in\calh_0$), and such that $\lambda_0$ sends inclusions 
to the identity.  

Fix $P\in\calp$ which is fully normalized in $\calf$, and let 
$\delta_{P,P}\:N_S(P)\Right2{}\Aut_{\callq}(P)$ be the homomorphism of 
Proposition \ref{deltaPQ}.  Then $\Im(\delta_{P,P})$ is a Sylow 
$p$-subgroup of $\Aut_{\callq}(P)$, since $\Aut_S(P)\in\sylp{\autf(P)}$ by 
axiom (I).  We identify $N_S(P)$ as a subgroup of $\Aut_{\callq}(P)$ to 
simplify notation.  Then 
    $$ \Aut_{\callq}(P)/O^p(\Aut_{\callq}(P)) \cong 
    N_S(P)\big/\bigl(N_S(P)\cap O^p(\Aut_{\callq}(P))\bigr) = N_S(P)/N_0, $$
where by Lemma \ref{G:hyperfocal}, $N_0$ is the subgroup generated by all 
commutators $[g,x]$ for $g\in{}Q\le{}N_S(P)$, and 
$x\in{}N_{\Aut_{\callq}(P)}(Q)$ of order prime to $p$. In this situation, 
conjugation by $x$ lies in $\autf(Q)$ by Proposition \ref{deltaPQ}(d), and 
thus $[g,x]=g{\cdot}c_x(g)^{-1}\in\hfocal{\calf}S$.  We conclude that 
$N_0\le\hfocal{\calf}{S}$, and hence that the inclusion of $N_S(P)$ into 
$S$ extends to a homomorphism 
	$$ \lambda_P\: \Aut_{\callq}(P) \Right5{} S/\hfocal{\calf}{S}. $$

We claim that condition ($*$) in Lemma \ref{ext.L} holds for $\lambda_0$ 
and $\lambda_P$.  To see this, fix $P\lneqq{}Q\le{}S$ such that $P\nsg{}Q$ 
and $Q$ is fully normalized in $N_\calf(P)$, and fix 
$\alpha\in\Aut_{\callq}(P)$ and $\beta\in\Aut_{\callq}(Q)$ such that 
$\alpha=\beta|_P$.  We must show that 
$\lambda_P(\alpha)=\lambda_0(\beta)$. Upon replacing $\alpha$ by 
$\alpha^k$ and $\beta$ by $\beta^k$ for some appropriate $k\equiv1$ (mod 
$p$), we can assume that both automorphisms have order a power of $p$.  
Since $Q$ is fully normalized, $\Aut_{N_S(P)}(Q)$ is a Sylow subgroup of 
$\Aut_{N_{\calf}(P)}(Q)$.  Hence (since any two Sylow $p$-subgroups of a 
finite group $G$ are conjugate by an element of $O^p(G)$), there is an 
automorphism $\widebar{\gamma}\in{}O^p(\Aut_{N_{\callq}(P)}(Q))$ such that 
$\widebar{\gamma}\beta\widebar{\gamma}^{-1}=\widebar{\delta}_Q(g)$ for 
some $g\in{}N_S(Q)\cap{}N_S(P)$.  In particular, 
$\lambda_0(\widebar{\gamma})=1$ since it is a composite of automorphisms 
of order prime to $p$.  Set $\gamma=\widebar{\gamma}|_P$; then 
$\gamma\in{}O^p(\Aut_{\callq}(P))$ and hence $\lambda_P(\gamma)=1$. Using 
axiom (C)$_q$ and Lemma \ref{L-props}, we see that 
$\gamma\alpha\gamma^{-1}=\widebar{\delta}_P(g)$; and thus (since 
$\lambda_P(\gamma)=1$ and $\lambda_0(\widebar{\gamma})=1$) that
    $$ \lambda_0(\beta)=\lambda_0(\widebar{\delta}_Q(g))=g
    =\lambda_P(\widebar{\delta}_P(g))=\lambda_P(\alpha). $$

Thus, by Lemma \ref{ext.L}, we can extend $\lambda_0$ to a functor defined 
on $\call^{\calh}$.  Upon continuing this procedure, we obtain a functor 
$\lambda$ defined on all of $\callq$.
\end{proof}

For any finite group $G$, $\pi_1(BG\pcom)\cong G/O^p(G)$.  (This is implicit 
in \cite[\S VII.3]{BK}, and shown explicitly in \cite[Proposition 
A.2]{BLO1}.)  Hence our main result in this section can be thought of as 
the hyperfocal theorem for $p$-local finite groups.

\begin{Thm}[Hyperfocal subgroup theorem for $p$-local finite groups] 
\label{pi1(hfg)}
Let $\SFL$ be a $p$-local finite group.  Then
    $$ \pi_1(\nv{\call}\pcom) \cong
    S/\hfocal{\calf}{S}. $$
More precisely, the natural map $\tau\:S\rTo\pi_1(|\call|\pcom)$ 
is surjective, and $\Ker(\tau)=\hfocal{\calf}S$. 
\end{Thm}

\begin{proof}  Let $\lambda\:\call\Right2{}\calb(S/\hfocal{\calf}S)$ be 
the functor of Proposition \ref{ext-function}, and let $|\lambda|$ be the 
induced map between geometric realizations.  Since 
$|\calb(S/\hfocal{\calf}S)|$ is the classifying space of a finite $p$-group 
and hence $p$-complete, $|\lambda|$ factors through the $p$-completion 
$|\call|\pcom$.  Consider the following commutative diagram:
	\begin{diagram}[w=45pt]
	S & \rTo^{j} & \pi_1(|\call|) \\
	& \rdTo_{\tau} & \dTo^{(-)\pcom} & \rdTo^{\pi_1(|\lambda|)} \\
	&& \pi_1(|\call|\pcom) & \rTo^{\bar\lambda=}_{\pi_1(|\lambda|\pcom)} 
	& S/\hfocal{\calf}S \rlap{\,.} 
	\end{diagram}
Here, $\tau$ is surjective by \cite[Proposition 1.12]{BLO2}.  Also, by 
construction, the composite $\bar\lambda\circ\tau=\pi_1(|\lambda|)\circ{}j$ 
is the natural projection.  Thus $\Ker(\tau)\le\hfocal{\calf}S$, and it 
remains to show the opposite inclusion. 

Fix $g\in{}P\le{}S$ and $\alpha\in\autf(P)$ such that $\alpha$ has order 
prime to $p$; we want to show that $g^{-1}\alpha(g)\in\Ker(\tau)$. Since 
$\ker(\tau)$ is closed under $\calf$-conjugacy (Proposition 
\ref{conj-pi1L}(d)), for any $\varphi\in\isof(P,P')$, $g^{-1}\alpha(g)$ is 
in $\Ker(\tau)$ if $\varphi(g^{-1}\alpha(g))$ is in $\Ker(\tau)$.  In 
particular, since $\varphi(g^{-1}\alpha(g))= 
\varphi(g)^{-1}{\cdot}(\varphi\alpha\varphi^{-1})(\varphi(g))$, we can 
assume that $P$ is fully centralized in $\calf$.  Then, upon extending 
$\alpha$ to an automorphism of $P{\cdot}C_S(P)$, which can be assumed also 
to have order prime to $p$ (replace it by an appropriate power if 
necessary), we can assume that $P$ is $\calf$-centric.  In this case, by 
Proposition~\ref{conj-pi1L}(c), $j(g)$ and $j(\alpha(g))$ are conjugate in 
$\pi_1(|\call|)$ by an element of order prime to $p$, and hence are equal 
in $\pi_1(|\call|\pcom)$ since this is a $p$-group.  This shows that 
$g^{-1}\alpha(g)\in\Ker(\tau)$, and finishes the proof of the theorem.
\end{proof}

The following result, which will be useful in Section 5, is of a similar 
nature, but much more elementary.

\begin{Prop} \label{pi1|L|->|F|}
Fix a $p$-local finite group $\SFL$, and let $\callq$ be its associated
quasicentric linking system.  Then the induced maps 
	$$ \pi_1(|\call|)\Right4{}\pi_1(|\calfc|)
	\qquad\textup{and}\qquad
	\pi_1(|\callq|)\Right4{}\pi_1(|\calfq|) $$
are surjective, and their kernels are generated by elements of $p$-power 
order.  
\end{Prop}

\begin{proof}  We prove this for $|\callq|$ and $|\calfq|$; a similar 
argument applies to $|\call|$ and $|\calfc|$.  Recall from Section 1 that 
we can regard $\pi_1(|\callq|)$ as the group generated by the loops 
$J(\alpha)$, for $\alpha\in\Mor(\callq)$, with relations given by 
composition of morphisms and making inclusion morphisms equal to 1. In a 
similar way, we regard $\pi_1(\nv{\calfq})$ as being generated by loops 
$J(\varphi)$ for $\varphi\in\Mor(\calfq)$. Since every morphism 
$\alpha\in\hom_{\calfq}(P,Q)$ has a lifting to a  morphism of $\callq$, 
the map $\pi_\#\:\pi_1(|\callq|)\rTo\pi_1(|\calfq|)$ induced by the 
projection functor $\pi\:\callq\rTo\calfq$ is an epimorphism.

Write $JC$ for the normal subgroup of $\pi_1(\nv{\callq})$ generated by 
the loops $j(g)=J(\delta_P(g))$, for all $\calf$-quasicentric subgroups 
$P\le{}S$ and all $g\in C_S(P)$.  In particular, $JC$ is 
generated by elements of $p$-power order.  Since $\pi_\#(j(g))=1$, $JC$ is 
contained in the kernel of $\pi_\#$, and we have a factorization 
$\widehat\pi_\#\: \pi_1(\nv{\callq})/JC \rTo\pi_1(\nv{\calfq})$.

Define an inverse  $s\: \pi_1(\nv{\calfq})\rTo\pi_1(\nv{\callq})/JC$ as 
follows. Given a morphism $\alpha\in \hom_{\calfq}(P,Q)$, choose a lifting 
$\bar\alpha$ in $\callq$, and set $s(\alpha) = [J(\bar\alpha)]$. If  
$\bar\alpha,\bar\alpha'\in \hom_{\callq}(P,Q)$ are two liftings of  
$\alpha$, then there is an isomorphism $\chi\:P'\to P$ in $\callq$ where 
$P'$ is fully centralized, and an element $g\in C_S(P')$, such that 
$\bar\alpha'=\bar\alpha\circ(\chi\circ\delta_{P'}(g)\circ\chi^{-1})$. 
Since $\chi\circ\delta_{P'}(g)\circ\chi^{-1}$ represents a loop that  
belongs to $JC$, $[J(\bar\alpha)]=[J(\bar\alpha')]$. Thus the definition 
of $s$ does not depend on the choice of $\bar\alpha$. It remains to show 
that $s$ preserves the relations among the generators. But we clearly have 
that $s([J(\incl_P^Q)])=[J(\iota_P^Q)]=1$.  Also, if $\alpha$ and $\beta$ are 
composable morphisms of $\calfq$, and $\bar\alpha$ and $\bar\beta$ are 
liftings to $\callq$, then $\bar\alpha\circ\bar\beta$ is a lifting of 
$\alpha\circ\beta$, and hence
$s(\alpha)s(\beta)=[J(\bar\alpha)]
[J(\bar\beta)]= [J(\bar\alpha\circ\bar\beta)]= s(\alpha\circ\beta)$.

This shows that $\pi_\#$ induces an isomorphim $\pi_1(\nv{\callq})/JC 
\cong\pi_1(\nv{\calfq})$. 
\end{proof}


\newsect{Subsystems with $p$-solvable quotient} 

In this section, we prove some general results about subsystems of 
saturated fusion systems and $p$-local finite subgroups:  subsystems with 
$p$-group quotient or quotient of order prime to $p$.  These will then be 
used in the next two sections to prove some more specific theorems.  

It will be convenient in this section to write ``$p'$-group'' for a finite 
group of order prime to $p$.  Recall that a $p$-solvable group is a group 
$G$ with normal series $1=H_0\nsg{}H_1\nsg\cdots\nsg{}H_k=G$ such that each 
$H_i/H_{i-1}$ is a $p$-group or a $p'$-group.  As one consequence of the 
results of this section, we show that for any $p$-local finite group 
$\SFL$, and any homomorphism $\theta$ from $\pi_1(|\call|)$ to a finite 
$p$-solvable group, there is another $p$-local finite group $\SFL[_0]$ such 
that $|\call_0|$ is homotopy equivalent to the covering space of $|\call|$ 
with fundamental group $\Ker(\theta)$.

We start with some definitions.  Recall that for any finite group $G$, 
$O^{p'}(G)$ and $O^p(G)$ are the smallest normal subgroups of $G$ of index 
prime to $p$ and of $p$-power index, respectively.  Equivalently, 
$O^{p'}(G)$ is the subgroup generated by elements of $p$-power order in 
$G$, and $O^p(G)$ is the subgroup generated by elements of order prime to 
$p$ in $G$.

We want to identify the fusion subsystems of a given fusion system
which are analogous to subgroups of $G$ which contain $O^{p'}(G)$ or 
$O^p(G)$.  This motivates the following definitions.

\begin{Defi} \label{def-p'-index}
Let $\calf$ be a saturated fusion system over a $p$-group $S$.  Let 
$(S',\calf')\subseteq(S,\calf)$ be a saturated fusion subsystem; i.e., 
$S'\le{}S$ is a subgroup, and $\calf'\subseteq\calf$ is a subcategory 
which is a saturated fusion system over $S'$.  
\begin{enumerate}  
\item $(S',\calf')$ is \emph{of $p$-power index} in $(S,\calf)$ if 
$S'\ge\hfocal{\calf}{S}$, and $\Aut_{\calf'}(P)\ge O^p(\autf(P))$ for all 
$P\le{}S'$.  Equivalently, a saturated fusion subsystem 
$\calf'\subseteq\calf$ over $S'\ge\hfocal{\calf}S$ has $p$-power index if 
it contains all $\calf$-automorphisms of order prime to 
$p$ of subgroups of $S'$.

\item $(S',\calf')$ is \emph{of index prime to $p$} in $(S,\calf)$ if 
$S'=S$, and $\Aut_{\calf'}(P)\ge O^{p'}(\autf(P))$ for all $P\le{}S$.  
Equivalently, a saturated fusion subsystem $\calf'\subseteq\calf$ over $S$ 
has index prime to $p$ if it contains all $\calf$-automorphisms of 
$p$-power order.
\end{enumerate}
\end{Defi}

This terminology has been chosen for simplicity.  Subsystems ``of $p$-power 
index'' or ``of index prime to $p$'' are really analogous to subgroups 
$H\le{}G$ which contain \emph{normal subgroups} of $G$ of $p$-power index 
or index prime to $p$, respectively.  For example, if $\calf$ is a 
saturated fusion system over $S$, and $\calf'\subseteq\calf$ is the fusion 
system of $S$ itself (i.e., the minimal fusion system over $S$), then 
$\calf'$ is not in general a subsystem of index prime to $p$ under the 
above definition, despite the inclusion $\calf'\subseteq\calf$ being 
analogous to the inclusion of a Sylow $p$-subgroup in a group.  

Over the next three sections, we will classify all saturated fusion 
subsystems of $p$-power index, or of index prime to $p$, in a given 
saturated fusion system.  In both cases, there will be a minimal such 
subsystem, denoted $\Opf$ or $\Oppf$; and the saturated fusion subsystems 
of the given type will be in bijective correspondence with the subgroups 
of a given $p$-group or $p'$-group.  

The following terminology will be useful for describing some of the 
categories we have to work with.  

\begin{Defi} \label{D:subcats}
Fix a finite $p$-group $S$.  
\begin{enumerate} 
\item A \emph{restrictive category over $S$} is a category $\calf$ such 
that $\Ob(\calf)$ is the set of subgroups of $S$, such that all morphisms 
in $\calf$ are group monomorphisms between the subgroups, and with the 
following additional property:  for each $P'\le{}P\le{}S$ and 
$Q'\le{}Q\le{}S$, and each $\varphi\in\Hom_{\calf}(P,Q)$ such that 
$\varphi(P')\le{}Q'$, $\varphi|_{P'}\in\Hom_{\calf}(P',Q')$.  

\item A restrictive category $\calf$ over $S$ is \emph{normalized} by an 
automorphism $\alpha\in\autf(S)$ if for each $P,Q\le{}S$, and each 
$\varphi\in\Hom_{\calf}(P,Q)$, 
$\alpha\varphi\alpha^{-1}\in\Hom_{\calf}(\alpha(P),\alpha(Q))$.

\item For any restrictive category $\calf$ over $S$ and any subgroup 
$A\le\Aut(S)$, $\gen{\calf,A}$ is the smallest restrictive category over 
$S$ which contains $\calf$ together with all automorphisms in $A$ and 
their restrictions. 
\end{enumerate}
\end{Defi}

By definition, any restrictive category is required to contain all 
inclusion homomorphisms (restrictions of $\Id_S$).  The main difference 
between a restrictive category over $S$ and a fusion system over $S$ 
is that the restrictive subcategory need not contain $\calf_S(S)$.  

\mynote{B: one additional difference is that group isomorphisms in a 
restrictive category need not be isomorphisms in the category.  We could 
include this condition in the definition, but it never seems relevant to 
what we need.}

When $\calf$ is a fusion system over $S$, then an automorphism 
$\alpha\in\Aut(S)$ normalizes $\calf$ if and only if it is \emph{fusion 
preserving} in the sense used in \cite{BLO1}.  For the purposes of this 
paper, it will be convenient to use both terms, in different situations.

We next define the following subcategories of a given fusion system $\calf$.

\begin{Defi} \label{D:oppf:opf}
Let $\calf$ be any fusion system over a $p$-group $S$. 
\begin{enumerate} 
\item $\opf\subseteq\calf$ denotes the smallest restrictive subcategory of 
$\calf$ whose morphism set contains $O^{p}(\autf(P))$ for all 
subgroups $P\le{}S$. 

\item $\oppf\subseteq\calf$ denotes the smallest restrictive subcategory of 
$\calf$ whose morphism set contains $O^{p'}(\autf(P))$ for all 
subgroups $P\le{}S$. 
\end{enumerate}
\end{Defi}

By definition, for subgroups $P,Q$, a morphism $\varphi\in\homf(P,Q)$
lies in $\oppf$ if and only if it is a composite of morphisms 
which are restrictions of elements of $O^{p'}(\autf(R))$ for 
subgroups $R\le{}S$.  Morphisms in $\opf$ are 
described in a similar way.

The subcategory $\opf$ is not, in general, a fusion system --- and this is 
why we had to define restrictive categories.  The subcategory $\oppf$ is 
always a fusion system (since $\Aut_S(P)\le{}O^{p'}(\autf(P))$ for all 
$P\le{}S$), but it is not, in general, saturated.  The subscripts ``$*$'' 
have been put in as a reminder of these facts, and as a contrast with the 
notation $\Opf$ and $\Oppf$ which will be used to denote certain minimal 
saturated fusion systems.  

We now check some of the basic properties of the subcategories $\opf$ 
and $\oppf$, stated in terms of these definitions.

\begin{Lem}\label{F=F'.N_F(S)}
The following hold for any fusion system $\calf$ over a 
$p$-group $S$. 
\begin{enumerate}  
\item  $\opf$ and $\oppf$ are normalized by $\autf(S)$.
\item  If $\calf$ is saturated, then 
$\calf=\gen{\opf,\autf(S)}=\gen{\oppf,\autf(S)}$.
\item If $\calf'\subseteq\calf$ is any restrictive subcategory normalized 
by $\autf(S)$ and such that $\calf=\gen{\calf',\autf(S)}$, then for each 
$P,Q\le{}S$ and $\varphi\in\homf(P,Q)$, there are morphisms 
$\alpha\in\autf(S)$, $\varphi'\in\Hom_{\calf'}(\alpha(P),Q)$, and 
$\varphi''\in\Hom_{\calf'}(P,\alpha^{-1}Q)$, such that 
$\varphi=\varphi'\circ\alpha|_P=\alpha\circ\varphi''$.
\end{enumerate}
\end{Lem}

\begin{proof}  To simplify notation in the following proofs, for 
$\alpha\in\autf(S)$, we write $\alpha$ in composites of morphisms between 
subgroups of $S$, rather than specifying the appropriate restriction each 
time.

\noindent\textbf{(a) }  Set $\calf'=\opf$ or $\oppf$.  
Let $\alpha\in\autf(S)$ and let $\varphi\in\Hom_{\calf'}(P,Q)$. We 
must show that $\alpha\varphi\alpha^{-1}
\in\Hom_{\calf'}(\alpha(P),\alpha(Q))$.  By definition of $\calf'$, there 
are subgroups 
	\[ P=P_0,P_1,\dots,P_k=\varphi(P)\le{}Q, \] 
subgroups $\widebar{P}_1,\dots,\widebar{P}_k$, and automorphisms 
$\chi_i\in{}O^{p}(\autf(\widebar{P}_i))$ (if $\calf'=\opf$) or 
$\chi_i\in{}O^{p'}(\autf(\widebar{P}_i))$ (if $\calf'=\oppf$), such that 
$P_{i-1},P_i\le\widebar{P}_i$, $\chi_i(P_{i-1})=P_i$, and
	$$ \varphi = \incl_{P_k}^{Q}\circ (\chi_k|_{P_{k-1}}) \circ
	\cdots \circ (\chi_2|_{P_1}) \circ(\chi_1|_{P_0}) . $$
Let $P'_i=\alpha(P_i)$, $\widebar{P}'_i=\alpha(\widebar{P}_i)$, and 
$\chi'_i=\alpha\chi_i\alpha^{-1}\in{}O^{p}(\autf(\widebar{P}'_i))$ or
$O^{p'}(\autf(\widebar{P}'_i))$. Then 
	$$ \alpha\varphi\alpha^{-1} =
	\incl_{P'_k}^{\alpha(Q)}\circ (\chi'_k|_{P'_{k-1}}) \circ
	\cdots \circ (\chi'_2|_{P'_1}) \circ(\chi'_1|_{P'_0}) $$
is in $\Hom_{\calf'}(\alpha(P),\alpha(Q))$, as required.

\smallskip

\noindent\textbf{(b1) } Now assume that $\calf$ is saturated.  We show 
here that $\calf=\gen{\oppf,\autf(S)}$; i.e., that each 
$\varphi\in\Mor(\calf)$ is a composite of morphisms in $\oppf$ and in 
$\autf(S)$.  By Alperin's fusion theorem for saturated fusion systems 
(Theorem \ref{Alp.fusion}(a)), it suffices to prove this when 
$\varphi\in\autf(P)$ for some $P\le{}S$ which is $\calf$-centric and fully 
normalized.  The result clearly holds if $P=S$.  So we can assume that 
$P\lneqq{}S$, and also assume inductively that the lemma holds for every 
automorphism of any subgroup $P'\le{}S$ such that $|P'|>|P|$.  

Consider the subgroup
	$$ K = \varphi \Aut_S(P) \varphi^{-1} =
	\{ \varphi c_g\varphi^{-1} \in \autf(P) \,|\, g\in{}N_S(P) \}. $$
Then $K$ is a $p$-subgroup of $\autf(P)$, and hence $K\le 
\Aut_{\oppf}(P)$.  Since $P$ is fully normalized in $\calf$, 
$\Aut_S(P)\in\sylp{\autf(P)}$.  Thus $K$ and $\Aut_S(P)$ are both Sylow 
$p$-subgroups of $\Aut_{\oppf}(P)$, and there is some 
$\chi\in\Aut_{\oppf}(P)$ such that $\chi{}K\chi^{-1}\le\Aut_S(P)$.  In 
other words,
	$$ N_{\chi\varphi} \defeq 
	\{ g\in{}N_S(P) \,|\, (\chi\varphi) c_g(\chi\varphi)^{-1} \in 
	\Aut_S(P) \} = N_S(P). $$
So by condition (II) in Definition \ref{sat.Frob.}, $\chi\circ\varphi$ 
can be extended to a homomorphism $\psi\in\homf(N_S(P),S)$.  By the 
induction hypothesis (and since $N_S(P)\gneqq{}P$), $\psi$, and hence 
$\chi\circ\varphi$, are composites of morphisms in $N_{\calf}(S)$ and in 
$\oppf$.  Also, $\chi\in\Aut_{\oppf}(P)$ by assumption, and hence 
$\varphi$ is a composite of morphisms in these two subcategories.

\smallskip

\noindent\textbf{(b2) } To see that $\calf=\gen{\opf,\autf(S)}$ when 
$\calf$ is saturated, it again suffices to restrict to the case where 
$\varphi\in\autf(P)$ for some $P\le{}S$ which is $\calf$-centric and fully 
normalized.  But in this case, $\autf(P)$ is generated by its Sylow 
$p$-subgroup $\Aut_S(P)\le\Aut_{N_\calf(S)}(P)$, together with 
$O^{p}(\autf(P))\le\Aut_{\opf}(P)$. 

\smallskip

\noindent\textbf{(c) } Since $\calf=\gen{\calf',\autf(S)}$, every morphism 
in $\calf$ is a composite of morphisms in $\calf'$ and restrictions of 
automorphisms in $\autf(S)$.  Assume $\varphi\in\homf(P,Q)$ is the 
composite
	$$ P = P_0 \Right2{\alpha_0} Q_0 \Right2{\varphi_1} P_1 
	\Right2{\alpha_1} Q_1 \Right2{} \ldots \Right2{} P_n 
	\Right2{\alpha_n} Q_n = Q, $$
where $\alpha_i\in\autf(S)$, $\alpha_i(P_i)=Q_i$, and 
$\varphi_i\in\Hom_{\calf'}(Q_{i-1},P_i)$.  Write 
$\alpha_{j,i}=\alpha_j\cdots\alpha_i$ for any $i\le{}j$, and set
$\alpha=\alpha_{n,0}=\alpha_n\cdots\alpha_0$.  Then 
	$$ \varphi = \alpha\circ 
	(\alpha_{n-1,0}{}^{-1}\varphi_n\alpha_{n-1,0}) \cdots
	(\alpha_{0,0}{}^{-1}\varphi_1\alpha_{0,0}) 
	= (\alpha_{n,n}\varphi_n\alpha_{n,n}{}^{-1})\cdots
	(\alpha_{n,1}\varphi_1\alpha_{n,1}{}^{-1})\circ\alpha; $$
where each $\alpha_{n,i}\varphi_i\alpha_{n,i}{}^{-1}$ and each 
$\alpha_{i,0}{}^{-1}\varphi_{i+1}\alpha_{i,0}$ is a morphism in $\calf'$ 
since $\calf'$ is normalized by $\autf(S)$.  
\end{proof}

The following lemma will also be needed later.  

\begin{Lem} \label{P0:F0-centric}
Let $\calf$ be a saturated fusion system over a $p$-group $S$.  Fix a 
normal subgroup $S_0\nsg{}S$ which is strongly $\calf$-closed; i.e., no 
element of $S_0$ is $\calf$-conjugate to any element of $S{\sminus}S_0$.  
Let $(S_0,\calf_0)$ be a saturated fusion subsystem of $(S,\calf)$.  Then 
for any $P\le S$ which is $\calf$-centric and $\calf$-radical, $P\cap{}S_0$ 
is $\calf_0$-centric.  
\end{Lem}

\begin{proof}  Assume $P\le S$ is $\calf$-centric and $\calf$-radical, and 
set $P_0=P\cap{}S_0$ for short.  Choose a subgroup $P'_0\le{}S_0$ which is 
$\calf$-conjugate to $P_0$ and fully normalized in $\calf$.  In particular, 
by (I), $P'_0$ is fully centralized in $\calf$.  By Lemma 
\ref{N->N}, there is $\varphi\in\homf(N_S(P_0),N_S(P'_0))$ such that 
$\varphi(P_0)=P'_0$.  Set $P'=\varphi(P)$; thus $P'$ is $\calf$-conjugate to $P$ 
(so $P'$ is also $\calf$-centric and $\calf$-radical), and 
$P'_0=P'\cap{}S_0$ since $S_0$ is strongly closed.  For any 
$P''_0\le{}S_0$ which is $\calf_0$-conjugate to $P_0$ (hence 
$\calf$-conjugate to $P'_0$), there is a morphism 
$\psi\in\homf(P''_0{\cdot}C_S(P''_0),P'_0{\cdot}C_S(P'_0))$ such that 
$\psi(P''_0)=P'_0$ (by axiom (II)), and then 
$\psi(C_{S_0}(P''_0))\le{}C_{S_0}(P'_0)$.  So if 
$C_{S_0}(P'_0)=Z(P'_0)$, then $C_{S_0}(P''_0)=Z(P''_0)$ for all $P''_0$ 
which is $\calf_0$-conjugate to $P_0$, and $P_0$ is $\calf_0$-centric.  
 
We are thus reduced to showing that $C_{S_0}(P'_0)=Z(P'_0)$; and without 
loss of generality, we can assume that $P'=P$ and $P'_0=P_0$.  Since
$S_0$ is strongly closed, every $\alpha\in\autf(P)$ leaves $P_0$ 
invariant.  Let $\autf^0(P)\le\autf(P)$ be the subgroup of elements which 
induce the identity on $P_0$ and on $P/P_0$.  This is a normal subgroup of 
$\autf(P)$ since all elements of $\autf(P)$ leave $P_0$ invariant, and is 
also a $p$-subgroup by Lemma \ref{gorenstein}.  Thus 
$\autf^0(P)\le{}O_p(\autf(P))$, and hence $\autf^0(P)\le\Inn(P)$ since $P$ 
is $\calf$-radical.  

We want to show that $C_{S_0}(P_0)=Z(P_0)$.  Fix any $x\in{}C_{S_0}(P_0)$, 
and assume first that the coset $x{\cdot}Z(P_0)\in C_{S_0}(P_0)/Z(P_0)$ is 
fixed by the conjugation action of $P$.  Thus $x\in{}S_0$, $[x,P_0]=1$, 
and $[x,P]\le{}Z(P_0)$, so $c_x\in\autf^0(P)\le\Inn(P)$, and 
$xg\in{}C_S(P)$ for some $g\in{}P$.  Since $P$ is $\calf$-centric, this 
implies that $xg\in{}P$, so $x\in{}C_{S_0}(P_0)\cap{}P=Z(P_0)$.  In other 
words, $[C_{S_0}(P_0)/Z(P_0)]^P=1$, so $C_{S_0}(P_0)/Z(P_0)=1$, and thus 
$C_{S_0}(P_0)=Z(P_0)$.  
\end{proof}

The motivation for the next definition comes from considering the 
situation which arises when one is given a saturated fusion system $\calf$ 
with an associated quasicentric linking system $\callq$, and a functor 
$\callq\rTo\calb(\Gamma)$ which sends inclusions to the identity (or 
equivalently a homomorphism from $\pi_1(|\call^q|)$ to $\Gamma$) for some 
group $\Gamma$. Such a functor is equivalent to a function 
$\widebar{\Theta}\:\Mor(\callq)\rTo\Gamma$ which sends composites to 
products and sends inclusions to the identity; and for any $H\le\Gamma$, 
$\widebar{\Theta}^{-1}(H)=\Mor(\callq_H)$ for some subcategory 
$\callq_H\subseteq\callq$ with the same objects.  Let 
$\calfq_H\subseteq\calfq$ be the image of $\callq_H$ under the canonical 
projection; then in some sense (to be made precise later), $\calfq_H$ is a 
subsystem of $\calfq$ of index $[\Gamma{:}H]$.  

What we now need is to make sense of such ``inverse image subcategories'' 
of the fusion system $\calf$, when we are not assuming that we have an 
associated linking system.  Let $\sset(\Gamma)$ denote the set of 
\emph{nonempty subsets} of $\Gamma$. Given a function $\widebar{\Theta}$ 
as above, there is an obvious associated function 
$\Theta\:\Mor(\calfq)\rTo\sset(\Gamma)$, which sends a morphism 
$\alpha\in\Mor(\calfq)$ to $\widebar{\Theta}(\pi^{-1}(\alpha))$.  Here, 
$\pi$ denotes the natural projection from $\callq$ to $\calfq$.  Moreover, 
$\widebar{\Theta}$ also induces a homomorphism 
$\theta=\widebar{\Theta}\circ\delta_S$ from $S$ to $\Gamma$.  The maps 
$\Theta$ and $\theta$ are closely related to each other, and satisfy 
certain properties, none of which depend on the existence (or choice) of a 
quasicentric linking system associated to $\calf$. In fact, we will see 
that the data encoded in such a pair of functions, if they satisfy the 
appropriate conditions, suffices to describe precisely what is meant by 
``inverse image subcategories'' of fusion systems; and to show that under 
certain restrictions, these categories are (or generate) saturated fusion 
subsystems.

\begin{Defi} \label{D:J,j}
Let $\calf$ be a saturated fusion system over a $p$-group $S$, and let 
$\calf_0\subseteq\calfq$ be any full subcategory such that $\Ob(\calf_0)$ 
is closed under $\calf$-conjugacy.  A \emph{fusion mapping triple} for 
$\calf_0$ consists of a triple $(\Gamma,\theta,\Theta)$, where $\Gamma$ is 
a group, $\theta\:S\Right1{}\Gamma$ is a homomorphism, and 
	$$ \Theta\: \Mor(\calf_0) \Right4{} \sset(\Gamma), $$
is a map which satisfies the following conditions for all 
subgroups $P,Q,R\le{}S$ which lie in $\calf_0$:
\begin{enumerate}\renewcommand{\labelenumi}{\textup{(\roman{enumi})}}

\item For all $P\Right1{\varphi}Q\Right1{\psi}R$ in $\calf_0$, and all 
$x\in\Theta(\psi)$, $\Theta(\psi\varphi)=x{\cdot}\Theta(\varphi)$.

\item If $P$ is fully centralized in $\calf$, then 
$\Theta(\Id_P)=\theta(C_S(P))$.  

\item If $\varphi=c_g\in\homf(P,Q)$, where $g\in N_S(P,Q)$, then 
$\theta(g)\in\Theta(\varphi)$.

\item For all $\varphi\in\homf(P,Q)$, all $x\in\Theta(\varphi)$, and all 
$g\in{}P$, $x\theta(g)x^{-1}=\theta(\varphi(g))$.
\end{enumerate}
For any fusion mapping triple $(\Gamma,\theta,\Theta)$ and any $H\le\Gamma$, 
we 
let $\lfx{\Theta}H\subseteq\calf$ be the smallest restrictive subcategory 
which contains all $\varphi\in\Mor(\calfq)$ such that 
$\Theta(\varphi)\cap{}H\ne\emptyset$.  Let 
$\lf{\Theta}H\subseteq\lfx{\Theta}H$ be the full subcategory whose objects 
are the subgroups of $\theta^{-1}(H)$.  
\end{Defi}

When $\theta$ is the trivial homomorphism (which is always the 
case when $\Gamma$ is a $p'$-group), then a fusion mapping triple 
$(\Gamma,\theta,\Theta)$ on a subcategory $\calf_0\subseteq\calfq$ is 
equivalent to a functor from $\calf_0$ to $\calb(\Gamma)$; i.e., 
$\Theta(\varphi)$ contains just one element for all 
$\varphi\in\Mor(\calf_0)$.  By (i), it suffices to show this for identity 
morphisms; and by (ii), $|\Theta(\Id_P)|=1$ if $P$ is fully centralized.  
For arbitrary $P$, it then follows from (i), together with the assumption 
(included in the definition of $\sset(\Gamma)$) that 
$\Theta(\varphi)\ne\emptyset$ for all $\varphi$.

The following additional properties of fusion mapping triples will be 
needed.

\begin{Lem} \label{P:Theta-theta}
Fix a saturated fusion system $\calf$ over a $p$-group $S$, let $\calf_0$ 
be a full subcategory such that $\Ob(\calf_0)$ is closed under 
$\calf$-conjugacy, and let $(\Gamma,\Theta,\theta)$ be a fusion mapping 
triple for $\calf_0$.  Then the following hold for all 
$P,Q,R\in\Ob(\calf_0)$.
\begin{enumerate}\renewcommand{\labelenumi}{\textup{(\roman{enumi})}}%
\setcounter{enumi}{4}
\item $\Theta(\Id_P)$ is a subgroup of $\Gamma$, and $\Theta$ restricts to 
a homomorphism 
	$$ \Theta_P\:\autf(P)\rTo N_\Gamma(\Theta(\Id_P))/\Theta(\Id_P). $$ 
Thus $\Theta_P(\alpha)=\Theta(\alpha)$ (as a coset of $\Theta(\Id_P)$) 
for all $\alpha\in\autf(P)$.  

\item For all $P\Right1{\varphi}Q\Right1{\psi}R$ in $\calf_0$, and all 
$x\in\Theta(\varphi)$, $\Theta(\psi\varphi)\supseteq\Theta(\psi){\cdot}x$, 
with equality if $\varphi(P)=Q$.  In particular, if $P\le{}Q$, then 
$\Theta(\psi|_P)\supseteq\Theta(\psi)$.  

\item Assume $S\in\Ob(\calf_0)$.  Then for any $\varphi\in\homf(P,Q)$, any 
$\alpha\in\autf(S)$, and any $x\in\Theta(\alpha)$, 
$\Theta(\alpha\varphi\alpha^{-1})=x\Theta(\varphi)x^{-1}$, where 
$\alpha\varphi\alpha^{-1}\in\homf(\alpha(P),\alpha(Q))$. 
\end{enumerate}
\end{Lem}

\begin{proof}  \textbf{(v) } 
By (i), for any $\alpha,\beta\in\autf(P)$ and any $x\in\Theta(\alpha)$, 
$\Theta(\alpha\beta)=x{\cdot}\Theta(\beta)$.  When applied with 
$\alpha=\beta=\Id_P$, this shows that $\Theta(\Id_P)$ is a subgroup of 
$\Gamma$.  (Note here that $\Theta(\Id_P)\ne\emptyset$ by definition of 
$\sset(\Gamma)$.)  When applied with $\beta=\alpha^{-1}$, this shows that 
$x^{-1}\in\Theta(\alpha^{-1})$ if $x\in\Theta(\alpha)$.  Hence 
$\Theta(\alpha)=x{\cdot}\Theta(\Id_P)$ implies that 
$\Theta(\alpha^{-1})=\Theta(\Id_P){\cdot}x^{-1}$ and 
$\Theta(\alpha^{-1})=x^{-1}{\cdot}\Theta(\Id_P)$, and thus that each 
$\Theta(\alpha)$ is a right coset as well as a left coset.  Thus 
$\Theta(\alpha)\subseteq{}N_\Gamma(\Theta(\Id_P))$ for all 
$\alpha\in\autf(P)$, and the induced map $\Theta_P$ is clearly a 
homomorphism.

\noindent\textbf{(vi) }  By (i), 
$\Theta(\psi\varphi)\supseteq\Theta(\psi){\cdot}\Theta(\varphi)$ for 
any pair of composable morphisms $\varphi,\psi$ in $\calf_0$.  In 
particular, $\Theta(\psi\varphi)\supseteq\Theta(\psi){\cdot}x$ if 
$x\in\Theta(\varphi)$.  If $\varphi$ is an isomorphism, then 
$1\in\Theta(\Id_{P})=x{\cdot}\Theta(\varphi^{-1})$ by (v) and (i), so 
$x^{-1}\in\Theta(\varphi^{-1})$.  This gives the inclusions
	$$ \Theta(\psi)=\Theta(\psi\varphi\varphi^{-1}) 
	\supseteq \Theta(\psi\varphi){\cdot}x^{-1}
	\supseteq \Theta(\psi)xx^{-1}, $$
and hence these are both equalities.  The last statement is the special 
case where $P\le{}Q$ and $\varphi=\incl_P^Q$; $1\in\Theta(\incl_P^Q)$ by 
(iii). 

\noindent\textbf{(vii) }  For $x\in\Theta(\alpha)$, 
$\Theta(\alpha\varphi)=x{\cdot}\Theta(\varphi) = 
\Theta(\alpha\varphi\alpha^{-1}){\cdot}x$ by (i) and (vi). 
\end{proof}

We are now ready to prove the main result about fusion mapping triples.

\begin{Prop}\label{F:p-solv.quot.}
Let $\calf$ be a saturated fusion system over a finite $p$-group $S$.  Let 
$(\Gamma,\theta,\Theta)$ be any fusion mapping triple on $\calfq$, where 
$\Gamma$ is a $p$-group or a $p'$-group, and 
	$$ \theta \: S \Right3{} \Gamma 
	\qquad\textup{and}\qquad
	\Theta\: \Mor(\calfq) \Right3{} \sset(\Gamma), $$
Then the following hold for any subgroup $H\le\Gamma$, where we set 
$S_H=\theta^{-1}(H)$.
\begin{enumerate}

\item $\lf{\Theta}H$ is a saturated fusion system over $S_H$.

\item If $\Gamma$ is a $p$-group, then a subgroup 
$P\le{}S_H$ is $\lf{\Theta}H$-quasicentric if and only if it is 
$\calf$-quasicentric.  Also, $\lfx{\Theta}{}\supseteq\opf$.

\item If $\Gamma$ is a $p'$-group, then $S_H=S$.  A subgroup $P\le{}S$ is 
$\lf{\Theta}H$-centric (fully centralized in $\lf{\Theta}H$, fully 
normalized in $\lf{\Theta}H$) if and only if it is $\calf$-centric (fully 
centralized in $\calf$, fully normalized in $\calf$).  Also, 
$\lfx{\Theta}{}\supseteq\oppf$.
\end{enumerate}
\end{Prop}

\begin{proof}  Throughout the proof, (i)--(vii) refer to the conditions in 
Definition \ref{D:J,j} and Lemma \ref{P:Theta-theta}.  If $X$ is any set 
of morphisms of $\calf$, then we let $\Theta(X)$ be the union of the sets 
$\Theta(\alpha)$ for $\alpha\in{}X$. Condition (v) implies that for any 
$P\le{}S$ and any subgroup $A\le\autf(P)$, $\Theta(A)$ is a subgroup of 
$\Gamma$.  

We first prove the following two additional properties of these 
subcategories:
\begin{enumerate}
\item[(1)] For each pair of subgroups $P,Q\le{}S$, 
and each $\varphi\in\homf(P,Q)$, there are $Q'\le{}S$, 
$\varphi'\in\Hom_{\lfx{\Theta}{}}(P,Q')$, and $\alpha\in\Aut_\calf(S)$ such 
that $\alpha(Q')=Q$ and $\varphi=(\alpha|_{Q'})\circ\varphi'$.  
\item[(2)] For all $P\le{}S$ there exists $P'\le{}S$ which is fully 
normalized in $\calf$, and 
$\varphi\in\Hom_{\lfx{\Theta}{}}(N_S(P),N_S(P'))$ such that 
$\varphi(P)=P'$.  
\end{enumerate}

By (vii), for any $\alpha\in\autf(S)$, any $x\in\Theta(\alpha)$, and any 
$\varphi\in\Hom_{\calf}(P,Q)$, $\Theta(\alpha\varphi\alpha^{-1})= 
x\Theta(\varphi)x^{-1}$.  In particular, 
$\Theta(\alpha\varphi\alpha^{-1})=1$ if and only if $\Theta(\varphi)=1$, 
and thus $\lfx{\Theta}{}$ is normalized by $\autf(S)$.  

Let $\oppfc\subseteq\oppf$ and $\opfc\subseteq\opf$ be the full 
subcategories whose objects are the $\calf$-centric subgroups of $S$.  By 
(v), for each $\calf$-quasicentric subgroup $P\le{}S$, 
$\Aut_{\lfx{\Theta}{}}(P)=\Ker(\Theta_P)$ where $\Theta_P$ is a 
homomorphism to a subquotient of $\Gamma$.  Hence 
$\Aut_{\lfx{\Theta}{}}(P)$ contains $O^{p'}(\autf(P))$ if $\Gamma$ is a 
$p'$-group or $O^p(\autf(P))$ if $\Gamma$ is a $p$-group.  This shows that
$\lfx{\Theta}{}$ contains either $\oppfc$ (if $\Gamma$ is a $p'$-group) or 
$\opfc$ (if $\Gamma$ is a $p$-group).  Thus 
$\calf^c\subseteq\gen{\lfx{\Theta}{},\autf(S)}$ by Lemma 
\ref{F=F'.N_F(S)}(b), and so $\calf=\gen{\lfx{\Theta}{},\autf(S)}$ since 
$\calf$ is the smallest restrictive category over $S$ which contains 
$\calf^c$ by Theorem \ref{Alp.fusion}(a) (Alperin's fusion theorem). Point 
(1) now follows from Lemma \ref{F=F'.N_F(S)}(c).  

To see (2), recall that by Lemma \ref{N->N}, if $P''$ is any subgroup 
$\calf$-conjugate to $P$ and fully normalized in $\calf$, then there is a 
morphism $\varphi'\in\homf(N_S(P),N_S(P''))$ such that $\varphi'(P)=P''$.  
By (1), $\varphi'=\alpha\circ\varphi$ for some $\alpha\in\autf(S)$ and 
some $\varphi\in\Hom_{\lfx{\Theta}{}}(N_S(P),N_S(\alpha^{-1}(P'')))$.  
Also, the subgroup $P'=\alpha^{-1}(P'')$ is fully normalized in $\calf$ 
since $P''$ is.

\smallskip

\noindent\textbf{(b) }  Assume $\Gamma$ is a $p$-group, and fix 
$H\le{}\Gamma$.  Consider the following set of subgroups of $S_H$:
	\begin{multline*}  
	\calq_H = \bigl\{ P\le S_H \,\big|\, 
	\forall\,P' \textup{ $\calf_H$-conjugate to $P$, }
	\forall\,P'\le Q\le P'{\cdot}C_{S_H}(P'),\    \\
	\forall\,\alpha\in\Aut_{\calf_H}(Q) \textup{ such that 
	$\alpha|_{P'}=\Id$, } |\alpha| \textup{ is a power of $p$} 
	\bigr\}. 
	\end{multline*}
This set clearly contains all $\calf_H$-centric subgroups of $S_H$.  By 
Lemma \ref{L:qcentric}, if $\calf_H$ is saturated, then $\calq_H$ is 
precisely the set of $\calf_H$-quasicentric subgroups.  We prove (b) 
here with ``$\calf_H$-quasicentric'' replaced by ``element of $\calq_H$''. 
This is all that will be needed for the proof that $\calf_H$ is saturated;
and once that is shown then (b) (in its original form) will follow. 

Assume $P\le{}S_H$ and $P\notin\calq_H$.  Fix $P'$ which is 
$\lf{\Theta}H$-conjugate to $P$, $Q\le{}P'{\cdot}C_{S_H}(P')$ which 
contains $P'$, and $\Id_Q\ne\alpha\in\Aut_{\lf{\Theta}H}(Q)$ of order 
prime to $p$, such that $\alpha|_{P'}=\Id_{P'}$.  Then $P'$, and hence 
$P$, is not $\calf$-quasicentric by Lemma \ref{L:qcentric}(a).  

We have shown that if $P$ is $\calf$-quasicentric, then $P\in\calq_H$, and 
it remains to check the converse. Assume $P$ is not $\calf$-quasicentric:  
fix $P'\le{}Q\le{}P'{\cdot}C_S(P')$ and $\alpha\in\autf(Q)$ as in Lemma 
\ref{L:qcentric}(b).  In particular, $Q$ is $\calf$-centric, and hence 
$\calf$-quasicentric.  Set $Q_1=Q\cap{}S_1$ (where $S_1=\theta^{-1}(1)$).  
Since $1\in\Theta(\alpha)$ by (v) (and since $|\alpha|$ is prime to $p$), 
$\theta(g)=\theta(\alpha(g))$ for all $g\in{}Q$ by (iv), and thus 
$\alpha(Q_1)=Q_1$ and (since $g^{-1}\alpha(g)\in\Ker(\theta)=S_1$ for 
$g\in{}Q$) $\alpha$ induces the identity on $Q/Q_1$.  Since $|\alpha|$ is 
not a power of $p$, it cannot be the identity on both $Q_1$ and $Q/Q_1$ 
(Lemma \ref{gorenstein}), and hence $\alpha|_{Q_1}\ne\Id_{Q_1}$.  Thus 
$P'\le{}P'Q_1\le{}P'{\cdot}C_{S_H}(P')$, $\alpha|_{P'Q_1}$ is a nontrivial 
automorphism of $P'Q_1$ of order prime to $p$ whose restriction to $P'$ is 
the identity.  

Finally, by (1), $P$ is $\lf{\Theta}H$-conjugate to a subgroup $P''$ for 
which there is some $\varphi\in\isof(P',P'')$ which is the restriction of 
some $\widebar{\varphi}\in\autf(S)$. Since $P\le{}S_H$ and is 
$\lf{\Theta}H$-conjugate to $P''$, $P''\le{}S_H$ by (iv) (applied with 
$\alpha\in\Iso_{\lf{\Theta}H}(P,P'')$ and $x\in\Theta(\alpha)\cap{}H$).  
Set $Q''=\widebar{\varphi}(P'Q_1)=P''{\cdot}\widebar\varphi(Q_1)\le{}S_H$ 
and $\alpha''=\widebar\varphi\alpha\widebar\varphi^{-1}\in\autf(Q'')$.  
Since $\Gamma$ is a $p$-group and $|\alpha''|$ is prime to $p$, (v) 
implies that $\Theta(\alpha'')=\Theta(\Id_Q)$, and hence that 
$1\in\Theta(\alpha'')$ and hence $\alpha''\in\Aut_{\lf{\Theta}H}(Q'')$.  
But then $P''\notin\calq_H$, and hence $P\notin\calq_H$.

It remains to show that $\lfx{\Theta}{}\supseteq\opf$; i.e., to show that
	\beq \Aut_{\lfx{\Theta}{}}(P)\ge O^p(\autf(P)) \tag{3} \eeq
for each $P\le{}S$.  If $P\le{}S$ is $\calf$-quasicentric, then (3) holds 
by (v):  $\Aut_{\lfx{\Theta}{}}(P)$ is the kernel of a homomorphism from 
$\autf(P)$ to a $p$-group.  If $P$ is not $\calf$-quasicentric but is fully 
centralized in $\calf$, then every automorphism $\alpha\in\autf(P)$ of 
order prime to $p$ extends to an automorphism 
$\widebar{\alpha}\in\autf(P{\cdot}C_S(P))$, which (after replacing it by 
an appropriate power) can also be assumed to have order prime to $p$, and 
hence in $O^p(\autf(P{\cdot}C_S(P)))$.  Thus (3) holds for $P$, since it 
holds for the $\calf$-centric subgroup $P{\cdot}C_S(P)$.  Finally, by (1), 
every subgroup of $S$ is $\lfx{\Theta}{}$-conjugate to a subgroup which is 
fully centralized in $\calf$, and thus (3) holds for all $P\le{}S$.

\smallskip

\noindent\textbf{(c) } This point holds, in fact, without assuming that 
the fusion system $\calf$ be saturated.  Assume $\Gamma$ is a $p'$-group;
then $\theta$ is the trivial homomorphism, and so $S_H=\theta^{-1}(H)=S$ 
for all $H$.  Fix $H\le\Gamma$, and let $P$ be any 
subgroup of $S$.  Since each $\lf{\Theta}H$-conjugacy class is contained 
in some $\calf$-conjugacy class, any subgroup which is fully centralized 
(fully normalized) in $\calf$ is also fully centralized (fully  
normalized) in $\lf{\Theta}H$.  By the same reasoning, any $\calf$-centric 
subgroup $P\le S$ is also $\lf{\Theta}H$-centric.

Conversely, assume $P$ is not fully centralized in $\calf$, and let $P'$ 
be a  subgroup in the $\calf$-conjugacy class of $P$ such that 
$|C_S(P')|>|C_S(P)|$.  Fix some $\varphi\in\isof(P,P')$.  By (1), there are 
a subgroup $P''\le{}S$ and isomorphisms $\alpha\in\autf(S)$ and 
$\varphi'\in\Iso_{\lf{\Theta}{}}(P,P'')$, such that $\alpha^{-1}(P'')=P'$ 
and $\varphi=(\alpha|_{P''})\circ\varphi'$.  Thus $P''$ is 
$\lf{\Theta}H$-conjugate to $P$, and is $\calf$-conjugate to $P'$ via a 
restriction of an $\calf$-automorphism of $S$.  Hence 
$|C_S(P'')|=|C_S(P')|>|C_S(P)|$, and this shows that $P$ is not fully 
centralized in $\lf{\Theta}H$. A similar argument shows that if $P$ is not 
fully normalized in $\calf$ (or not $\calf$-centric), then it is not fully 
normalized in $\lf{\Theta}H$ (or not $\lf{\Theta}H$-centric).

The proof that $\lfx{\Theta}{}\supseteq\oppf$ is identical to the proof of 
the corresponding result in (b).

\smallskip

\noindent\textbf{(a) } Fix $H\le{}\Gamma$.  Clearly, $\lf{\Theta}H$ is a 
fusion system over $S_H$; we must prove it is saturated.  By (b) or (c), 
each $\lf{\Theta}H$-centric subgroup of $S_H$ is $\calf$-quasicentric.  
Hence by Theorem \ref{Alp.fusion}(b), it suffices to prove conditions (I) 
and (II) in Definition \ref{sat.Frob.} for $\calf$-quasicentric subgroups 
$P\le{}S_H$.  Thus we will be working only with subgroups $P\le{}S_H$ for 
which $\Theta$ is defined on $\homf(P,S)$.

\smallskip

\noindent\textbf{Proof of (I): }  Assume $\Gamma$ is a $p'$-group; thus 
$S_H=S$.  If $P$ is fully normalized in $\lf{\Theta}H$, then by (c), $P$ is 
fully normalized in $\calf$, hence it is fully centralized in $\calf$ and 
in $\lf{\Theta}H$.  Also, $\Aut_S(P)\in\sylp{\autf(P)}$, and hence 
$\Aut_S(P)$ is also a Sylow $p$-subgroup of $\Aut_{\lf{\Theta}H}(P)$.  This 
proves (I) in this case.

Now assume $\Gamma$ is a $p$-group.  Fix $P\le{}S_H$ which is 
$\calf$-quasicentric and fully normalized in $\lf{\Theta}H$, and let 
$\alpha\in\Hom_{\lfx{\Theta}{}}(N_S(P),N_S(P'))$ be as in (2).  In 
particular, $1\in\Theta(\alpha)$ by definition of $\calf_1^\bullet$; so by 
(iv), $\theta(\alpha(N_{S_H}(P)))=\theta(N_{S_H}(P))\le{}H$.  Hence 
$\alpha(N_{S_H}(P))\le{}S_H\cap{}N_S(P')=N_{S_H}(P')$, and this is an 
equality since $P$ is fully normalized in $\lf{\Theta}H$. In particular, 
this shows that $P'\le{}S_H$.  The conclusion of (I) holds for $P$ (i.e., 
$P$ is fully centralized in $\lf{\Theta}H$ and $\Aut_{S_H}(P)\in 
\sylp{\Aut_{\lf{\Theta}H}(P)}$) if the conclusion of (I) holds for $P'$.  
So we can assume that $P=P'$ is fully normalized in $\calf$ and in 
$\lf{\Theta}H$.  

Fix $Q\le{}S_H$ which is $\lf{\Theta}H$-conjugate to $P$ and fully 
centralized in $\lf{\Theta}H$, and choose an isomorphism 
$\psi\in\Iso_{\lf{\Theta}H}(P,Q)$.  After applying (2) again, we can 
assume that $Q$ is also fully normalized in $\calf$, and hence also fully 
centralized.  Hence $\psi$ extends to some 
$\psi_1\in\homf(P{\cdot}C_S(P),Q{\cdot}C_S(Q))$, which is an isomorphism 
since $P$ and $Q$ are both fully centralized in $\calf$.  Fix 
$h\in\Theta(\psi^{-1})\cap{}H$ and $g\in\Theta(\psi_1)$.  Then 
$gh\in\Theta(\Id_Q)=\theta(C_S(Q))$ by (i) and (ii).  Let $a\in{}C_S(Q)$ 
be such that $\theta(a)=gh$, and set 
	$$ \psi_2=c_a^{-1}\circ\psi_1
	\in\isof(P{\cdot}C_S(P),Q{\cdot}C_S(Q)). $$  
Thus $\psi_2|_P=\psi_1|_P=\psi$ since $a\in{}C_S(Q)$, and 
$h^{-1}=\theta(a)^{-1}g\in\Theta(\psi_2)$ by (iii) and (i).  By (iv), for 
all $g\in{}P{\cdot}C_S(P)$, $\theta(\psi_2(g))=h^{-1}\theta(g)h$, and hence
$\psi_2(g)\in{}H$ if and only if $g\in{}H$.  Thus $\psi_2$ sends 
$C_{S_H}(P)$ onto $C_{S_H}(Q)$.  Since $Q$ is fully centralized in 
$\lf{\Theta}H$, so is $P$. 

It remains to show that $\Aut_{S_H}(P)\in\sylp{\Aut_{\lf{\Theta}H}(P)}$.  
Set $P_\Theta=\Theta(\Id_P)=\theta(C_S(P))$ for short (see (ii)).  By (v), 
$\Theta$ restricts to a homomorphism 
	$$ \Theta_P \: \autf(P) \Right5{} N_\Gamma(P_\Theta)/P_\Theta. $$
Set
	$$ \widehat{H} = \bigl\{ hP_\Theta \,\big|\, 
	h\in{}H\cap{}N_\Gamma(P_\Theta) \bigr\}
	\le N_\Gamma(P_\Theta)/P_\Theta; $$
then $\Aut_{\lf{\Theta}H}(P)=\Theta_P{}^{-1}(\widehat{H})$ by definition 
of $\lf{\Theta}H$.  For any 
$\varphi=c_a\in\Aut_S(P)\cap\Aut_{\lf{\Theta}H}(P)$ (where 
$a\in{}N_S(P)$), $\Theta(\varphi)=\theta(a){\cdot}\theta(C_S(P))$ by (ii) 
and (iii); and thus $\Theta(\varphi)$ contains some $h\in{}H$ where 
$h=\theta(b)$ for some $b\in{}a{\cdot}C_S(P)$.  Hence $b\in{}S_H$, and 
$\varphi=c_b\in\Aut_{S_H}(P)$.  This shows that
	$$ \Aut_{S_H}(P) = \Aut_S(P)\cap\Aut_{\lf{\Theta}H}(P)
	=\Aut_S(P)\cap\Theta_P^{-1}(\widehat{H}). $$
Also, $\Theta_P(\Aut_S(P))=\Theta_P(\autf(P))$, since 
$\Aut_S(P)\in\sylp{\autf(P)}$ and $\Gamma$ is a $p$-group, and hence
	$$ [\Aut_S(P):\Aut_{S_H}(P)] = [\Im(\Theta_P):\widehat{H}] = 
	[\autf(P):\Aut_{\lf{\Theta}H}(P)]. $$
Since $\Aut_S(P)\in\sylp{\autf(P)}$ (i.e., $[\autf(P):\Aut_S(P)]$ is prime 
to $p$), this implies that $\Aut_{S_H}(P)\in\sylp{\Aut_{\lf{\Theta}H}(P)}$.  

\smallskip

\noindent\textbf{Proof of (II): }  Fix a morphism 
$\varphi\in\Iso_{\lf{\Theta}H}(P,Q)$, for some $\calf$-quasicentric 
$P,Q\le{}S_H$ such that $Q=\varphi(P)$ is fully centralized in 
$\lf{\Theta}H$, and set $N_\varphi= 
\{g\in{}N_{S_H}(P)\,|\,\varphi{}c_g\varphi^{-1}\in\Aut_{S_H}(Q)\}$.  By 
(2), there is a subgroup $Q'$ fully normalized in $\calf$, and a morphism 
$\psi\in\Hom_{\lfx{\Theta}{}}(N_S(Q),N_S(Q'))$, such that $\psi(Q)=Q'$. By 
condition (II) for the saturated fusion system $\calf$, there is 
$\widebar{\varphi}_1\in\homf(N_\varphi,N_S(Q'))$ such that 
$\widebar{\varphi}_1|_P=\psi\circ\varphi$.  Fix some
	$$ x\in\Theta(\widebar{\varphi}_1) \subseteq\Theta(\psi\varphi)
	=\Theta(\varphi) $$
(where the last equality holds since $1\in\Theta(\psi)$).  Since 
$\varphi\in\Mor(\lf{\Theta}H)$, there is $h\in{}H$ such that 
$h\in\Theta(\varphi)$, and thus 
	$$ hx^{-1} \in \Theta(\Id_{Q'}) = \theta(C_S(Q')). $$
Fix $a\in{}C_S(Q')$ such that $hx^{-1}=\theta(a)$, and set 
	$$ \widebar{\varphi}_2 =c_a\circ\widebar{\varphi}_1 
	\in\homf(N_\varphi,N_S(Q')). $$
Then by (i), $h\in\Theta(\widebar{\varphi}_2)$, so $\widebar{\varphi}_2 
\in\Hom_{\lf{\Theta}H}(N_\varphi,N_{S_H}(Q'))$; and 
$\widebar{\varphi}_2|_P=\widebar{\varphi}_1|_P$ since $a\in{}C_S(Q')$.  
Since $Q$ is fully centralized in $\lf{\Theta}H$, $\psi$ sends $C_{S_H}(Q)$ 
isomorphically onto $C_{S_H}(Q')$; and hence (by definition of 
$N_\varphi$) $\widebar{\varphi}_2(N_\varphi)\le\Im(\psi)$.  So 
$\widebar{\varphi}\defeq\psi^{-1}\circ\widebar{\varphi}_2$ sends 
$N_\varphi$ into $N_{S_H}(Q)$ and extends $\varphi$.
\end{proof}

We next extend Proposition \ref{F:p-solv.quot.} to a result about linking 
systems and $p$-local finite groups.  The main point of the following 
theorem is that for any $p$-local finite group $\SFL$ and any epimorphism 
$\pi_1(|\call|)\Onto2{\widehat{\theta}}\Gamma$, where $\Gamma$ is a finite 
$p$-group or $p'$-group, there is another $p$-local finite group 
$\SFL[_H]$ for each subgroup $H\le{}\Gamma$, such that $|\call_H|$ is 
homotopy equivalent to the covering space of $|\call|$ whose fundamental 
group is $\widehat{\theta}^{-1}(H)$.

Recall that for any $p$-local finite group $\SFL$ with associated 
quasicentric linking system $\callq$, 
$j\:S\rTo\pi_1(|\call|)\cong\pi_1(|\callq|)$ denotes the homomorphism 
induced by the distinguished monomorphism $\delta_S\:S\rTo\Aut_\call(S)$, 
and $J\:\callq\rTo\calb(\pi_1(|\callq|))$ is the functor which sends 
morphisms to loops as defined in Section 1.  Let $\widehat{\theta}$ be a 
homomorphism from $\pi_1(|\callq|)$ to a group $\Gamma$, and set 
	$$ \theta=\widehat{\theta}\circ{}j\in\Hom(S,\Gamma) 
	\qquad\textup{and}\qquad
	\widehat\Theta=\calb(\widehat{\theta})\circ{}J\:
	\callq \Right5{} \calb(\Gamma). $$
Note that these depend on a choice of a compatible set of inclusions 
$\{\iota_P^Q\}$ for $\callq$ (since $J$ depends on such a choice). For any 
subgroup $H\le\Gamma$, let $\llx{\widehat{\Theta}}H\subseteq\callq$ be the 
subcategory with the same objects and with morphism set 
$\widehat{\Theta}^{-1}(H)$; and let 
$\callq_H\subseteq\llx{\widehat{\Theta}}H$ be the full subcategory 
obtained by restricting to subgroups of $S_H\defeq\theta^{-1}(H)$. 
Finally, let $\lf{\widehat{\Theta}}H$ be the fusion system over $S_H$ 
generated by $\pi(\callq_H)\subseteq\calfq$ and restrictions of morphisms, 
and let $\call_H\subseteq\callq_H$ be the full subcategory on those 
objects which are $\calf_H$-centric.  

\begin{Thm}  \label{L:p-p'-quot.}
Let $\SFL$ be a $p$-local finite group, let $\callq$ be its associated 
quasicentric linking system, and let $\pi\:\callq\rTo\calf$ be the 
projection.  Assume a compatible set of inclusions $\{\iota_P^Q\}$ has
been chosen for $\callq$.  Fix a finite group $\Gamma$ which is a
$p$-group or a $p'$-group, and a surjective homomorphism
	$$ \widehat{\theta}\: \pi_1(|\callq|) \Onto5{} \Gamma. $$
Set $\theta=\widehat{\theta}\circ{}j\:S\Right1{}\Gamma$.  Fix 
$H\le{}\Gamma$, and set $S_H=\theta^{-1}(H)$.  Then $\SFL[_H]$ is also a 
$p$-local finite group, and (via the inclusion of $\call_H$ into $\callq$) 
$|\call_H|$ is homotopy equivalent to the covering space of 
$|\callq|\simeq|\call|$ with fundamental group $\widehat{\theta}^{-1}(H)$. 
\end{Thm}

\begin{proof}  Define $\Theta\:\Mor(\calf^q)\rTo\sset(\Gamma)$ by setting 
$\Theta(\alpha)=\widehat{\Theta}(\pi^{-1}(\alpha))$, where 
$\widehat{\Theta}=\calb(\widehat{\theta})\circ{}J$ as above.  Then 
$\Theta$ and $\theta$ satisfy hypotheses (i)--(iv) of Definition 
\ref{D:J,j}: points (i) and (ii) follow from (A)$_q$, while (iii) follows 
from Proposition \ref{deltaPQ} and (iv) from (C)$_q$.  Thus 
$(\Gamma,\theta,\Theta)$ is a fusion mapping triple on $\calfq$, and 
$\calf_H$ is a saturated fusion system over $S_H$ by Proposition 
\ref{F:p-solv.quot.}.  

Let $\opf$ and $\oppf$ be the categories of Definition \ref{D:oppf:opf}, 
and let $\oppf^q\subseteq\oppf$ and $\opf^q\subseteq\opf$ be the full 
subcategories whose objects are the $\calf$-quasicentric subgroups of $S$. 
By Proposition \ref{F:p-solv.quot.}(b,c), $\pi(\call_1^\bullet)$ contains 
$\oppf^q$ (if $\Gamma$ is a $p'$-group) or $\opf^q$ (if $\Gamma$ is a 
$p$-group).  Hence in either case, by Lemma \ref{F=F'.N_F(S)}(b), all 
morphisms in $\callq$ are composites of morphisms in $\call_1^\bullet$ and 
restrictions of morphisms in $\Aut_\call(S)$.  Since 
$\call_1^\bullet=\widehat{\Theta}^{-1}(1)$, by definition, and since 
$\widehat{\Theta}(\alpha)=\widehat{\Theta}(\beta)$ whenever $\alpha$ is a 
restriction of $\beta$, this shows that $\widehat{\Theta}$  restricts to a 
surjection of $\Aut_\call(S)$ onto $\Gamma$.  In particular, this implies 
that 
	\beq \textup{$\forall$ $P\in\Ob(\callq)$ and $\forall$ 
	$g\in{}\Gamma$, \ 
	$\exists$ $P'\le{}S$ and $\alpha\in\Iso_{\callq}(P,P')$ such that 
	$\widehat{\Theta}(\alpha)=g$;} \tag{1} \eeq
where in fact, $\alpha$ can always be chosen to be the restriction of an 
automorphism of $S$.  

We start by proving that $\callq_H$ is a linking system associated to 
$\calf_H$ (for its set of objects), and hence that $\call_H$ is a centric 
linking system.  Let $P\le{}S_H$ be a $\calf$-quasicentric subgroup, and 
choose $g\in{}P{\cdot}C_S(P)$.  By construction, 
$\widehat{\Theta}(\delta_S(g))=\theta(g)$, and 
$\widehat{\Theta}(\iota_P)=1$.  In particular, the inclusion morphisms are 
in $\call_H$. Also, $\iota_P\circ\delta_P(g)=\delta_S(g)\circ\iota_P$ by 
definition of an inclusion morphism (Definition \ref{D:incl}).  Hence 
$\widehat{\Theta}(\delta_P(g))=\theta(g)$ in this situation; and in 
particular $\delta_P(g)\in\Aut_{\callq_H}(P)$ if and only if $g\in{}S_H$. 
Thus $\delta_P$ restricts to a distinguished monomorphism 
$P{\cdot}C_{S_H}(P)\Right1{}\Aut_{\callq_H}(P)$ for $\callq_H$ and axiom 
(D)$_q$ is satisfied.  Moreover, if $f,f'\in\Mor_{\callq_H}(P,Q)$ are such 
that $\pi(f)=\pi(f')$ in $\homf(P,Q)$, and $g\in{}C_S(P)$ is the unique 
element such that $f'=f\circ\delta_P(g)$ (using axiom (A)$_q$ for 
$\callq$), then $\widehat{\Theta}(\delta_P(g))=1$, and hence $g\in{}S_H$.  
This shows that axiom (A)$_q$ for $\callq_H$ holds.  Axioms (B)$_q$ and 
(C)$_q$ for $\callq_H$ also follow immediately from the same properties 
for $\callq$.  

We have now shown that $\SFL[_H]$ is a $p$-local finite group.  It remains 
to show that $|\call_H|$ is homotopy equivalent to a certain covering space 
of $|\callq|\simeq|\call|$.  We show this by first choosing certain full 
subcategories $\call^x\subseteq\callq$ and $\call^x_H\subseteq\callq_H$ 
such that $|\call^x|\simeq|\callq|\simeq|\call|$ and 
$|\call^x_H|\simeq|\callq_H|\simeq|\call_H|$, and then proving directly 
that $|\call^x_H|$ is a covering space of $|\call^x|$.

If $\Gamma$ is a $p'$-group, then for any $P\le{}S_H=S$, $P$ is 
$\calf_H$-centric if and only if it is $\calf$-centric (Proposition 
\ref{F:p-solv.quot.}(c)).  By the above remarks, $\call_H$ is a centric 
linking system associated to $\calf_H$.  Set $\call^x=\call$ and 
$\call^x_H=\call_H$ in this case.

If $\Gamma$ is a $p$-group, then for any $P\le{}S_H$, $P$ is 
$\calf_H$-quasicentric if and only if it is $\calf$-quasicentric 
(Proposition \ref{F:p-solv.quot.}(b)).  So by what was just shown, 
$\callq_H$ is a quasicentric linking system associated to $\calf_H$ which 
extends $\call_H$.  Let $\call^x\subseteq\callq$ be the full subcategory 
whose objects are those subgroups $P\le{}S$ such that $P\cap{}S_H$ is 
$\calf$-quasicentric.  Set $S_1=\Ker(\theta)$, and let 
$\calf_1\subseteq\calf$ be the saturated fusion system over $S_1$ defined 
in Proposition \ref{F:p-solv.quot.}.  The definition of $\theta$ as a 
restriction of $\widehat{\Theta}$ ensures that $\theta(g)$ and 
$\theta(g')$ are $\Gamma$-conjugate whenever $g,g'$ are $\calf$-conjugate; 
in particular, no element of $S_1$ is $\calf$-conjugate to any element of 
$S{\sminus}S_1$.  Hence by Lemma \ref{P0:F0-centric}, for each $P\le{}S$ 
which is $\calf$-centric and $\calf$-radical, $P\cap{}S_1$ is 
$\calf_1$-centric, hence $\calf$-quasicentric, so $P\cap{}S_H$ is 
$\calf$-quasicentric, and thus $P\in\Ob(\call^x)$.  So the inclusion of 
$|\call^x|$ in $|\callq|$ is a homotopy equivalence by Proposition 
\ref{L-props}.  

Still assuming $\Gamma$ is a $p$-group, let $\call^x_H\subseteq\call^x$ be 
the subcategory with the same objects, where 
$\Mor(\call^x_H)=\widehat{\Theta}^{-1}(H)$.  For each 
$P\in\Ob(\call^x_H)=\Ob(\call^x)$, $P\cap{}S_H$ is $\calf$-quasicentric by 
assumption, hence $\calf_H$-quasicentric; and by Proposition \ref{deltaPQ}, 
each $\varphi\in\Mor_{\call^x_H}(P,Q)$ restricts to a unique morphism 
$\varphi_H\in\Mor_{\call^x_H}(P\cap{}S_H,Q\cap{}S_H)$.  These restrictions 
define a deformation retraction from $|\call^x_H|$ to $|\callq_H|$, and 
thus the inclusion of categories induces a homotopy equivalence 
$|\call^x_H|\simeq|\callq_H|\simeq|\call_H|$.  

Thus in both cases, we have chosen categories $\call^x_H\subseteq\call^x$ 
with the same objects, where $\call^x$ is a full subcategory of $\callq$ 
and $\Mor(\call^x_H)=\Mor(\call^x)\cap\widehat{\Theta}^{-1}(H)$, and where 
$|\call^x|\simeq|\call|$ and $|\call^x_H|\simeq|\call_H|$.  Let 
$\cale_\Gamma(\Gamma/H)$ be the category with object set $\Gamma/H$, and 
with a morphism $\b{g}$ from $aH$ to $gaH$ for each $g\in{}\Gamma$ and 
$aH\in{}\Gamma/H$.  Thus 
$\Aut_{\cale_\Gamma(\Gamma/H)}(1{\cdot}H)\cong{}H$, and 
$|\cale_\Gamma(\Gamma/H)|=EG/H\simeq{}BH$.  

Let $\widetilde{\call}$ be the pullback category in the following square:
	\begin{diagram}[w=40pt]
	\widetilde\call & \rTo & \cale_\Gamma(\Gamma/H) \\
	\dTo && \dTo \\
	\call^x & \rTo & \calb(\Gamma) \rlap{\,,}
	\end{diagram}
Thus $\Ob(\widetilde{\call})=\Ob(\call^x)\times{}\Gamma/H$, and 
$\Mor(\widetilde{\call})$ is the set of pairs of morphisms in $\call^x$ and 
$\cale_\Gamma(\Gamma/H)$ which get sent to the same morphism in 
$\calb(\Gamma)$.  Then $\call^x_H$ 
can be identified with the full subcategory of $\widetilde{\call}$ with 
objects the pairs $(P,1{\cdot}H)$ for $P\in\Ob(\call^x)$.  By (1), each 
object in $\widetilde{\call}$ is isomorphic to an object in $\call^x_H$, 
and so $|\call^x_H|\simeq|\widetilde{\call}|$.  By construction, 
$|\widetilde{\call}|$ is the covering space over $|\call^x|$ with 
fundamental group $\theta^{-1}(H)$.  Since 
$|\widetilde{\call}|\simeq|\call_H|$ and $|\call^x|\simeq|\call|$, this 
finishes the proof of the last statement. 
\end{proof}

The following is an immediate corollary to Theorem \ref{L:p-p'-quot.}.  

\begin{Cor}  \label{L:p-solv.quot.}
For any $p$-local finite group $\SFL$, any finite $p$-solvable 
group $\Gamma$, and any homomorphism
	$$ \widehat{\theta}\: \pi_1(|\call|) \Right5{} \Gamma, $$
there is a $p$-local finite group $\SFL[_0]$ such that $|\call_0|$ is 
homotopy equivalent to the covering space of $|\call|$ with fundamental 
group $\Ker(\widehat{\theta})$.  Furthermore, this can be chosen such that 
$S_0$ is a subgroup of $S$ and $\calf_0$ is a subcategory of $\calf$.
\end{Cor}


For any $p$-local finite group $\SFL$, since $\calf$ is a finite category, 
Corollary \ref{L:p-solv.quot.} implies that there is a unique maximal 
$p$-solvable quotient group of $\pi_1(|\call|)$, which is finite.  In 
contrast, if we look at arbitrary finite quotient groups of the 
fundamental group, they can be arbitrarily large.

As one example, consider the case where $p=2$, $S\in\Syl_2(A_6)$ (so 
$S\cong{}D_8$), $\calf=\calf_S(A_6)$, and $\call=\call^c_S(A_6)$.  It is 
not hard to show directly, using Van Kampen's theorem, that 
$\pi_1(|\call|)\cong\Sigma_4\underset{S}{\textup{\Large{$*$}}}\Sigma_4$:  
the amalgamated free product of two copies of $\Sigma_4$ intersecting in 
$S$, where each of the two subgroups $C_2^2$ in $D_8$ is normalized by one 
of the $\Sigma_4$.  Thus $\pi_1(|\call|)$ surjects onto any finite group 
$\Gamma$ which is generated by two copies of $\Sigma_4$ intersecting in 
the same way.  This is the case when $\Gamma=A_6$, and also when 
$\Gamma=PSL_2(q)$ for any $q\equiv\pm9$ (mod $16$).  However, the kernel 
of any such homomorphism defined on $\pi_1(|\call|)$ is torsion free (and 
infinite), and hence cannot be the fundamental group of the geometric 
realization of any centric linking system.  In fact, in this case, there 
is no nontrivial homomorphism from $\pi_1(|\call|)$ to a finite 
$2$-solvable group.


\bigskip

\newsect{Fusion subsystems and extensions of $p$-power index} 

Recall that for any saturated fusion system $\calf$ over a $p$-group $S$, 
we defined $\hfocal{\calf}{S}\nsg S$ to be the subgroup generated by all 
elements of the form $g^{-1}\alpha(g)$, for $g\in{}P\le{}S$, and 
$\alpha\in\autf(P)$ of order prime to $p$. In this section, we classify 
all saturated fusion subsystems of $p$-power index in a given saturated 
fusion system $\calf$ over $S$, and show that there is a one-to-one 
correspondence between such subsystems and the subgroups of $S$ which 
contain $\hfocal{\calf}{S}$.  In particular, there is a unique minimal 
subsystem $O^{p}(\calf)$ of this type, which is a fusion system over 
$\hfocal{\calf}{S}$. We then look at extensions, and describe the 
procedure for finding all larger saturated fusion systems of which $\calf$ 
is a fusion subsystem of $p$-power index.

\bigskip

\newsub{Subsystems of $p$-power index}
In this subsection, we classify all saturated fusion subsystems of 
$p$-power index in a given saturated fusion system $\calf$, and show that 
there is a unique minimal subsystem $O^{p}(\calf)$ of this type. We also 
show that there is a bijective correspondence between subgroups $T$ of the 
finite $p$-group
	$$ \gpf \defeq S/\hfocal{\calf}S, $$
and fusion subsystems $\calf_T$ of $p$-power index in $\calf$.  We have 
already seen that $\gpf=S/\hfocal{\calf}S$ is isomorphic 
to $\pi_1(|\call|\pcom)$ for any centric linking system $\call$ associated 
to $\calf$, and our result also shows that all (connected) covering 
spaces of $|\call|\pcom$ are realized as classifying spaces of subsystems 
of $p$-power index.  The index of $\calf_T$ in $\calf$ can then be 
defined to be the index of $T$ in $\gp{\calf}$, or equivalently the 
covering degree of the covering space.

The main step in doing this is to construct a fusion mapping triple 
$(\gpf,\theta,\Theta)$ for $\calfq$, where $\theta$ is the canonical 
surjection of $S$ onto $\gpf$.  This construction parallels very closely 
the construction in Proposition \ref{ext-function} of a functor from 
$\callq$ to $\calb(\gpf)$, when $\call$ is a linking system associated to 
$\calf$.  In fact, we could in principal state and prove the two results 
simultaneously, but the extra terminology which that would require seemed 
to add more complications than would be saved by combining the two.

The following lemma provides a very general, inductive tool for 
constructing explicit fusion mapping triples.  

\begin{Lem} \label{ext.F}
Fix a saturated fusion system $\calf$ over a $p$-group $S$.  Let $\calh_0$ 
be a set of $\calf$-quasicentric subgroups of $S$ which is closed under 
$\calf$-conjugacy and overgroups.  Let $\calp$ be an $\calf$-conjugacy 
class of $\calf$-quasicentric subgroups maximal among those not in 
$\calh_0$, set $\calh=\calh_0\cup\calp$, and let 
$\calf^{\calh_0}\subseteq\calf^{\calh}\subseteq\calfq$ be the full 
subcategories with these objects.  Fix a group $\Gamma$ and a homomorphism 
$\theta\:S\rTo\Gamma$, and let 
	$$ \Theta_0 \: \Mor(\calf^{\calh_0}) \Right4{} \sset(\Gamma) $$
be such that $(\Gamma,\theta,\Theta_0)$ is a fusion mapping triple for  
$\calf^{\calh_0}$. Let $P\in\calp$ be fully normalized in $\calf$, and fix 
a homomorphism 
	$$ \Theta_P\: \autf(P) \Right5{} 
	N_\Gamma(\theta(C_S(P)))/\theta(C_S(P)) $$
such that the following two conditions hold:
\begin{enumerate} 
\item[\textup{(\texttt{+})}]  $x\theta(g)x^{-1}=\theta(\alpha(g))$ for all 
$g\in{}P$, $\alpha\in\autf(P)$, and $x\in\Theta_P(\alpha)$. 
\item[\textup{($*$)}]  For all $P\lneqq{}Q\le{}S$ such that $P\nsg{}Q$ and 
$Q$ is fully normalized in $N_\calf(P)$, and for all $\alpha\in\autf(P)$ 
and $\beta\in\autf(Q)$ such that $\alpha=\beta|_P$, 
$\Theta_P(\alpha)\supseteq\Theta_0(\beta)$.
\end{enumerate}
Then there is a unique extension of $\Theta_0$ to a fusion mapping triple 
$(\Gamma,\theta,\Theta)$ for $\calf^{\calh}$ such that 
$\Theta(\alpha)=\Theta_P(\alpha)$ for all $\alpha\in\autf(P)$.  
\end{Lem}

\begin{proof}  Note that (\texttt{+}) is just point (v) of Lemma 
\ref{P:Theta-theta} applied to the subgroup $P$, while ($*$) is just point 
(vi) applied to restrictions to $P$.  So both of these conditions are 
necessary if we want to be able to extend $\Theta_0$ and $\Theta_P$ to a 
fusion mapping triple for $\calf^\calh$.

The uniqueness of the extension is an immediate consequence 
of Alperin's fusion theorem, in the form of Theorem \ref{Alp.fusion}(a).  
The proof of existence is almost identical to the proof of Lemma 
\ref{ext.L}, so we just sketch it here briefly.

We first show that we can replace ($*$) by the following (a priori 
stronger) statement:
\begin{enumerate} 
\item[\textup{($**$)}]  for all $Q,Q'\le{}S$ which strictly contain $P$, 
and for all $\beta\in\homf(Q,Q')$ and $\alpha\in\autf(P)$ such that 
$\alpha=\beta|_P$, $\Theta_P(\alpha)\supseteq\Theta_0(\beta)$.
\end{enumerate}
It suffices to show this when $P$ is normal in $Q$ and $Q'$, since 
otherwise we can replace $Q$ and $Q'$ by $N_Q(P)$ and $N_{Q'}(P)$.  In 
this case, $\beta\in\Hom_{N_\calf(P)}(Q,Q')$, and by Theorem 
\ref{Alp.fusion}(a) (Alperin's fusion theorem), it is a composite of 
restrictions of automorphisms of subgroups fully normalized in 
$N_\calf(P)$.  So it suffices to prove ($**$) when $\beta$ is such an 
automorphism, and this is what is assumed in ($*$).

Now fix any morphism $\varphi\in\homf(P_1,Q)$ which lies in 
$\calf^{\calh}$ but not in $\calf^{\calh_0}$; thus $P_1\in\calp$.  Set 
$P_2=\varphi(P_1)\le{}Q$, and let $\varphi'\in\isof(P_1,P_2)$ be the 
``restriction'' of $\varphi$.  By Lemma \ref{N->N}, there are isomorphisms 
$\widebar{\varphi}_i\in\Iso_{\calfq}(N_S(P_i),N_i)$, for some 
$N_i\le{}N_S(P)$ containing $P$, which restrict to isomorphisms 
$\varphi_i\in\Iso_{\calfq}(P_i,P)$.  Fix elements 
$x_i\in\Theta_0(\widebar{\varphi}_i)$.  Set 
$\psi=\varphi_2\circ\varphi'\circ\varphi_1^{-1}\in\Aut_{\calfq}(P)$.  Thus
$\varphi'=\varphi_2^{-1}\circ\psi\circ\varphi_1$, and we define
	$$ \Theta(\varphi) = \Theta(\varphi') = 
	x_2^{-1} {\cdot} \Theta_P(\psi) {\cdot} x_1. $$
This is independent of the choice of $x_i$, since 
$\Theta_0(\widebar{\varphi}_i)\subseteq\theta(C_S(P)){\cdot}x_i$ by axioms 
(i) and (ii) in the definition of a fusion mapping triple.  It is 
independent of the choice of $\varphi_1$ and $\varphi_2$ by the same 
argument as was used in the proof of Lemma \ref{ext.L} (and this is where 
we need point ($**$)).  Conditions (i)--(iv) are easily checked.  For 
example, (iv) --- the condition that $x\theta(g)x^{-1}=\theta(\alpha(g))$ 
whenever $g\in{}P_1$, $\varphi\in\homf(P_1,P_2)$, and 
$x\in\Theta(\varphi)$ --- holds when $\varphi$ can be extended to a larger 
subgroup since $(\Gamma,\theta,\Theta_0)$ is already a fusion mapping 
triple, holds for $\varphi\in\autf(P)$ by (\texttt{+}), and thus holds in 
the general case since $\Theta(\varphi)$ was defined via a composition of 
such morphisms.  Thus $(\Gamma,\theta,\Theta)$ is a fusion mapping triple 
on $\calf^{\calh}$.
\end{proof}

The construction of a fusion mapping triple to $\gpf$ in the following 
lemma is a first application of Lemma \ref{ext.F}.  Another application 
will be given in the next section. 

\begin{Lem} \label{Theta:p-quot.}
Let $\calf$ be a saturated fusion system over a $p$-group $S$, and let 
	$$ \theta\:S\Right5{} \gpf = S/\hfocal{\calf}S $$ 
be the projection.  Then there is a fusion mapping triple 
$(\gpf,\theta,\Theta)$ on $\calfq$. 
\end{Lem}

\begin{proof} The function $\Theta$ will be constructed inductively, using 
Lemma \ref{ext.F}.  Let $\calh_0\subseteq\Ob(\calfq)$ be a subset 
(possibly empty) which is closed under $\calf$-conjugacy and overgroups.  
Let $\calp$ be an $\calf$-conjugacy class of $\calf$-quasicentric 
subgroups maximal among those not in $\calh_0$, set 
$\calh=\calh_0\cup\calp$, and let 
$\calf^{\calh_0}\subseteq\calf^{\calh}\subseteq\calfq$ be the full 
subcategories with these objects.  Assume we have already constructed a 
fusion mapping triple $(\gpf,\theta,\Theta_0)$ for $\calf^{\calh_0}$.  

We recall the notation of Lemma \ref{G:hyperfocal}.  If $G$ is any finite 
group, and $S\in\sylp{G}$, then 
	$$ O_G^p(S) \defeq \bigl\langle [g,x] \,\big|\, g\in{}P\le{}S,\ 
	x\in{}N_G(P) \textup{ of order prime to $p$} \bigr\rangle. $$
By Lemma \ref{G:hyperfocal}, $O_G^p(S)=S\cap{}O^p(G)$, and hence 
$G/O^p(G)\cong{}S/O_G^p(S)$.  

Fix $P\in\calp$ which is fully normalized in $\calf$.  Let $N_0$ be the 
subgroup generated by commutators $[g,x]$ for $g\in{}N_S(P)$ and 
$x\in{}N_{\autf(P)}(N_S(P))$ of order prime to $p$.  Then 
$\Aut_S(P)\in\sylp{\autf(P)}$ and $\Aut_{N_0}(P)=O_{\autf(P)}^p(\Aut_S(P))$, 
and by Lemma \ref{G:hyperfocal},
	$$ \autf(P)\big/O^p(\autf(P)) \cong Aut_S(P)/Aut_{N_0}(P)\cong 
	N_S(P)\big/\gen{N_0,C_S(P)}. $$
Also, $N_0\le\hfocal{\calf}{S}$, and so the inclusion of $N_S(P)$ into $S$ 
induces a homomorphism 
	$$ \Theta_P\: \autf(P) \Onto4{} 
	\underset{\cong{}N_S(P)/\gen{N_0,C_S(P)}}
	{\autf(P)/O^p(\autf(P))} 
	\Right4{} 
	\underset{\cong{}N_S(C_S(P){\cdot}S_0)/(C_S(P){\cdot}S_0)}
	{N_{\gpf}(\theta(C_S(P)))/\theta(C_S(P))}. $$
Here, we write $S_0=\hfocal{\calf}S$ for short; thus $\gpf=S/S_0$ and 
$\theta(C_S(P))=C_S(P){\cdot}S_0/S_0$.  

Point (\texttt{+}) in Lemma \ref{ext.F} holds by the construction of 
$\Theta_P$.  So it remains only to prove that condition ($*$) in Lemma 
\ref{ext.F} holds.  

To see this, fix $P\lneqq{}Q\le{}S$ such that $P\nsg{}Q$ 
and $Q$ is fully normalized in $N_\calf(P)$, and fix 
$\alpha\in\Aut_{\calfq}(P)$ and $\beta\in\Aut_{\calfq}(Q)$ such that 
$\alpha=\beta|_P$.  We must show that 
$\Theta_P(\alpha)\supseteq\Theta_0(\beta)$. Upon replacing $\alpha$ by 
$\alpha^k$ and $\beta$ by $\beta^k$ for some appropriate $k\equiv1$ (mod 
$p$), we can assume that both automorphisms have order a power of $p$.  
Since $Q$ is fully normalized, $\Aut_{N_S(P)}(Q)$ is a Sylow subgroup of 
$\Aut_{N_{\calf}(P)}(Q)$; and hence there are automorphisms 
$\widebar{\gamma}\in\Aut_{\calf}(Q)$ and $\gamma\in\Aut_{\calf}(P)$ of 
order prime to $p$ such that $\gamma=\widebar{\gamma}|_P$ and 
$\widebar{\gamma}\beta\widebar{\gamma}^{-1}=c_g|_Q$ for some 
$g\in{}N_S(Q)\cap{}N_S(P)$.  Then $\gamma\alpha\gamma^{-1}=c_g|_P$.  Also, 
$1\in\Theta_P(\gamma)$ and $1\in\Theta_0(\widebar{\gamma})$, since both 
automorphisms have order prime to $p$.  So
    $$ \Theta_0(\beta)=\Theta(c_g|_Q) = g{\cdot}\theta(C_S(Q))
    \subseteq g{\cdot}\theta(C_S(P)) =
    \Theta_P(c_g|_P)=\Theta_P(\alpha). $$

Thus, by Lemma \ref{ext.L}, we can extend $\Theta_0$ to a fusion mapping 
triple on $\calf^{\calh}$.  Upon continuing this procedure, we obtain a 
fusion mapping triple defined on all of $\calfq$.
\end{proof}

We now apply Lemma \ref{Theta:p-quot.} to classify fusion subsystems of 
$p$-power index.  Recall that by definition, if $\calf$ is a saturated 
fusion system over $S$ and $\calf_0\subseteq\calf$ is a fusion subsystem 
over $S_0\le{}S$, then $\calf$ has $p$-power index if and only if 
$S_0\ge\hfocal{\calf}S$, and $\Aut_{\calf_0}(P)\ge{}O^p(\autf(P))$ for all 
$P\le{}S_0$.

\begin{Thm} \label{F:p-quot.}
Fix a saturated fusion system $\calf$ over a $p$-group $S$.  Then for 
each subgroup $T\le{}S$ containing $\hfocal{\calf}S$, there is 
a unique saturated fusion system $\calf_T\subseteq\calf$ over $T$ with 
$p$-power index.  For each such $T$, $\calf_T$ has the properties:
\begin{enumerate}  
\item a subgroup $P\le{}T$ is $\calf_T$-quasicentric if 
and only if it is $\calf$-quasicentric; and 
\item for each pair $P,Q\le{}T$ of $\calf$-quasicentric subgroups,
	$$ \Hom_{\calf_T}(P,Q) = \bigl\{\varphi\in\homf(P,Q)\,\big|\, 
	\Theta(\varphi)\cap(T/\hfocal{\calf}S)\ne\emptyset \bigr\}. $$
\end{enumerate}
Here, 
	$$ \Theta\: \Mor(\calfq) \Right5{} \sset(\gpf) =
	\sset(S/\hfocal{\calf}S) $$
is the map of Lemma \ref{Theta:p-quot.}.  
\end{Thm}

\begin{proof}  Let $\calf_T\subseteq\calf$ be the fusion system over $T$ 
defined on $\calf$-quasicentric subgroups by the formula in (b), and then 
extended to arbitrary subgroups by taking restrictions and composites.  
(This is the fusion system denoted $\calf_{T/\hfocal{\calf}S}$ in 
Proposition \ref{F:p-solv.quot.}, but we simplify the notation here.)  
By Proposition \ref{F:p-solv.quot.}(a,b) (applied with $\Gamma=\gpf$ 
and $H=T/\hfocal{\calf}S$), $\calf_T$ is saturated, a subgroup 
$P\le{}T$ is $\calf_T$-quasicentric if and only if it is 
$\calf$-quasicentric, and $\Aut_{\calf_T}(P)\ge{}O^p(\autf(P))$ for all 
$P\le{}T$.  

Now let $\calf'_T\subseteq\calf$ be another saturated subsystem over the 
same subgroup $T$ which also has $p$-power index.  We claim that 
$\calf'_T=\calf_T$, and thus that $\calf_T$ is the unique subsystem with 
these properties.  By assumption, for each $P\le{}T$, $\Aut_{\calf_T}(P)$ 
and $\Aut_{\calf'_T}(P)$ both contain $O^p(\autf(P))$, and hence each is 
generated by $O^p(\autf(P))$ and any one of its Sylow $p$-subgroups.  So
if $P$ is fully normalized in $\calf_T$ and $\calf'_T$ both, then 
	\beq \Aut_{\calf_T}(P) = \gen{O^p(\autf(P)),\Aut_{T}(P)} = 
	\Aut_{\calf'_T}(P). \tag{1} \eeq

In particular, $\Aut_{\calf_T}(T)=\Aut_{\calf'_T}(T)$.  Set 
$p^k=|T|$, fix $0\le{}m<k$, and assume inductively that 
$\Hom_{\calf_T}(P,Q)=\Hom_{\calf'_T}(P,Q)$ for all $P,Q\le{}T$ of order 
$>p^m$.  By Alperin's fusion theorem for saturated fusion systems (Theorem 
\ref{Alp.fusion}(a)), if $|P|=p^m$, $|Q|\ge{}p^m$, and $P\ne{}Q$, then 
all morphisms in $\Hom_{\calf_T}(P,Q)$ and $\Hom_{\calf'_T}(P,Q)$ are 
composites of restrictions of morphisms between subgroups of order $>p^m$, 
and hence $\Hom_{\calf_T}(P,Q)=\Hom_{\calf'_T}(P,Q)$ by the induction 
hypothesis.  In particular, two subgroups of order $p^m$ are 
$\calf_T$-conjugate if and only if they are $\calf'_T$-conjugate.  So for 
any $P\le{}T$ of order $p^m$, $P$ is fully normalized in $\calf_T$ if 
and only if it is fully normalized in $\calf'_T$.  In either case, 
$\Aut_{\calf_T}(P)=\Aut_{\calf'_T}(P)$:  by (1) if $P$ is fully 
normalized, and by Alperin's fusion theorem again (and the induction 
hypothesis) if it is not.  
\end{proof}

We can now define $\Opf$ as the \emph{minimal} fusion subsystem of $\calf$ 
of $p$-power index:  the unique fusion subsystem over $\hfocal{\calf}S$ of 
$p$-power index.  The next theorem will show that when $\calf$ has an 
associated linking system $\call$, then $\Opf$ has an associated linking 
system $\Opl$, and that $|\Opl|\pcom$ is the universal cover of 
$|\call|\pcom$.  

\begin{Thm} \label{L:p-quot.}
Fix a $p$-local finite group $\SFL$.  Then for each subgroup $T\le{}S$ 
containing $\hfocal{\calf}S$, there is a unique $p$-local finite subgroup 
$(T,\calf_T,\call_T)$ such that $\calf_T$ has $p$-power index in $\calf$, 
and such that $\callq_T=\pi^{-1}(\calfq_T)$ where $\pi$ is the usual 
projection of $\callq$ onto $\calfq$.  Furthermore, $|\call_T|$ is 
homotopy equivalent, via the inclusion of $|\callq_T|\simeq|\call_T|$ into 
$|\callq|\simeq|\call|$, to a covering space of $|\call|$ of degree 
$[S:T]$.  Hence the classifying space $|\call_T|\pcom$ of 
$(T,\calf_T,\call_T)$ is homotopy equivalent to the covering space of 
$|\call|\pcom$ with fundamental group $T/\hfocal{\calf}S$.
\end{Thm}

\begin{proof}  Assume a compatible set of inclusions $\{\iota_P^Q\}$ has 
been chosen for $\callq$.  By Proposition \ref{ext-function}, there is a 
functor $\lambda\:\callq\rTo\calb(\gpf)$ which sends inclusions to the 
identity, and such that $\lambda(\delta_S(g))=g$ for all $g\in{}S$.  Hence 
by Theorem \ref{L:p-p'-quot.}, $(T,\calf_T,\call_T)$ is a $p$-local finite 
group, and $|\call_T|$ is a covering space of $|\call|$.  Also, if we 
write $\call_1$ for the linking system over $\hfocal{\calf}S$, then the 
fibration sequences
	$$ |\call_1|\rTo|\call|\rTo B\gpf
	\qquad\textup{and}\qquad
	|\call_1|\rTo|\call_T|\rTo B(T/\hfocal{\calf}S) $$
are still fibration sequences after $p$-completion \cite[II.5.2(iv)]{BK}, 
and hence $|\call_T|\pcom$ is the covering space of $|\call|\pcom$ with 
fundamental group $T/\hfocal{\calf}S$. The uniqueness follows from Theorem 
\ref{F:p-quot.}.
\end{proof}

Thus there is a bijective correspondence between fusion subsystems of 
$(S,\calf)$, or $p$-local finite subgroups of $\SFL$, of $p$-power index, 
and subgroups of $S/\hfocal{\calf}S\cong\pi_1(|\call|\pcom)$.  The 
classifying spaces of the $p$-local finite subgroups of $\SFL$ of 
$p$-power index are (up to homotopy) just the covering spaces of the 
classifying space of $\SFL$.

\bigskip

\newcommand{\calm}{\mathcal{M}}

\newsub{Extensions of $p$-power index}
We next consider the opposite problem:  how to construct extensions 
of $p$-power index of a given $p$-local finite group.  In the course of 
this construction, we will see that the linking system really is needed to 
construct an extension of the fusion system.  The following definition will 
be useful.  

\begin{Defi} \label{D:autfus}
Fix a saturated fusion system $\calf$ over a $p$-group $S$.  An 
automorphism $\alpha\in\Aut(S)$ is \emph{fusion preserving} if it 
normalizes $\calf$; i.e., if it induces an automorphism of the category 
$\calf$ by sending $P$ to $\alpha(P)$ and $\varphi\in\Mor(\calf)$ to 
$\alpha\varphi\alpha^{-1}\in\Mor(\calf)$.  Let $\Aut\fus(S,\calf)\le\Aut(S)$ 
denote the group of all fusion preserving automorphisms, and set 
	$$ \Out\fus(S,\calf)=\Aut\fus(S,\calf)/\Aut_\calf(S). $$
\end{Defi}

We first describe the algebraic data needed to determine extensions of 
$p$-power index.  Fix a $p$-local finite group $\SFL$, let $\callq$ be the 
associated quasicentric linking system, and let $\{\iota_P^Q\}$ be a 
compatible set of inclusions.  Then for any $g\in{}S$, $g$ acts on the set
$\Mor(\callq)$ by composing on the left or right with $\delta_S(g)$ and 
its restrictions.  Thus for any $\varphi\in\Mor_{\callq}(P,Q)$, we set 
	$$ g\varphi=\delta_{Q,gQg^{-1}}(g)\circ\varphi
	\in\Mor_{\callq}(P,gQg^{-1}) $$
and
	$$ \varphi{}g=\varphi\circ\delta_{g^{-1}Pg,P}(g)\in
	\Mor_{\callq}(g^{-1}Pg,Q). $$
This defines natural left and right actions of $S$ on the set 
$\Mor(\callq)$.  The resulting conjugation action 
$\varphi\mapsto{}g\varphi{}g^{-1}$ extends to an action on the category, 
where $g$ sends an object $P$ to $gPg^{-1}$.  The functor 
$\pi\:\callq\rTo\calfq$ is equivariant with respect to the conjugation 
action of $S$ on $\callq$ and the action of $\Inn(S)\le\Aut\fus(S,\calf)$ 
on $\calf$.

If $\SFL[_0]$ is contained in $\SFL$ with $p$-power index, and 
$S_0\nsg{}S$, then the $S$ action on $\call$ clearly restricts to an 
$S$-action on $\call_0$.  The following theorem provides a converse to 
this.  Given a $p$-local finite group $\SFL[_0]$, an extension $S$ of 
$S_0$, and an $S$-action on $\call$ which satisfies certain obvious 
compatibility conditions, this data always determines a $p$-local finite 
group which contains $\SFL[_0]$ with $p$-power index.  

\begin{Thm} \label{SFL_0<SFL}
Fix a $p$-local finite group $\SFL[_0]$, and assume that a compatible set of 
inclusions $\{\iota_P^Q\}$ has been chosen for $\call_0$.  Fix a $p$-group 
$S$ such that $S_0\nsg{}S$ and $\Aut_S(S_0)\le\Aut\fus(S_0,\calf_0)$, and 
an action of $S$ on $\call_0$ which: 
\begin{enumerate} 
\item extends the conjugation action of $S_0$ on $\call_0$; 
\item makes the canonical monomorphism 
$\delta_{S_0}\:S_0\Right2{}\Aut_{\call_0}(S_0)$ $S$-equivariant; 

\item makes the projection $\pi\:\call_0\Right2{}\calf_0$ $S$-equivariant 
with respect to the $\Aut_S(S_0)$-action on $\calf_0$; and 

\item sends inclusion morphisms in $\call_0$ to inclusion morphisms. 
\end{enumerate}
Then there is a $p$-local finite group $\SFL$ such that 
$\calf\supseteq\calf_0$, $\callq\supseteq\call_0$, the conjugation action 
of $S$ on $\callq$ restricts to the given $S$-action on $\call_0$, and 
$\SFL[_0]$ is a subgroup of $p$-power index in $\SFL$.
\end{Thm}

\begin{proof}  Set
	$$ \calh_0 = \Ob(\call_0)=\{P\le{}S_0\,|\,\textup{$P$ is 
	$\calf_0$-centric} \}
	\qquad\textup{and}\qquad
	\calh = \{P\le{}S \,|\, P\cap{}S_0\in\calh_0 \} . $$
To simplify notation, for any $P\le{}S$, we write $P_0=P\cap{}S_0$.  For 
$g\in{}S$ and $\varphi\in\Mor_{\call_0}(P,Q)$, we write 
$g\varphi{}g^{-1}\in\Mor_{\call_0}(gPg^{-1},gQg^{-1})$ for the given action 
of $g$ on $\varphi$.  By (a), when $g\in{}S_0$, this agrees with the 
morphism $g\varphi{}g^{-1}$ already defined.

\smallskip

\newcommand\5[1]{[\![#1]\!]}

\noindent\textbf{Step 1: }  We first define categories 
$\call_1\supseteq\call_0$ and $\calf_1\supseteq\calf_0$, where 
$\Ob(\calf_1)=\Ob(\calf_0)$ and $\Ob(\call_1)=\calh_0$.  Set
	$$ \Mor(\call_1) = S \times_{S_0} \Mor(\call_0) 
	= \bigl(S \times \Mor(\call_0)\bigr)\big/{\sim}, $$
where $(gg_0,\varphi)\sim(g,g_0\varphi)$ for $g\in{}S$, $g_0\in{}S_0$, and 
$\varphi\in\Mor(\call_0)$.  If $\varphi\in\Mor_{\call_0}(P,Q)$, then 
$\5{g,\varphi}\in\Mor_{\call_1}(P,gQg^{-1})$ denotes the equivalence class 
of the pair $(g,\varphi)$.  Composition is defined by
	$$ \5{g,\varphi}\circ\5{h,\psi}=\5{gh,h^{-1}\varphi{}h\circ\psi}. $$
Note that if $\varphi\in\Mor_{\call_0}(P,Q)$, then 
$h^{-1}\varphi{}h\in\Mor_{\call_0}(h^{-1}Ph,h^{-1}Qh)$.  To show that this 
is well defined, we note that for all $g,h\in{}S$, $g_0,h_0\in{}S_0$, and 
$\varphi,\psi\in\Mor(\call_0)$ with appropriate domain and range,
	\begin{multline*}  
	\5{gg_0,\varphi}\circ\5{hh_0,\psi}
	=\5{gg_0hh_0,h_0^{-1}(h^{-1}\varphi{}h)h_0\circ\psi}
	=\5{gh{\cdot}(h^{-1}g_0h),(h^{-1}\varphi{}h)\circ{}h_0\psi} \\
	=\5{gh,(h^{-1}g_0h){\cdot}(h^{-1}\varphi{}h)\circ{}h_0\psi} 
	=\5{gh,h^{-1}(g_0\varphi)h\circ{}h_0\psi}
	=\5{g,g_0\varphi}\circ\5{h,h_0\psi}\,.
	\end{multline*}
Here, the second to last equality follows from assumptions (b) and (d).  

Let $\calf_1$ be the smallest fusion system over $S_0$ which contains 
$\calf_0$ and $\Aut_S(S_0)$.  By assumption, 
$\Aut_S(S_0)\le\Aut\fus(S_0,\calf_0)$.  Thus for each $g\in{}S$, $c_g$ 
normalizes the fusion system $\calf_0$:  for each 
$\varphi\in\Mor(\calf_0)$ there is $\varphi'\in\Mor(\calf_0)$ such that 
$\varphi\circ{}c_g=c_g\circ\varphi'$.  Hence each morphism in $\calf_1$ 
has the form $c_g\circ\varphi$ for some $g\in{}S$ and 
$\varphi\in\Aut(\calf_0)$.  Define 
	$$ \pi_{\call_1} \: \call_1 \Right4{} \calf_1 $$
by sending $\pi_{\call_1}(\5{g,\varphi})=c_g\circ\pi_0(\varphi)$, where 
$\pi_0$ denotes the natural projection from $\call_0$ to $\calf_0$.  This 
is a functor by (c).

For all $P,Q\in\calh_0$, define
	$$ \widehat{\delta}_{P,Q} \: N_S(P,Q) \Right5{} 
	\Mor_{\call_1}(P,Q) $$
by setting $\widehat{\delta}_{P,Q}(g)=\5{g,\iota_{P}^{g^{-1}Qg}}$.  
This extends the canonical monomorphism $\delta_{P,Q}$ defined from 
$N_{S_0}(P,Q)$ to $\Mor_{\call_0}(P,Q)$.  To simplify the notation below, 
we sometimes write $\q{x}=\widehat{\delta}_{P,Q}(x)$ for $x\in{}N_S(P,Q)$.

\smallskip

\noindent\textbf{Step 2: }  We next construct categories $\call_2$ and 
$\calf_2$, both of which have object sets $\calh$, and which contain 
$\call_1$ and the restriction of $\calf_1$ to $\calh_0$, respectively.  
Afterwards, we let $\calf$ be the fusion system over $S$ generated by 
$\calf_2$ and restrictions of morphisms.  

Before doing this, we need to know that the following holds for each 
$P,Q\in\calh_0$ and each $\psi\in\Mor_{\call_1}(P,Q)$:
	\beq \forall\, x\in{}N_S(P) \textup{ there is at most one }
	y\in{}N_S(Q) \textup{ such that } 
	\q{y}\circ\psi=\psi\circ\q{x}. \tag{1} \eeq
Since $\psi$ is the composite of an isomorphism and an inclusion (and the 
claim clearly holds if $\psi$ is an isomorphism), it suffices to prove 
this when $P\le{}Q$ and $\psi$ is the inclusion.  We will show that we 
must have $y=x$ in that case.  By definition, 
	$$ \iota_P^Q\circ\q{x} = \5{1,\iota_P^Q}\circ\5{x,\Id_P} = 
	\5{x,\iota_P^{x^{-1}Qx}}
	\quad\textup{and}\quad
	\q{y}\circ\iota_P^Q = \5{y,\iota_P^Q} $$
(since conjugation sends inclusions to inclusions).  So if 
$\q{y}\circ\iota_P^Q=\iota_P^Q\circ\q{x}$, then there exists $g_0\in S_0$ 
such that $y=xg_0$ and $\iota_P^Q=g_0^{-1}\iota_P^{x^{-1}Qx}$, and thus 
	$$ \delta_{P,Q}(1)=\iota_P^Q
	=g_0^{-1}\iota_P^{x^{-1}Qx} = 
	\delta_{x^{-1}Qx,Q}(g_0^{-1}) \circ \delta_{P,x^{-1}Qx}(1) =
	\delta_{P,Q}(g_0^{-1}). $$ 
Hence $g_0=1$ by the injectivity of $\delta_{P,Q}$ in Proposition 
\ref{deltaPQ}; and so $\q{y}\circ\iota_P^Q=\iota_P^Q\circ\q{x}$ only if 
$x=y\in{}N_S(Q)$.  

Now let $\call_2$ be the category with $\Ob(\call_2)=\calh$, and where 
for all $P,Q\in\calh$,
	$$ \Mor_{\call_2}(P,Q) = 
	\bigl\{ \psi\in\Mor_{\call_1}(P_0,Q_0) \,\big|\,
	\forall\, x\in{}P\ \exists\, y\in{}Q \textup{ such that }
	\psi\circ\q{x}=\q{y}\circ\psi \bigr\}. $$
Let 
	$$ \widehat{\delta}_{P,Q} \: 
	\underset{\subseteq{}N_S(P_0,Q_0)}{N_S(P,Q)} \Right5{} 
	\underset{\subseteq\Mor_{\call_1}(P_0,Q_0)}{\Mor_{\call_2}(P,Q)} $$
be the restriction of $\widehat{\delta}_{P_0,Q_0}$.  Let $\calf_2$ be the 
category with $\Ob(\calf_2)=\calh$, and where
\begin{small}  
	$$ \Mor_{\calf_2}(P,Q) = 
	\bigl\{ \varphi\in\Hom(P,Q) \,\big|\, 
	\exists\,\psi\in\Mor_{\call_2}(P,Q) \textup{ such that }
	\psi\circ\q{x}=\q{\varphi(x)}\circ\psi\ \forall\, x\in{}P
	\bigr\}. $$
\end{small}
Let $\pi\:\call_2\rTo\calf_2$ be the functor which sends 
$\psi\in\Mor_{\call_2}(P,Q)$ to the homomorphism $\pi(\psi)(x)=y$ if 
$\psi\circ\q{x}=\q{y}\circ\psi$ (uniquely defined by (1)).  Let $\calf$ be 
the fusion system over $S$ generated by $\calf_2$ and restriction of 
homomorphisms.  

Set $\Gamma=S/S_0$.  Let
	$$ \widehat\theta\:\call_2\Right1{}\calb(\Gamma) $$
be the functor defined by setting $\widehat{\theta}(\5{g,\varphi})=gS_0$.  In 
particular, $\bigl(\theta^{-1}(1)\bigr)|_{\calh_0}=\call_0$.  

\smallskip

\noindent\textbf{Step 3: } 
We next show that each $P\in\calh$ is $\calf$-conjugate to a subgroup $P'$ 
such that $P'_0$ is fully normalized in $\calf_0$.  Moreover, we show that 
$P'$ can be chosen so that the following holds:
	\beq \forall\, g\in{}S \textup{ such that } gP_0g^{-1} 
	\textup{ is $\calf_0$-conjugate to $P_0$, \ 
	$g{\cdot}S_0\cap{}N_S(P'_0)\ne\emptyset$.} \tag{2} \eeq


To see this, let $\calp_{\textup{fn}}$ be the set of all $S_0$-conjugacy 
classes $[P'_0]$ of subgroups $P'_0\le{}S_0$ which are $\calf_0$-conjugate 
to $P_0$ and fully normalized in $\calf_0$.  (If $P'_0$ is fully 
normalized in $\calf_0$, then so is every subgroup in $[P'_0]$.)  Let 
$S'\subseteq{}S$ be the subset of elements $g\in{}S$ such that 
$gP_0g^{-1}$ is $\calf_0$-conjugate to $P_0$.  In particular, 
$S'\ge{}N_S(P_0)\ge{}P$.  Since each $g\in{}S$ acts 
on $\calf_0$ --- two subgroups $Q,Q'\le{}S_0$ are $\calf_0$-conjugate if 
and only if $gQg^{-1}$ and $gQ'g^{-1}$ are $\calf_0$-conjugate --- $S'$ is 
a subgroup of $S$.  


For all $g\in{}S'$ and $[P'_0]\in\calp_{\textup{fn}}$, $gP'_0g^{-1}$ is 
$\calf_0$-conjugate to $gP_0g^{-1}$ and hence to $P_0$, and is fully 
normalized since $g$ normalizes $S_0$.  Thus $S'/S_0$ acts on 
$\calp_{\textup{fn}}$, and this set has order prime to $p$ by Proposition 
\ref{|Rep(P,S)|}.  So we can thus choose a subgroup $P'_0\le{}S_0$ such 
that $[P'_0]\in\calp_{\textup{fn}}$ and is fixed by $S'$.  In particular, 
$P'_0$ is fully normalized in $\calf_0$.  Also, for each $g\in{}S'$, 
some element of $g{\cdot}S_0$ normalizes $P'_0$ (since 
$[P'_0]=[gP'_0g^{-1}]$) --- and this proves (2).  

Now consider the set $\Rep_{\calf_0}(P_0,S_0)= 
\Hom_{\calf_0}(P_0,S_0)/\Inn(S_0)$.  Since $P\le{}S'$, the group $P/P_0$ 
acts on this set by conjugation (i.e., $gP_0\in{}P/P_0$ acts on 
$[\varphi]$, for $\varphi\in\Hom_{\calf_0}(P_0,S_0)$, by sending it to 
$[c_g\varphi{}c_g^{-1}]$).  In particular, since the $S_0$-conjugacy class 
$[P'_0]$ is invariant under conjugation by $P/P_0$, this group leaves 
invariant the subset $X\subseteq\Rep_{\calf_0}(P_0,S_0)$ of all conjugacy 
classes $[\varphi]$ of homomorphisms such that $[\Im(\varphi)]=[P'_0]$.  
Fix any $\varphi\in\Iso_{\calf_0}(P_0,P'_0)$ (recall that subgroups in 
$\calp_{\textup{fn}}$ are $\calf_0$-conjugate to $P_0$).  Every element of 
$X$ has the form $[\alpha\varphi]$ for some 
$\alpha\in\Aut_{\calf_0}(P'_0)$, and $[\alpha\varphi]=[\beta\varphi]$ if 
and only if $\alpha\beta^{-1}\in\Aut_{S_0}(P'_0)$.  Thus 
$|X|=|\Aut_{\calf_0}(P'_0)|\big/|\Aut_{S_0}(P'_0)|$, and is prime to $p$ 
since $P'_0$ is fully normalized in $\calf_0$.  We can thus choose 
$\varphi_0\in\Hom_{\calf_0}(P_0,S_0)$ such that $\varphi_0(P_0)=P'_0$ and 
$[\varphi_0]$ is invariant under the $P/P_0$-action.  

Fix $\psi\in\Iso_{\call_0}(P_0,P'_0)$ such that $\pi(\psi)=\varphi_0$.  
Since $[\varphi_0]$ is $P/P_0$-invariant, for each $x\in{}P$, there is 
some $y\in{}x{\cdot}S_0$ such that 
$c_y\circ\varphi_0=\varphi_0\circ{}c_x$.  By (A)$_q$ (and since $\psi$ is 
an isomorphism), there is $y'\in{}y{\cdot}C_{S_0}(P'_0)$ such that 
$\q{y'}\circ\psi=\psi\circ\q{x}$ in $\call_1$.  This element $y'$ is 
unique by (1); and upon setting $\varphi(x)=y'$ we get a homomorphism 
$\varphi\in\homf(P,S)$ which extends $\varphi_0$.  Set $P'=\varphi(P)$; 
then $P'$ is $\calf$-conjugate to $P$ and $P'_0=P'\cap{}S_0$.  

\smallskip


\noindent\textbf{Step 4: }  In Step 5, we will prove that $\calf$ is 
saturated, using \cite[Theorem 2.2]{bcglo1}.  Before that theorem can be 
applied, a certain technical condition must be checked.

Assume that $P$ is $\calf$-centric, but not in $\calh$.  By Step 3, $P$ is 
$\calf$-conjugate to some $P'$ such that $P'_0$ is fully normalized in 
$\calf_0$.  Thus $P'_0$ is fully centralized in $\calf_0$ and not 
$\calf_0$-centric, which implies that $C_{S_0}(P'_0)\nleq{}P'_0$.  Then 
$P'$ acts on $C_{S_0}(P'_0){\cdot}P'_0/P'_0$ with fixed subgroup 
$QP'_0/P'_0\ne1$ for some $Q\le{}C_{S_0}(P'_0)$, and $[Q,P']\le{}P'_0$ 
since $P'/P'_0$ centralizes $QP'_0/P'_0$.  Hence $Q\nleq{}P'_0$, and 
$Q\le{}N_S(P')$ since $[Q,P']\le{}P'$.  For any $x\in{}Q{\sminus}P'_0$, 
$[c_x]\ne1\in\Out(P')$ ($C_S(P')\le{}P'$ since $P$ is $\calf$-centric), 
but $c_x$ induces the identity on $P'_0$ (since $Q\le{}C_{S_0}(P'_0)$) and 
on $P'/P'_0$ (since $x\in{}S_0$).  Hence $[c_x]\in{}O_p(\outf(P'))$ by 
Lemma \ref{gorenstein}.  This shows that
	\beq P \textup{ $\calf$-centric, } P\notin\calh\
	\Longrightarrow\ \exists\, P' \textup{ $\calf$-conjugate to $P$, } 
	\Out_S(P')\cap O_p(\outf(P'))\ne1. \tag{3} \eeq


\smallskip

\noindent\textbf{Step 5: } We next show that $\calf$ is saturated, and 
also (since it will be needed in the proof of (II)) that axiom (A)$_q$ 
holds for $\call_2$.  By \cite[Theorem 2.2]{bcglo1} (the stronger form of 
Theorem \ref{centr->sat}(b)), it suffices to prove that the subgroups in 
$\calh$ satisfy the axioms for saturation.  Note in particular that 
condition ($*$) in \cite[Theorem 2.2]{bcglo1} is precisely what is shown 
in (3).

\smallskip

\textbf{Proof of (I): } 
Fix a subgroup $P\in\calh$ which is fully normalized in $\calf$.  Let 
$S'/S_0\le{}S/S_0$ be the stabilizer of the $\calf_0$-conjugacy class of 
$P_0$.  By Step 3, $P$ is $\calf$-conjugate to a subgroup $P'$ such that 
$P'_0$ is fully normalized in $\calf_0$; and such that for each 
$g\in{}S'$, some element of $g{\cdot}S_0$ normalizes $P'_0$ (see (2)).  
Hence there are short exact sequences
	$$ 1 \Right3{} \Aut_{\call_0}(P'_0) \Right4{} \Aut_{\call_1}(P'_0) 
	\Right4{} S'/S_0 \Right1{} 1 $$
	$$ 1 \Right3{} N_{S_0}(P'_0) \Right4{} N_{S}(P'_0) 
	\Right4{} S'/S_0 \Right1{} 1. $$
We consider $P'\le{}N_S(P'_0)$ as subgroups of $\Aut_{\call_1}(P'_0)$ via 
$\widehat{\delta}_{P'_0,P'_0}$.  Then 
	$$ [\Aut_{\call_1}(P'_0):N_S(P'_0)] =
	[\Aut_{\call_0}(P'_0):N_{S_0}(P'_0)] $$
is prime to $p$ since $P'_0$ is fully normalized in $\calf_0$, and hence 
$N_S(P'_0)\in\sylp{\Aut_{\call_1}(P'_0)}$.  Fix $\psi\in\Aut_{\call_1}(P'_0)$ 
such that $\psi^{-1}N_S(P'_0)\psi$ contains a Sylow $p$-subgroup of the 
group $N_{\Aut_{\call_1}(P'_0)}(P')$ (in particular, 
$P'\le\psi^{-1}N_S(P'_0)\psi$), and set 
$P''=\psi{}P'\psi^{-1}\le{}N_S(P'_0)$.  Then 
$\psi\in\Iso_{\call_2}(P',P'')$ by 
definition of $\call_2$.  In particular, $P''$ is $\calf$-conjugate to $P$, 
and $P''_0=P'_0$.  Also, $N_S(P'_0)$ contains a Sylow $p$-subgroup of 
$N_{\Aut_{\call_1}(P'_0)}(P'')$, 
	$$ \Aut_{\call_2}(P'') = N_{\Aut_{\call_1}(P'_0)}(P'')
	\qquad\textup{and}\qquad
	N_S(P'') = N_{N_{S}(P'_0)}(P''), $$
and it follows that $N_S(P'')\in\sylp{\Aut_{\call_2}(P'')}$. 

Now, $\Aut_{\call_2}(P)\cong\Aut_{\call_2}(P'')$ since they are 
$\calf$-conjugate, and $|N_S(P)|\ge|N_S(P'')|$ since $P$ is fully 
normalized.  Thus $N_S(P)\in\sylp{\Aut_{\call_2}(P)}$, and hence 
$\Aut_S(P)\in\sylp{\autf(P)}$.  Also, $N_S(P)$ contains the kernel of the 
projection from $\Aut_{\call_2}(P)$ to $\autf(P)$; i.e., $C_S(P)$ is 
isomorphic to this kernel.  For all $Q$ which is $\calf$-conjugate to $P$, 
$C_S(Q)$ is isomorphic to a subgroup of the same kernel, so 
$|C_S(Q)|\le|C_S(P)|$, and thus $P$ is fully centralized in $\calf$.  

\smallskip

\noindent\textbf{Proof of (A)$_q$ for $\call_2$: }  It is clear from the 
construction that for any $P,Q\in\calh$, $C_S(P)$ acts freely on 
$\Mor_{\call_2}(P,Q)$ via $\widehat{\delta}_{P,P}$.  So it remains to show 
that when $P$ is fully centralized in $\calf$, then for all 
$\psi,\psi'\in\Mor_{\call_2}(P,Q)$ such that $\pi(\psi)=\pi(\psi')$, there 
is some $x\in{}C_S(P)$ such that $\psi'=\psi\circ\q{x}$.  Since every 
morphism in $\call_2$ is the composite of an isomorphism followed by an 
inclusion, it suffices to show this when $\psi$ and $\psi'$ are 
isomorphisms.  But in this case, $\psi^{-1}\psi'\in\Aut_{\call_2}(P)$ lies 
in the kernel of the map to $\autf(P)$.  We have just seen, in the proof 
of (I), that this implies there is some $x\in{}C_S(P)$ such that 
$\q{x}=\psi^{-1}\psi'$, so $\psi'=\psi\circ\q{x}$, and this is what we 
wanted to prove.


\smallskip

\textbf{Proof of (II): } Fix $\varphi\in\Hom_{\calf}(P,S)$, where 
$\varphi(P)$ is fully centralized in $\calf$.  Set $P'=\varphi(P)$.  
By definition of $\calf_2\subseteq\calf$ and of $\call_2$, there is some 
$\psi\in\Iso_{\call_2}(P,P')\subseteq\Iso_{\call_1}(P_0,P'_0)$ such that 
$\psi\circ\q{g}= \q{\varphi(g)}\circ\psi$ for all $g\in{}P$.  Upon 
replacing $\psi$ by $\q{x}\circ\psi$ for some appropriate $x\in{}S$ (and 
replacing $\varphi$ by $c_x\circ\varphi$ and $P'$ by $xP'x^{-1}$), we can 
assume that $\psi\in\Iso_{\call_0}(P_0,P'_0)$ and 
$\varphi|_{P_0}\in\Hom_{\calf_0}(P_0,S_0)$.  

Consider the subgroups
	\begin{align*} 
	N_\varphi &= \{x\in{}N_S(P) \,|\, \varphi c_x\varphi^{-1} 
	\in \Aut_S(P') \} \\
	N' = N_{\varphi|P_0}\cap{}N_\varphi &= 
	\bigl\{x\in{}N_\varphi\cap{}S_0 \,\big|\, 
	(\varphi c_x\varphi^{-1})|_{P_0} \in \Aut_{S_0}(P'_0) \bigr\} .
	\end{align*}
We will see shortly that $N'=N_\varphi\cap{}S_0$.  By (II) applied to the 
saturated fusion system $\calf_0$, there is 
$\widebar{\varphi}_0\in\Hom_{\calf_0}(N',S)$ which extends 
$\varphi|_{P_0}$, and it lifts to $\widebar{\psi}\in\Hom_{\call_0}(N',S)$. 
By (A) (applied to $\call_0$), there is $z\in{}Z(P_0)$ such that 
$\psi=(\widebar{\psi}|_{P_0})\circ\q{z}$.  Upon replacing $\widebar{\psi}$ 
by $\widebar{\psi}\circ\q{z}$ (and $\widebar{\varphi}_0$ by 
$\widebar{\varphi}_0\circ{}c_z$), we can assume that 
$\psi=\widebar{\psi}|_{P_0}$.  

For any $x\in{}N_\varphi\cap{}S_0$, 
$\varphi{}c_x\varphi^{-1}=c_y\in\autf(P')$ for some $y\in{}N_S(P')$.  
Hence by (A)$_q$ (for $\call_2$), $\psi\q{x}\psi^{-1}=\q{yz}$ for 
some unique $z\in{}C_S(P')$.  Thus $\q{yz}\in\Aut_{\call_0}(P'_0)$, and by 
definition of the distinguished monomorphisms for $\call_1$, this is 
possible only if $yz\in{}S_0$.  Thus $\varphi{}c_x\varphi^{-1}=c_{yz}$ 
where $yz\in{}N_{S_0}(P')$, and so $x\in{}N'$.  This shows that 
$N'=N_\varphi\cap{}S_0$.

Define $\widebar{\varphi}\in\Hom(N_\varphi,S)$ by the relation 
$\varphi(x)=y$ if $\psi\q{x}\psi^{-1}=\q{y}$.  In particular, this implies 
that $\q{yx^{-1}}=\psi\circ(x\psi{}x^{-1})^{-1}$ is a morphism in 
$\call_0$, and hence that $yx^{-1}\in{}S_0$.  Also, 
$\widebar{\varphi}|_P=\varphi$ by the original assumption on $\psi$. It 
remains to show that $\widebar{\varphi}\in\homf(N_\varphi,S)$.  To do 
this, it suffices to show that $\widebar{\psi}\circ\q{x}= 
\q{y}\circ\widebar{\psi}$ in $\Mor_{\call_2}(N',S)$ for all $x$, where 
$y=\widebar{\varphi}(x)$.  Equivalently, we must show that
	\beq \widebar{\psi} = \q{yx^{-1}} \circ 
	(x\widebar{\psi}x^{-1}). \tag{4} \eeq
Since $yx^{-1}\in{}S_0$, both sides in (4) are in $\call_0$, and they are 
equal after restriction to $P_0$.  Hence they are equal as morphisms 
defined on $N_\varphi\cap{}S_0$ by \cite[Lemma 3.9]{bcglo1}, and this 
finishes the proof.

\smallskip

\noindent\textbf{Step 6: }  We next check that $\calf_0$ has $p$-power 
index in $\calf$.  For any $P\le{}S$ and any $\alpha\in\autf(P)$ of order 
prime to $P$, $\alpha$ induces the identity on $P/P_0$ by construction, 
and hence $x^{-1}\alpha(x)\in{}S_0$ for all $x\in{}P$.  This shows that 
$S_0\ge\hfocal{\calf}{S}$ (see Definition \ref{O^p(F)}).  Also, by 
construction, for all $P\le{}S_0$, $\Aut_{\calf_0}(P)$ is normal of 
$p$-power index in $\autf(P)$, and thus contains $O^p(\autf(P))$.  This 
proves that the fusion subsystem $\calf_0$ has $p$-power index in $\calf$ 
in the sense of Definition \ref{def-p'-index}.  

\smallskip

\noindent\textbf{Step 7: }  It remains to construct a quasicentric linking 
system $\callq$ which contains $\call_2$ as a full subcategory, and which 
is associated to $\calf$.  Note first that the axioms of Definition 
\ref{L^q} are all satisfied by $\call_2$:  axiom (A)$_q$ holds by Step 5, 
while axioms (B)$_q$, (C)$_q$, and (D)$_q$ follow directly from the 
construction in Steps 1 and 2.  

Let $\call^c_2\subseteq\call_2$ be the full subcategory whose objects are 
the set $\calh^c\subseteq\calh$ of subgroups in $\calh$ which are 
$\calf$-centric.  We first construct a centric linking system 
$\call\supseteq\call^c_2$ associated to $\calf$.  For any set $\calk$ of 
$\calf$-centric subgroups of $S$, let $\orb^c(\calf)$ be the orbit 
category of $\calf$, let $\orb^\calk(\calf)\subseteq\orb^c(\calf)$ be the 
full subcategory with object set $\calk$, and let $\calz_\calf^{\calk}$ 
be the functor $P\mapsto{}Z(P)$ on $\orb^\calk(\calf)$ which sends (see 
Definition \ref{D:orbit} for more detail).  
By \cite[Proposition 3.1]{BLO2}, when $\calk$ is 
closed under $\calf$-conjugacy and overgroups, the obstruction to the 
existence of a linking system with object set $\calk$ lies in 
$\higherlim{}3(\calz_\calf^{\calk})$, and the 
obstruction to its uniqueness lies in 
$\higherlim{}2(\calz_\calf^{\calk})$.  Furthermore, by \cite[Proposition 
3.2]{BLO2}, if $\calp$ is an $\calf$-conjugacy class of $\calf$-centric 
subgroups maximal among those not in $\calk$, and $\calk'=\calk\cup\calp$, 
then the inclusion of functors induces an isomorphism between the higher 
limits of $\calz_\calf^{\calk}$ and $\calz_\calf^{\calk'}$ if certain 
groups $\Lambda^*(\outf(P);Z(P))$ vanish for $P\in\calp$.  By (3), 
$O_p(\outf(P))\ne1$ for any $\calf_0$-centric subgroup $P\notin\calh^c$, 
and hence $\Lambda^*(\outf(P);Z(P))=0$ for such $P$ by \cite[Proposition 
6.1(ii)]{JMO}.  So these obstructions all vanish, and there is a centric 
linking system $\call\supseteq\call_2^c$ associated to $\calf$.

Now let $\callq$ be the quasicentric linking system associated to $\SFL$.  
For each $P\in\calh=\Ob(\call_2)$, $P\cap{}S_0$ is $\calf_0$-centric by 
definition, hence is $\calf$-quasicentric by Theorem \ref{F:p-quot.}(a), 
and thus $P$ is also $\calf$-quasicentric.  Also, $\calh$ contains all 
$\calf$-centric $\calf$-radical subgroups by (3), and $\calh$ is closed 
under $\calf$-conjugacy and overgroups (by definition of $\calf$ and 
$\calh$).  Hence by \cite[Proposition 3.12]{bcglo1}, since $\call_2^c$ is 
a full subcategory of $\call$ by construction, $\call_2$ is isomorphic to 
a full subcategory of $\callq$ in a way which preserves the projection 
functors and distinguished monomorphisms.  So $\call_0$ can also be 
identified with a linking subsystem of $\callq$, and this finishes the 
proof.  
\end{proof}

We now prove a topological version of Theorem \ref{SFL_0<SFL}. By Theorem 
\ref{L:p-quot.}, if $\SFL[_0]\subseteq\SFL$ is an inclusion of $p$-local 
finite groups of $p$-power index, where $S_0\nsg{}S$ and $\Gamma=S/S_0$, 
then there is a fibration sequence 
$|\call_0|\pcom\Right2{}|\call|\pcom\Right2{}B\Gamma$.  So it is natural 
to ask whether the opposite is true:  given a fibration sequence whose 
base is the classifying space of a finite $p$-group, and whose fiber is 
the classifying space of a $p$-local finite group, is the total space also 
the classifying space of a $p$-local finite group?  The next proposition 
shows that this is, in fact, the case.

Before stating the proposition, we first define some categories which will 
be needed in its proof.  Fix a space $Y$, a $p$-group $S$, and a map 
$f\:BS\Right2{}Y$.  For $P\le{}S$, we regard $BP$ as a subspace of $BS$; 
all of these subspaces contain the basepoint $*\in{}BS$.  We define three 
categories in this situation, $\calf_{S,f}(Y)$, $\call_{S,f}(Y)$, and 
$\calm_{S,f}(Y)$, all of which have as objects the subgroups of $S$.  Of 
these, the first two are discrete categories, while $\calm_{S,f}(Y)$ has a 
topology on its morphism sets.  Morphisms in $\calf_{S,f}(Y)$ are defined 
by setting
	$$ \Mor_{\calf_{S,f}(Y)}(P,Q) = \bigl\{\varphi\in\Hom(P,Q) 
	\,\big|\, f|_{BP}\simeq f|_{BQ}\circ B\varphi \bigr\}; $$
we think of this as the fusion category of $Y$ (with respect to $S$ and $f$).  

Next define
	\begin{multline*} 
	\Mor_{\calm_{S,f}(Y)}(P,Q) = \bigl\{(\varphi,H) \,\big|\, 
	\textup{$\varphi\in\Hom(P,Q)$, \ $H\:BP\times[0,t]\rTo Y$,} \\
	\textup{$t\ge0$, \ $H|_{BP\times0}=f|_{BP}$, \
	$H|_{BP\times{}t}=f|_{BQ}\circ{}B\varphi$} \bigr\} \,.
	\end{multline*}
Thus a morphism in $\calm_{S,f}(Y)$ has the form $(\varphi,H)$ 
where $H$ is a \emph{Moore homotopy} in $Y$.  Composition is defined by 
	$$ (\psi,K)\circ(\varphi,H)=
	(\psi\varphi,(K\circ(B\varphi\times\Id)){\cdot}H), $$
where if $H$ and $K$ are homotopies parameterized by $[0,t]$ and $[0,s]$, 
respectively, then $(K\circ(B\varphi\times\Id)){\cdot}H$ is the composite 
homotopy parameterized by $[0,t+s]$.  

Let $\calp(Y)$ be the category of Moore paths in $Y$, and let 
$\Res_*\:\calm_{S,f}(Y)\rTo\calp(Y)$ be the functor which sends each 
object to $f(*)$, and sends a morphism $(\varphi,H)$ to the path obtained 
by restricting $H$ to the basepoint $*\in{}BS$.  Define a map 
$\textup{ev}$ from $|\calp(Y)|$ to $Y$ as follows.  For any $n$-simplex 
$\Delta^n$ in $|\calp(Y)|$, indexed by a composable sequence of paths 
$\phi_1,\dots,\phi_n$ where $\phi_i$ is defined on the interval 
$[0,t_i]$, let $\textup{ev}|_{\Delta^n}$ be the composite
	$$ \Delta^n \Right5{\lambda(t_1,\dots,t_n)} [0,t_1+\ldots+t_n] 
	\Right5{\phi_n\cdots\phi_1} Y, $$
where $\lambda(t_0,\dots,t_n)$ is the affine map which sends the $i$-th 
vertex to $t_1+\ldots+t_{i}$.  The category $\calm_{S,f}(Y)$ thus comes 
equipped with an ``evaluation function''
	$$ \textup{eval}\: |\calm_{S,f}(Y)| \Right4{|\Res_*|} 
	|\calp(Y)| \Right4{\textup{ev}} Y. $$

Now set 
	$$ \Mor_{\call_{S,f}(Y)}(P,Q) = 
	\pi_0\bigl(\Mor_{\calm_{S,f}(Y)}(P,Q)\bigr). $$
We think of $\call_{S,f}(Y)$ as the linking category of $Y$.  Also, for 
any set $\calh$ of subgroups of $S$, we let $\call_{S,f}^{\calh}(Y)$ and 
$\calm_{S,f}^{\calh}(Y)$ denote the full subcategories of $\call_{S,f}(Y)$ 
and $\calm_{S,f}(Y)$ with object set $\calh$.  For more about the fusion 
and linking categories of a space, see \cite[\S7]{BLO2}.  

\begin{Thm} \label{|L|->X->Bpi}
Fix a $p$-local finite group $\SFL[_0]$, a $p$-group $\Gamma$, and a 
fibration $X\Right2{v}B\Gamma$ with fiber $X_0\simeq|\call_0|\pcom$.  Then 
there is a $p$-local finite group $\SFL$ such that $S_0\nsg{}S$, 
$\calf_0\subseteq\calf$ is a fusion subsystem of $p$-power index, 
$S/S_0\cong{}\Gamma$, and $X\simeq|\call|\pcom$.
\end{Thm}

\begin{proof}  Let $*$ denote the base point of $B\Gamma$, and assume 
$X_0=v^{-1}(*)$.  Fix a homotopy equivalence 
$f\:|\call_0|\pcom\Right1{}X_0$, regard $BS_0$ as a subspace of 
$|\call_0|$, and set $f_0=f|_{BS_0}\:BS_0\Right1{}X_0$, also regarded as a 
map to $X$.  Let $\calh_0$ be the set of $\calf_0$-centric subgroups of 
$S_0$.  

\smallskip

\noindent\textbf{Step 1: }  By \cite[Proposition 7.3]{BLO2}, 
	$$ \calf_0\cong \calf_{S_0,f_0}(X_0) \qquad\textup{and}\qquad
	\call_0 \cong \call_{S_0,f_0}^{\calh_0}(X_0). $$ 
We choose the inclusions $\iota_P^Q\in\Mor_{\call_0}(P,Q)$ (for $P\le{}Q$ 
in $\calh_0$) to correspond to the morphisms $(\incl_P^Q,[c])$ in 
$\call_{S_0,f_0}(X_0)$, where $c$ is the constant homotopy $f_0|_{BP}$.  
Set
	$$ \calf_1 = \calf_{S_0,f_0}(X) \qquad\textup{and}\qquad
	\call_1 = \call_{S_0,f_0}^{\calh_0}(X), $$ 
where $f_0$ is now being regarded as a map $BS_0\Right2{}X$.  
The inclusion $X_0\subseteq{}X$ makes $\calf_0$ into a subcategory of 
$\calf_1$ and $\call_0$ into a subcategory of $\call_1$.  

For all $P\le{}S_0$ which is fully centralized in $\calf_0$, 
$\map(BP,X_0)_{f_0|BP}\simeq|C_{\call_0}(P)|\pcom$ by \cite[Theorem 
6.3]{BLO2}, where $C_{\call_0}(P)$ is a linking system over the 
centralizer $C_S(P)$.  Since $P\in\calh_0$ (i.e., $P$ is 
$\calf_0$-centric) if and only if it is fully centralized and 
$C_{S_0}(P)=Z(P)$, this shows that
	\beq \textup{$P\in\calh_0$ $\Longleftrightarrow$ 
	$\map(BP,X_0)_{f_0|_{BP}}\simeq{}BZ(P)$.} \tag{1} \eeq
If $P$ and $P'$ are $\calf_1$-conjugate, and 
$\varphi\in\Iso_{\calf_1}(P,P')$, then the homotopy between $f_0|_{BP}$ 
and $f_0|_{BP'}\circ{}B\varphi$ as maps from $BP$ to $X$ induces, using 
the homotopy lifting property for the fibration $v$, a homotopy between 
$w\circ{}f_0|_{BP}$ and $f_0|_{BP'}\circ{}B\varphi$ (as maps from $BP$ to 
$X_0$), where $w\:X_0\rTo^{\simeq}X_0$ is the homotopy equivalence induced 
by lifting some loop in $B\Gamma$.  Since $w$ is a homotopy equivalence 
and $B\varphi$ is a homeomorphism, this shows that the mapping spaces 
$\map(BP,X_0)_{f_0|_{BP}}$ and $\map(BP',X_0)_{f_0|_{BP'}}$ are homotopy 
equivalent, and hence (by (1)) that $P'\in\calh_0$ if 
$P'\in\calh_0$.  Thus for all $P,P'\le{}S_0$,
	\beq \textup{$P$ $\calf_1$-conjugate to $P'$ and $P\in\calh_0$ 
	$\Longrightarrow$ $P'\in\calh_0$.} \tag{2} \eeq

Fix $S\in\sylp{\Aut_{\call_1}(S_0)}$.  We identify $S_0$ as a subgroup of 
$S$ via the distinguished monomorphism $\delta_{S_0}$ from $S_0$ to 
$\Aut_{\call_0}(S_0)\le\Aut_{\call_1}(S_0)$.  

\smallskip

\noindent\textbf{Step 2: }  For all $P\le{}S$, and for $Y=X_0$ or $X$, we 
define 
	$$ \map(BP,Y)_\Phi= \bigl\{f\:BP\rTo Y \,\big|\,
	f\simeq f_0\circ B\varphi,\ \ \varphi\in\Hom(P,S),\ \ 
	\varphi(P)\in\calh_0 \bigr\}. $$
Using (2), we see that the fibration sequence 
$X_0\Right2{}X\Right2{}B\Gamma$ induces a fibration sequence of mapping 
spaces
	\beq \map(BP,X_0)_\Phi \Right4{} \map(BP,X)_\Phi \Right4{} 
	\map(BP,B\Gamma)_{ct}\simeq{}B\Gamma, \tag{3} \eeq
where $\map(BP,B\Gamma)_{ct}$ is the space of null homotopic maps, and the 
last equivalence is induced by evaluation at the basepoint.  By (1), each 
connected component of $\map(BP,X_0)_\Phi$ is homotopy equivalent to 
$BZ(P)$, and hence the connected components of $\map(BP,X)_\Phi$ are 
also aspherical.  

For any morphism $(\varphi,[H])\in\Mor_{\call_1}(P,Q)$, where 
$\varphi\in\Hom_{\calf_1}(P,Q)$ and $[H]$ is the homotopy class of the 
path $H$ in $\map(BP,X)_\Phi$, restricting $v\circ{}H$ to the basepoint of 
$BP$ defines a loop in $B\Gamma$, and thus an element of $\Gamma$.  This 
defines a map from $\Mor(\call_1)$ to $\Gamma$ which sends 
composites to products, and thus a functor
	$$ \widehat{\theta}\: \call_1 \Right4{} \calb(\Gamma). $$
By construction (and the fibration sequence (3)), for any
$\psi\in\Mor(\call_1)$, we have $\psi\in\Mor(\call_0)$ if and only if 
$\widehat{\theta}(\psi)=1$.  

Using the homotopy lifting property in (3), we see that $\widehat{\theta}$ 
restricts to a surjection of $\Aut_{\call_1}(S_0)$ onto $\Gamma$, with 
kernel $\Aut_{\call_0}(S_0)$.  Since $S_0\in\sylp{\Aut_{\call_0}(S_0)}$, 
this shows that $\widehat{\theta}$ induces an isomorphism 
$S/S_0\cong\Gamma$, where $S$ is a Sylow $p$-subgroup of 
$\Aut_{\call_1}(S_0)$ which contains $S_0$.  Hence for any 
$\psi\in\Mor(\call_1)$, there is $g\in{}S$ such that 
$\widehat{\theta}(\psi)=\widehat{\theta}(g)$; and for any such $g$ there 
is a unique morphism $\psi_0\in\Mor(\call_0)$ such that 
$\psi=\delta(g)\circ\psi_0$ where $\delta(g)\in\Mor(\call_1)$ denotes the 
appropriate restriction of $g\in\Aut_{\call_1}(S_0)$.  In other words,
	\beq \Mor(\call_1) \cong S \times_{S_0} \Mor(\call_0). \tag{4} \eeq

\smallskip

\noindent\textbf{Step 3: }  The conjugation action of $S$ on 
$\Mor(\call_0)\subseteq\Mor(\call_1)$ defines an action of $S$ on 
$\call_0$, which satisfies the hypotheses of Theorem \ref{SFL_0<SFL}.  
(Note in particular that this action sends inclusions to inclusions, since 
they are assumed to be represented by constant homotopies.)  So we can now 
apply that theorem to construct a $p$-local finite group $\SFL$ which 
contains $\SFL[_0]$ with $p$-power index.  Let $\callq$ be the 
quasicentric linking system associated to $\SFL$.  By (4), the category 
$\call_1$ defined here is equal to the category $\call_1$ defined in the 
proof of Theorem \ref{SFL_0<SFL}; i.e., the full subcategory of $\callq$ 
with object set $\calh_0$. 

Let $\calh$ be the set of subgroups $P\le{}S$ 
such that $P\cap{}S_0\in\calh_0$, and let $\call_2\subseteq\callq$ be the 
full subcategory with object set $\calh$.  By Step 4 in the proof of Theorem 
\ref{SFL_0<SFL}, all $\calf$-centric $\calf$-radical subgroups of $S$ lie 
in $\calh$, and so $|\call|\simeq|\callq|\simeq|\call_2|$ by Proposition 
\ref{L-props}(a).  Also, $|\call_1|$ is a deformation retract of 
$|\call_2|$, where the retraction is defined by sending $P\in\calh$ to 
$P\cap{}S_0\in\calh_0$ (and morphisms are sent to their restrictions, 
uniquely defined by Proposition \ref{L-props}(b)).  Thus 
$|\call|\simeq|\call_1|$.  So by Theorem \ref{L:p-quot.}, we have a 
homotopy fibration sequence $|\call_0|\pcom\rTo|\call_1|\pcom\rTo 
B\Gamma$.  

\smallskip

\noindent\textbf{Step 4: } It remains to construct a homotopy equivalence 
$|\call_1|\pcom\Right1{}X$, which extends to a homotopy equivalence 
between the fibration sequences.  This is where we need to use the 
topological linking categories defined above.  Set
	$$ \calm_0 = \calm_{S_0,f_0}^{\calh_0}(X_0) 
	\qquad\textup{and}\qquad
	\calm_1 = \calm_{S_0,f_0}^{\calh_0}(X) $$
for short, and consider the following commutative diagram:
	\begin{diagram}[w=40pt] 
	|\call_0|\pcom & \lTo^{\proj_0}_{\simeq} & |\calm_0|\pcom 
	& \rTo^{\textup{eval}_0} & X_0 \\
	\dTo && \dTo && \dTo \\
	|\call_1|\pcom & \lTo^{\proj_1}_{\simeq} & |\calm_1|\pcom 
	& \rTo^{\textup{eval}_1} & X \rlap{\,.}
	\end{diagram} 
The vertical maps in the diagram are all inclusions.  Also, $X_0$ is 
$p$-complete by definition and $X$ by \cite[II.5.2(iv)]{BK}, so the 
evaluation maps defined above extend to the $p$-completed nerves 
$|\calm_i|\pcom$.  The maps $\proj_1$ and $\textup{eval}_1$ both commute 
up to homotopy with the projections to $B\Gamma$.  The projection maps 
$\proj_0$ and $\proj_1$ are both homotopy equivalences:  the connected 
components of the morphism spaces in $\calm_i$ are contractible since the 
connected components of the fiber and total space in (3) are aspherical 
(see Step 2).  

We claim that $\textup{eval}_0$ is homotopic to $f\circ\proj_0$ as maps to 
$X_0$.  By naturality, it suffices to check this on the uncompleted nerve 
$|\calm_0|$, and in the case where $X_0=|\call|\pcom$ and $f=\Id$.  But in 
this case, the only real difference between the maps is that 
$\textup{eval}_0$ sends all vertices of $|\calm_0|$ to the base point of 
$X_0=|\call|\pcom$, while $\proj_0$ sends the vertex for a subgroup 
$P\le{}S$ to the corresponding vertex in $|\call|\pcom$.  So the maps are 
homotopic via a homotopy which sends vertices to the base point along the 
edges of $|\call|$ corresponding to the inclusion morphisms.  In 
particular, this shows that $\textup{eval}_0$ is also a homotopy 
equivalence.  We thus have an equivalence between the fibration sequences
	\begin{diagram}[w=30pt] 
	|\call_0|\pcom & \rTo & |\call_1|\pcom & \rTo & B\Gamma \\
	\dTo<{\simeq} && \dTo && \dIgual \\
	|\call_0|\pcom & \rTo & X & \rTo & B\Gamma \rlap{\,,}
	\end{diagram}
and this proves that the two sequences are equivalent.
\end{proof}

Recall (Definition \ref{D:autfus}) that for any saturated fusion system 
$\calf$ over a $p$-group $S$, $\Aut\fus(S,\calf)$ denotes the group of all 
fusion preserving automorphisms of $S$.  The following corollary to 
Proposition \ref{SFL_0<SFL} describes how ``exotic'' fusion systems could 
potentially arise as extensions of $p$-power index; we still do not know 
whether the situation it describes can occur.

\begin{Cor} \label{exotic-ext}
Fix a finite group $G$, with Sylow subgroup $S\in\sylp{G}$.  Assume there 
is an automorphism $\widebar{\alpha}\in\Aut\fus(S,\calf_S(G))$ of 
$p$-power order, which is not the restriction to $S$ of an automorphism of 
$G$, and which moreover is not the restriction of an automorphism of $G'$ 
for any finite group $G'$ with $S\in\sylp{G'}$ and 
$\calf_S(G')=\calf_S(G)$.  Then there is a saturated fusion system 
$(\widehat{S},\widehat{\calf})\supseteq (S,\calf_S(G))$, such that 
$\calf_S(G)$ has $p$-power index in $\widehat{\calf}$, and such that 
$\widehat{\calf}$ is not the fusion system of any finite group.
\end{Cor}

\begin{proof}  By \cite[Theorem E]{BLO1}, together with \cite[Theorem 
A]{Oliver-odd} and \cite[Theorem A]{Oliver-2}, there is a short exact 
sequence
	$$ 0 \Right2{} \higherlim{\orb_S^c(G)}1(\calz_G) \Right5{}
	\Out\isotyp(\call_S^c(G)) \Right5{} \Out\fus(S,\calf_S(G)) 
	\Right2{} 0 $$
where $\Out\fus(S,\calf)=\Aut\fus(S,\calf)/\Aut_\calf(S)$, where 
$\orb_S^c(G)$ and $\calz_G$ are the category and functor of Definition 
\ref{D:orbit}(b), and where $\Out\isotyp(\call_S^c(G))$ is 
the group of ``isotypical'' automorphisms of $\call_S^c(G)$ modulo natural 
isomorphism (see the introduction of \cite{BLO1}, or \cite[Definition 
3.2]{BLO1}, for the definition).  Let $[\widebar{\alpha}]$ be the class of 
$\widebar{\alpha}$ in $\Out\fus(S,\calf_S(G))$; then $[\widebar{\alpha}]$ 
lifts to an automorphism $\alpha$ of the linking system $\call_S^c(G)$.  
Upon replacing $\alpha$ by some appropriate power, we can assume that it 
still has $p$-power order.  We can also assume, upon replacing $\alpha$ by 
another automorphism in the same conjugacy class if necessary, that the 
$\alpha$-action on $\Aut_{\call}(S)$ leaves invariant the subgroup 
$\delta_S(S)$; i.e., that the action of $\alpha$ on $\call_S^c(G)$ 
restricts to an action on $S$.

Set $\widehat{S}=S{\rtimes}\gen{x}$, where $|x|=|\alpha|$ and $x$ acts on 
$S$ via $\alpha$.  Then $\widehat{S}$ has an action on $\call_S^c(G)$ 
induced by the actions of $S$ and of $\alpha$, and this action satisfies 
conditions (a)--(d) in Theorem \ref{SFL_0<SFL}.  Let 
$\widehat{\calf}\supseteq\calf_S(G)$ be the saturated fusion system over 
$\widehat{S}$ constructed by the theorem.  

We claim that $\widehat{\calf}$ is not the fusion system 
of any finite group.  Assume otherwise:  assume $\widehat{\calf}$ is the 
fusion system of a group $\widehat{G}$.  Since $\calf$ has $p$-power index 
in $\widehat{\calf}$, $S\ge\hfocal{\widehat{\calf}}{\widehat{S}}$, and so 
$S\ge\widehat{S}\cap{}O^p(\widehat{G})$ by the hyperfocal subgroup 
theorem (Lemma \ref{G:hyperfocal}).  Set $G'=S{\cdot}O^p(\widehat{G})$.  
Then $G'\nsg\widehat{G}$ since $S\nsg\widehat{S}$; 
$\widehat{G}/G'\cong\widehat{S}/S'$, and thus $G'\nsg\widehat{G}$ has
$p$-power index and $S\in\sylp{G'}$.  By Theorem \ref{F:p-quot.}, 
there is a unique saturated fusion subsystem over $S$ of $p$-power index 
in $\widehat{\calf}$, and thus $\calf=\calf_S(G')$.  Also, 
$x\in\widehat{S}\le\widehat{G}$ acts on $G'$ via an automorphism whose 
restriction to $S$ is $\alpha$, and this contradicts the original 
assumption about $\widebar\alpha$.  
\end{proof}


\newsect{Fusion subsystems and extensions of index prime to $p$.} 

In this section, we classify all saturated fusion subsystems of index 
prime to $p$ in a given saturated fusion system $\calf$, and show that 
there is a unique minimal subsystem $\Oppf$ of this type.  More precisely, 
we show that there is a certain finite group of order prime to $p$ 
associated to $\calf$, denoted below $\gppf$, and a one-to-one 
correspondence between subgroups $T\le\gppf$ and fusion subsystems 
$\calf_T$ of index prime to $p$ in $\calf$.  The index of $\calf_T$ in 
$\calf$ can then be defined to be the index of $T$ in $\gppf$.

Conversely, we also describe extensions of saturated fusion
systems of index prime to $p$. Once more the terminology requires
motivation. Roughly speaking, an extension of index prime to $p$
of a given saturated fusion system $\calf$ is a saturated fusion
system $\calf'$ over the same $p$-group $S$, but where the
morphism set has been ``extended'' by an action of a group of
automorphisms whose order is prime to $p$. Here again, a
one-to-one correspondence statement is obtained, thus providing a
full classification.

\bigskip

\newsub{Subsystems of index prime to $p$}
We first classify all saturated fusion subsystems of index prime
to $p$ in a given saturated fusion system $\calf$.  The fusion
subsystems and associated linking systems will be constructed
using Proposition \ref{F:p-solv.quot.} and Theorem
\ref{L:p-p'-quot.}, respectively.  More precisely, Theorem
\ref{L:p-p'-quot.} has already told us that for any $p$-local
finite group $\SFL$, any surjection $\theta$ of $\pi_1(|\call|)$
onto a finite $p'$-group $\Gamma$, and any subgroup $H\le\Gamma$,
there is a $p$-local finite subgroup $\SFL[']$ such that
$|\call'|$ is homotopy equivalent to the covering space of
$|\call|$ with fundamental group $\theta^{-1}(H)$.  What is new in
this section is first, that we describe the ``universal''
$p'$-group quotient of $\pi_1(|\call|)$ as a certain quotient
group of $\outf(S)$; and second, we show that all fusion
subsystems of index prime to $p$ in $\calf$ (in the sense of
Definition \ref{def-p'-index}) are obtained in this way.

When applying Proposition \ref{F:p-solv.quot.} to this situation,
we need to consider fusion mapping triples
$(\Gamma,\theta,\Theta)$ on $\calfq$, where $\Gamma$ is finite of
order prime to $p$.  Since $\theta\in\Hom(S,\Gamma)$, it must then
be the trivial homomorphism.  In this case, conditions (i)--(iii)
in Definition \ref{D:J,j} are equivalent to requiring that there
is some functor $\widehat{\Theta}\:\calfq\rTo\calb(\Gamma)$ such
that $\Theta(\varphi)=\{\widehat{\Theta}(\varphi)\}$ for each
$\varphi\in\Mor(\calfq)$ (and condition (iv) is redundant).  So
instead of explicitly constructing fusion mapping triples, we
instead construct functors of this form.

We start with some definitions.  For a finite group $G$, one
defines $O^{p'}(G)$ to be the smallest normal subgroup of $G$ of
index prime to $p$, or equivalently the subgroup generated by
elements of $p$-power order in $G$.  These two definitions are
not, in general, equivalent in the case of an infinite group, (the
case $G=\Z$ being an obvious example).  We do, however, need to
deal with such subgroups here. The following definition is most
suitable for our purposes, and it is a generalization of the
finite case.

\newcommand{\oppg}[1]{O^{p'}(#1)}  

\begin{Defi}
For any group $G$ (possibly infinite), let $\oppg{G}$ be the
intersection of all normal subgroups in $G$ of finite index prime
to $p$.
\end{Defi}

In particular, under this definition, an epimorphism
$\alpha\:G\Onto2{}H$ with $\Ker(\alpha)\le{}\oppg{G}$ induces an
isomorphism $G/\oppg{G}\cong{}H/\oppg{H}$.  Thus by Proposition
\ref{pi1|L|->|F|}, for any $p$-local finite group $\SFL$, the
projections of $|\call|\simeq|\callq|$ onto $|\calfc|$ and
$|\calfq|$ induce isomorphisms
    $$ \pi_1(|\call|)\big/\oppg{\pi_1(|\call|)} \cong
    \pi_1(|\calfc|)\big/\oppg{\pi_1(|\calfc|)} \cong
    \pi_1(|\calfq|)\big/\oppg{\pi_1(|\calfq|)}. $$

Fix a saturated fusion system $\calf$ over a $p$-group $S$, and
define
    \[ \gpp{\calf} = \pi_1(|\calfc|)/\oppg{\pi_1(|\calfc|)}. \]
We will show that the natural functor
    \[ \varepsilon_{\calfc}\:\calfc \Right5{} \calb(\gpp{\calf}) \]
induces a bijective correspondence between subgroups of
$\gpp{\calf}$ and fusion subsystems of $\calf$ of index prime to
$p$.

Recall (Definition \ref{D:oppf:opf}) that for any saturated fusion
system $\calf$ over a $p$-group $S$, $\oppf\subseteq\calf$ is the
smallest fusion subsystem which contains the groups
$O^{p'}(\autf(P))$ for all $P\le{}S$; i.e., the smallest fusion
subsystem which contains all automorphisms in $\calf$ of $p$-power
order. Define
    $$ \outf^0(S) = \bigl\langle \alpha\in\outf(S) \,\big|\,
    \alpha|_P\in\Mor_{\oppf}(P,S), \textup{ some $\calf$-centric
    $P\le{}S$} \bigr\rangle. $$
Then $\outf^0(S)\nsg \outf(S)$, since $\oppf$ is normalized by
$\autf(S)$ (Lemma \ref{F=F'.N_F(S)}(a)).

\begin{Prop}  \label{Theta:p'-quot.} \label{hat{theta}}
There is a unique functor
    $$ \widehat{\theta}\: \calfc \Right7{} \calb(\outf(S)/\outf^0(S))
    $$
with the following properties:
\begin{enumerate}
\item $\widehat{\theta}(\alpha)=\alpha$ (modulo $\outf^0(S)$) for
all $\alpha\in\autf(S)$.

\item $\widehat{\theta}(\varphi)=1$ if $\varphi\in\Mor(\oppfc)$.
In particular, $\widehat{\theta}$ sends inclusion morphisms to the
identity.
\end{enumerate}
Furthermore, there is an isomorphism
    $$ \widebar{\theta}\: \gpp{\calf} =
    \pi_1(|\calfc|)/\oppg{\pi_1(|\calfc|)} \Right5{\cong}
    \outf(S)/\outf^0(S) $$
such that
$\widehat{\theta}=\calb{\widebar{\theta}}\circ\varepsilon^c_\calf$.
\end{Prop}

\begin{proof}  By Lemma \ref{F=F'.N_F(S)}(c), each morphism in $\calfc$
factors as the composite of the restriction of a morphism in
$\autf(S)$ followed by a morphism in $\oppfc$.  If
    $$ \varphi = \varphi_1\circ\alpha_1|_P = \varphi_2\circ\alpha_2|_P, $$
where $\varphi\in\homf(P,Q)$, $\alpha_i\in\autf(S)$, and
$\varphi_i\in\Hom_{\oppfc}(\alpha_i(P),Q)$, then we can assume
(after factoring out inclusions) that all of these are
isomorphisms, and hence
    $$ (\alpha_2\circ\alpha_1^{-1})|_P = \varphi_2^{-1}\circ \varphi_1
    \in \Iso_{\oppfc}(\alpha_1(P),\alpha_2(P)). $$
Thus $\alpha_2\circ\alpha_1^{-1}\in\outf^0(S)$; and so we can
define
    $$ \widehat{\theta}(\varphi) = [\alpha_1] =
    [\alpha_2] \in \outf(S)/\outf^0(S). $$
This shows that $\widehat{\theta}$ is well defined on morphisms
(and sends each object in $\calfc$ to the unique object of
$\calb(\outf(S)/\outf^0(S))$).  By Lemma \ref{F=F'.N_F(S)}(c)
again, $\widehat{\theta}$ preserves compositions, and hence is a
well defined functor. It satisfies conditions (a) and (b) above by
construction. The uniqueness of $\widehat{\theta}$ is clear.

It remains to prove the last statement.  Since
$\outf(S)/\outf^0(S)$ is finite of order prime to $p$,
$\pi_1(|\widehat{\theta}|)$ factors through a homomorphism
    $$ \bar{\theta}\: \pi_1(|\calfc|) / \oppg{\pi_1(|\calfc|)}
    \Right5{} \outf(S)/\outf^0(S). $$
The inclusion of $B\autf(S)$ into $|\calfc|$ (as the subcomplex
with one vertex $S$) induces a homomorphism
    $$ \tau\: \outf(S) \Right5{}
    \pi_1(|\calfc|)/\oppg{\pi_1(|\calfc|)}. $$
Furthermore, $\tau$ is surjective since
$\calf=\gen{\oppf,\autf(S)}$ (Lemma \ref{F=F'.N_F(S)}(b)), and
since any automorphism in $\oppf$ is a composite of restrictions
of automorphisms of $p$-power order. By (a), and since $\theta$
restricted to $\autf(S)$ is the projection onto $\outf(S)$, the
composite $\widebar{\theta}\circ\tau$ is the projection of
$\outf(S)$ onto the quotient group $\outf(S)/\outf^0(S)$.
Finally, $\outf^0(S)\le\Ker(\tau)$ by definition of $\outf^0(S)$,
and this shows that $\widebar{\theta}$ is an isomorphism.
\end{proof}

The following lemma shows that any fusion mapping triple on
$\calf^c$ can be extended uniquely to $\calf^q$. This will allow
us later to apply  Proposition \ref{F:p-solv.quot.} in order to
produce saturated fusion subsystems of index prime to $p$.

\begin{Lem} \label{extend:Fc-Fq}
Let $\calf$ be a saturated fusion system over a $p$-group $S$, and
let $(\Gamma,\theta,\Theta_0)$ be a fusion mapping triple on
$\calfc$. Then there is a unique extension
\[\Theta\colon\Mor(\calfq)\rTo\sset(\Gamma)\]
of $\Theta_0$, such that $(\Gamma,\theta,\Theta)$ is a fusion
mapping triple  on $\calfq$.
\end{Lem}

\begin{proof} We construct the extension $\Theta$ one $\calf$-conjugacy class at a time.
Thus, assume $\Theta$ has been defined on a set $\calh_0$ of
$\calf$-quasicentric subgroups of $S$ which is a union of
$\calf$-conjugacy classes, contains all $\calf$-centric subgroups,
and is closed under overgroups.  Let $\calp$ be maximal among
$\calf$-conjugacy classes of $\calf$-quasicentric subgroups not in
$\calh_0$; we show that $\Theta$ can be extended to
$\calh=\calh_0\cup\calp$.

Fix $P\in\calp$ which is fully normalized in $\calf$.  For each
$\alpha\in\autf(P)$, there is an extension
$\widebar{\alpha}\in\autf(P{\cdot}C_S(P))$ (axiom (II)), and we
define a map
	\[\Theta_P\colon \Aut_\calf(P)\rTo \sset(N_\Gamma(\theta(C_S(P))))\]
by $ \Theta_P(\alpha) = \Theta(\widebar{\alpha}){\cdot} \theta(C_S(P))$.
By Definition \ref{D:J,j} ((i) and (ii))
$\Theta(\widebar{\alpha})$ is a left coset of
$\theta(C_S(P{\cdot}C_S(P)))$, and by (iv), it is also a right
coset (where the right coset representative can be taken to be the
same as the one representing the left coset); hence
$\Theta_P(\alpha)$ is a left and a right coset of $\theta(C_S(P))$
(again with the same coset representative on both sides). If
$\widebar{\alpha}'\in\autf(P{\cdot}C_S(P))$ is any other extension
of $\alpha$, then by \cite[Lemma 3.8]{bcglo1}, there is some
$g\in{}C_S(P)$ such that
$\widebar{\alpha}'=c_g\circ\widebar{\alpha}$. Thus, by
Definition \ref{D:J,j} again, $\Theta(\widebar{\alpha}')  =
\Theta(c_g\circ\widebar{\alpha}) =
\theta(g)\Theta(\widebar{\alpha})$,  and
	\[\Theta(\widebar{\alpha}')\cdot\theta(C_S(P)) =
	\theta(g)\Theta(\widebar{\alpha})\cdot\theta(C_S(P)) =
	\Theta(\widebar{\alpha})\theta(\widebar{\alpha}(g))
	\cdot\theta(C_S(P)) =
	\Theta(\widebar{\alpha})\cdot\theta(C_S(P)),\]
and so the definition of $\Theta_P(\alpha)$ is independent of the choice 
of the extension $\widebar{\alpha}$. This shows that $\Theta_P$ is well 
defined.

Notice also that $\Theta_P$ clearly respects compositions, and
since $\Theta_P(\alpha) = x\cdot\theta(C_S(P)) =
\theta(C_S(P))\cdot x$, for some $x\in\Gamma$,  we conclude that
$x\in N_\Gamma(\theta(C_S(P)))$. Thus  $\Theta_P$ induces a
homomorphism
\[\Theta_P\colon \Aut_\calf(P)\rTo
N_\Gamma(\theta(C_S(P)))/\theta(C_S(P)).\] We now make use of
Lemma \ref{ext.F}, which gives sufficient conditions to the
existence of an extension of a fusion mapping triple.

If $\alpha\in\autf(P)$, and $x\in\Theta_P(\alpha)$, then
$x=y\cdot\theta(h)$ for some $h\in C_S(P)$ and
$y\in\Theta(\widebar{\alpha})$, where $\widebar{\alpha}$ is an
extension of $\alpha$ to $P\cdot C_S(P)$. Hence for any $g\in{}P$,
	\[x\theta(g)x^{-1}=y\cdot\theta(hgh^{-1})y^{-1} = 
	y\theta(g)y^{-1} = \theta(\widebar{\alpha}(g)) =
	\theta(\alpha(g)).\] 
This shows that point (+) of Lemma \ref{ext.F} holds, and so it remains to 
check ($*$).

Assume $P\lneqq{}Q\le{}S$, $P\nsg{}Q$, and let $\alpha\in\autf(P)$
and $\beta\in\autf(Q)$ be such that $\alpha=\beta|_P$.  Then
$Q{\cdot}C_S(P)\le{}N_\alpha$ in the terminology of axiom (II), so
$\alpha$ extends to another automorphism
$\gamma\in\autf(Q{\cdot}C_S(P))$, and
$\Theta_P(\alpha)=\Theta(\gamma){\cdot}\theta(C_S(P))$ by
definition of $\Theta_P$.  By \cite[Lemma 3.8]{bcglo1} again,
$\gamma|_Q=c_g\circ\beta$ for some $g\in{}C_S(P)$. Hence, by
Definition  \ref{D:J,j}, $\Theta(\gamma) = \Theta(c_g\circ\beta) =
\theta(g)\cdot\Theta(\beta)$, and so
	\[\Theta_P(\alpha) = 
	\theta(g)\cdot\Theta(\beta)\cdot\theta(C_S(P)) =
	\Theta(\beta)\theta(\beta(g))\cdot\theta(C_S(P)) =
	\Theta(\beta)\cdot\theta(C_S(P)).\] 
In particular, $\Theta_P(\alpha)\supseteq\Theta(\beta)$.  This shows that 
point ($*$) of Lemma \ref{ext.F} is satisfied as well, and thus, by the 
lemma, $\Theta$ can be extended to a fusion mapping triple on 
$\calf^\calh$.
\end{proof}

Recall that a fusion subsystem of index prime to $p$ in a
saturated fusion system $\calf$ over $S$ is a saturated subsystem
$\calf_0\subseteq\calf$ over the same $p$-group $S$, such that
$\Aut_{\calf_0}(P)\ge{}O^{p'}(\autf(P))$ for all $P\le{}S$.
Equivalently, $\calf_0\subseteq\calf$ has index prime to $p$ if
and only if it is saturated and contains the subcategory $\oppf$
of Definition \ref{D:oppf:opf}.  We are now ready to prove our
main result about these subsystems.

\begin{Thm}  \label{F:p'-quot.}
For any saturated fusion system $\calf$ over a $p$-group $S$,
there is a bijective correspondence between subgroups
    \[ H\le \gppf = \outf(S)/\outf^0(S), \]
and saturated fusion subsystems $\calf_H$ of $\calf$ over $S$ of
index prime to $p$ in $\calf$.  The correspondence is given by
associating to $H$ the fusion system generated by
$\widehat{\theta}^{-1}(\calb(H))$, where $\widehat{\theta}$ is the
functor of Proposition \ref{Theta:p'-quot.}.
\end{Thm}

\begin{proof}  Let $\calf_0\subseteq\calf$ be any saturated fusion 
subsystem over $S$ which contains $\oppf$. Then $\Out_\calf^0(S)\nsg 
\Out_{\calf_0}(S)$, and one can set $H=\Out_{\calf_0}(S)/\Out_\calf^0(S)$. 
We first show that a morphism $\varphi$ of $\calf^c$ is in $\calf_0$ if 
and only if $\widehat{\theta}(\varphi)\in H$. Clearly it suffices to prove 
this for isomorphisms in $\calf^c$. 

Let $P,Q\le S$ be $\calf$-centric, $\calf$-conjugate subgroups, and fix an 
isomorphism $\varphi\in\Iso_\calf(P,Q)$. By Lemma \ref{F=F'.N_F(S)}, we 
can write $\varphi=\psi\circ(\alpha|_P)$, where $\alpha\in\autf(S)$ and 
$\psi\in\Iso_{\oppf}(\alpha(P),Q)$. Then $\varphi$ is in $\calf_0$ if and 
only if $\alpha|_P$ is in $\calf_0$.  Also, by definition of 
$\widehat{\theta}$ (and of $H$), $\widehat{\theta}(\varphi)\in{}H$ if and 
only if $\alpha\in\Aut_{\calf_0}(S)$.  So it remains to prove that 
$\alpha|_P\in\Mor(\calf_0)$ if and only if $\alpha\in\Aut_{\calf_0}(S)$.

If $\alpha\in\Aut_{\calf_0}(S)$, then $\alpha|_P$ is also in $\calf_0$ by 
definition of a fusion system.  So it remains to prove the converse.  
Assume $\alpha|_P$ is in $\calf_0$.  The same argument as that used to 
prove Proposition \ref{F:p-solv.quot.}(c) shows that $\alpha(P)$ is 
$\calf_0$-centric, and hence fully centralized in $\calf_0$.  Since 
$\alpha|_P$ extends to an (abstract) automorphism of $S$, axiom (II) 
implies that it extends to some $\alpha_1\in\Hom_{\calf_0}(N_S(P),S)$.  
Since $P$ is $\calf$-centric, \cite[Proposition A.8]{BLO2} applies to show 
that $\alpha_1=(\alpha|_{N_S(P)})\circ{}c_g$ for some $g\in{}Z(P)$, and 
thus that $\alpha|_{N_S(P)}\in\Hom_{\calf_0}(N_S(P),S)$.  Also, 
$N_S(P)\gneqq{}P$ whenever $P\lneqq{}S$, and so we can continue this 
process to show that $\alpha\in\Aut_{\calf_0}(S)$.  This finishes the proof 
that $\calf_0=\widehat{\theta}^{-1}(H)$. 

Now fix a subgroup $H\le\outf(S)/\outf^0(S)$, and let $\calf_H$ be the 
smallest fusion system over $S$ which contains 
$\widehat{\theta}^{-1}(\calb(H))$.  We must show that $\calf_H$ is a 
saturated fusion subsystem of index prime to $p$ in $\calf$.  For 
$\calf$-centric subgroups $P,Q\le{}S$, $\Hom_{\calf_H}(P,Q)$ is the set of 
all morphisms $\varphi$ in $\homf(P,Q)$ such that 
$\widehat{\theta}(\varphi)\in{}H$.  In particular, 
$\calf_H\supseteq\oppf$, since morphisms in $\oppf$ are sent by 
$\widehat{\theta}$ to the identity element.

Define $\Theta\:\Mor(\calfc)\rTo\sset(\gppf)$ by setting
$\Theta(\varphi)=\{\widehat{\theta}(\varphi)\}$; i.e., the image
consists of subsets with one element.  Let
$\theta\in\Hom(S,\gppf)$ be the trivial (and unique) homomorphism.
Then $(\gppf,\theta,\Theta)$ is a fusion mapping triple on
$\calfc$.  By Lemma \ref{extend:Fc-Fq}, this can be extended to a
fusion mapping triple on $\calfq$; and hence $\calf_H$ is
saturated by Proposition \ref{F:p-solv.quot.}.

By Theorem \ref{Alp.fusion}(a) (Alperin's fusion theorem),
$\calf_H$ is the unique saturated fusion subsystem of $\calf$ with
the property that a morphism $\varphi\in\homf(P,Q)$ between
$\calf$-centric subgroups of $S$ lies in $\calf_H$ if and only if
$\widehat{\theta}(\varphi)\in{}H$. This shows that the
correspondence is indeed bijective.
\end{proof}

The next theorem describes the relationship between subgroups of
index prime to $p$ in a $p$-local finite group $\SFL$ and certain
covering spaces of $|\call|$.

\begin{Thm} \label{L:p'-quot.}
Fix a $p$-local finite group $\SFL$.  Then for each subgroup
$H\le\outf(S)$ containing $\outf^0(S)$, there is a unique
$p$-local finite subgroup $(S,\calf_H,\call_H)$, such that
$\calf_H$ has index prime to $p$ in $\calf$,
$\Out_{\calf_H}(S)=H$, and  $\call_H=\pi^{-1}(\calf_H)$ (where
$\pi$ is the usual functor from $\call$ to $\calf$). Furthermore,
$|\call_H|$ is homotopy equivalent, via its inclusion into
$|\call|$, to the covering space of $|\call|$ with fundamental
group $\widetilde{H}$:  where $\widetilde{H}\le\pi_1(|\call|)$ is
the subgroup such that
$\bar\theta(\widetilde{H}/\oppg{\pi_1(|\call|)})$ corresponds to
$H/\outf^0(S)$ under the isomorphism
    $$ \pi_1(|\call|)/\oppg{\pi_1(|\call|)} \cong
    \pi_1(|\calfc|)/\oppg{\pi_1(|\calfc|)}
    \RIGHT5{\widebar{\theta}}{\cong} \outf(S)/\outf^0(S) $$
of Proposition \ref{Theta:p'-quot.}.
\end{Thm}

\begin{proof}  By Theorem \ref{L:p-p'-quot.}, applied to the composite
functor
    $$ \call \Right5{\pi} \calfc \Right5{\widehat{\theta}}
    \calb(\outf(S)/\outf^0(S)), $$
$(S,\calf_H,\call_H)$ is a $p$-local finite group, and $|\call_H|$
is homotopy equivalent to the covering space of $|\call|$ with
fundamental group $\widetilde{H}$ defined above.  The uniqueness
follows from Theorem \ref{F:p'-quot.}.
\end{proof}

Theorem \ref{F:p'-quot.} shows, in particular, that any saturated
fusion system $\calf$ over $S$ contains a unique minimal saturated
subsystem $\Oppf$ of index prime to $p$:  the subsystem
$\calf_0\subseteq\calf$ with $\Out_{\calf_0}(S)=\outf^0(S)$.
Furthermore, if $\calf$ has an associated centric linking system
$\call$, then Theorem \ref{L:p'-quot.} shows that $\Oppf$ has an
associated linking system $\Oppl$, whose geometric realization
$|\Oppl|$ is homotopy equivalent to the covering space of
$|\call|$ with fundamental group $\oppg{\pi_1(|\call|)}$.

\bigskip


\newsub{Extensions of index prime to $p$}
It remains to consider the opposite problem:  describing the extensions of 
a given saturated fusion system of index prime to $p$.  As before, for a 
saturated fusion system $\calf$ over a $p$-group $S$, $\Aut\fus(S,\calf)$ 
denotes the group of fusion preserving automorphisms of $S$ (Definition 
\ref{D:autfus}).   Theorem \ref{p'-extn} below states that each subgroup 
of $\Out\fus(S,\calf)=\Aut\fus(S,\calf)/\autf(S)$ of order prime to $p$ 
gives rise to an extension of $\calf$.

The following lemma will be needed to compare the obstructions to the
existence and uniqueness of linking systems, in a fusion system
and in a fusion subsystem of index prime to $p$.  Recall the definition of 
the orbit category of a fusion system in Definition \ref{D:orbit}. 

\begin{Lem}  \label{H<|G}
Fix a saturated fusion system $\calf$ over a $p$-group $S$, and
let $\calf'\subseteq\calf$ be a saturated subsystem of index prime
to $p$. Assume $\Out_{\calf'}(S)\nsg\outf(S)$, and set
$\pi=\outf(S)/\Out_{\calf'}(S)$.  Then for any
$F:\orb^c(\calf)\rTo\zploc\mod$, there is a natural action of
$\pi$ on the higher limits of $F|_{\orb^c(\calf')}$, and
    $$ \higherlim{\orb^c(\calf)}*(F) \cong
    \Bigl[\higherlim{\orb^c(\calf')}*\bigl(F|_{\orb^c(\calf')}\bigr)
    \Bigr]^{\pi}. $$
\end{Lem}

\begin{proof}  Let $\orb^c(\calf)\mod$ denote the category of functors
$\orb^c(\calf)\op\rTo\Ab$.  For any $F$ in $\orb^c(\calf)\mod$,
$\outf(S)$ acts on $\prod_{P\in\Ob(\orb^c(\calf))}F(P)$ by letting
$\alpha\in\outf(S)$ send $F(\alpha(P))$ to $F(P)$ via the induced
map $\alpha^*$.  This restricts to an action of $\pi$ on
$\lambda(F)\defeq\invlim{\orb^c(\calf')}(F|_{\orb^c(\calf')})$;
and by definition of inverse limits,
    $$ \invlim{\orb^c(\calf)}(F) \cong
    \Bigl[\invlim{\orb^c(\calf')}(F|_{\orb^c(\calf')})\Bigr]^{\pi}
    = \bigl[\lambda(F)\bigr]^{\pi}. $$

Since $\calf'$ is a subsystem of index prime to $p$,
$\Out_{\calf'}(S)$ is a subgroup of $\Out_\calf^0(S)$. Hence, by a
slight abuse of notation, one has a functor
$\widehat{\theta}\:\calf\rTo\calb(\pi)$ given as the composite of
the functor $\widehat{\theta}$ of Proposition \ref{hat{theta}}
with projection to $\pi$.  For any $\Z[\pi]$-module $M$, regarded
as a functor on $\calb(\pi)$, we let $\beta(M)$ denote the
composite functor $M\circ\widehat{\theta}$.  Then
$\Z[\pi]\mod\rTo^{\beta}\orb^c(\calf)\mod$ is an exact functor,
and a left adjoint to
$\orb^c(\calf)\mod\rTo^{\lambda}\Z[\pi]\mod$.  In particular, the
existence of a left adjoint shows that $\lambda$ sends injective
objects in $\orb^c(\calf)\mod$ to injective $\Z[\pi]$-modules.

Thus, if $0\to{}F\to{}I_0\to{}I_1\to\cdots$ is an injective
resolution of $F$ in $\orb^c(\calf)\mod$, and $F$ takes values in
$\zploc$-modules, then
    \begin{align*}
    \higherlim{\orb^c(\calf)}*(F) &\cong
    H^*\bigl(0\to\lambda(I_0)_{(p)}^{\pi}\rTo
    \lambda(I_1)_{(p)}^{\pi} \rTo\cdots\bigr) \\
    &\cong \bigl[H^*\bigl(0\rTo\lambda(I_0)_{(p)}\to
    \lambda(I_1)_{(p)}\rTo\cdots\bigr)\bigr]^{\pi} \cong
    \Bigl[\higherlim{\orb^c(\calf')}*(F|_{\orb^c(\calf')})\Bigr]^{\pi}
    \,. \qedhere
    \end{align*}
\end{proof}

We are now ready to examine extensions of index prime to $p$.

\begin{Thm} \label{p'-extn}
Fix a saturated fusion system $\calf$ over a $p$-group $S$.  Let
    $$ \pi\le \Out\fus(S,\calf)\defeq\Aut\fus(S,\calf)/\autf(S) $$
be any subgroup of order prime to $p$, and let $\widetilde\pi$
denote the inverse image of $\pi$ in $\Aut\fus(S,\calf)$.
\begin{enumerate}
\item Let $\calf.\pi$ be the fusion system over $S$ generated (as
a category) by $\calf$ together with restrictions of automorphisms
in $\widetilde{\pi}$.  Then $\calf.\pi$ is a saturated fusion
system, which contains $\calf$ as a fusion subsystem of index
prime to $p$.

\item If $\calf$ has an associated centric linking system, then so does 
$\calf.\pi$.

\item Let $\call$ be a centric linking system associated to $\calf$.  
Assume that for each $\alpha\in\pi$, the action of $\alpha$ on $\calf$ 
lifts to an action on $\call$.  Then there is a unique centric linking 
system $\call.\pi$ associated to $\calf.\pi$ whose restriction to $\calf$ 
is $\call$.
\end{enumerate}
\end{Thm}

\begin{proof}  \noindent\textbf{(a) } By definition, every morphism in
$\calf.\pi$ is the composite of morphisms in $\calf$ and
restrictions of automorphisms of $S$ which normalize $\calf$. If
$\psi\in\hom_\calf(P,Q)$ and $\varphi\in\Aut\fus(S,\calf)$,  then
one has
	\beq \varphi|_{Q}\circ\psi =
	\varphi|_{Q}\circ\psi\circ\varphi^{-1}|_{\varphi(P)}\circ\varphi|_{P}
	= \psi'\circ\varphi|_{P}, \tag{1} \eeq
where $\psi'\in\homf(\varphi(P),\varphi(Q))$ (since $\varphi$ is fusion 
preserving). Hence each morphism in $\calf.\pi$ is the composite of
the restriction of a morphism in $\widetilde{\pi}$ followed by a
morphism in $\calf$. 

We next claim $\calf$ is a fusion subsystem of $\calf.\pi$ of index prime 
to $p$.  More precisely, we will show, for all $P\le{}S$, that 
$\Aut_{\calf}(P)\nsg\Aut_{\calf.\pi}(P)$ with index prime to $p$.  To see 
this, let $\pi_1\subseteq\widetilde{\pi}$ be the set of automorphisms 
$\varphi$ in $\widetilde{\pi}\le\Aut\fus(S,\calf)$ such that $\varphi(P)$ 
is $\calf$-conjugate to $P$, and let $\pi_0\subseteq\pi_1$ be the set of 
classes $\varphi$ such that $\varphi|_P\in\homf(P,S)$.  If $\varphi(P)$ 
and $\psi(P)$ are both $\calf$-conjugate to $P$, then $\psi(\varphi(P))$ 
is $\calf$-conjugate to $\psi(P)$ since $\psi$ is fusion preserving, and 
thus is $\calf$-conjugate to $P$.  This shows that $\pi_1$ is a subgroup 
of $\widetilde{\pi}$, and an argument using (1) shows that $\pi_0$ is also 
a subgroup.  By definition, $\pi_0\ge\autf(S)$, so $\pi_1/\pi_0$ has 
order prime to $p$ since $\pi$ does.  Using (1), define
	$$ \theta_P\: \Aut_{\calf.\pi}(P) \Right5{} \pi_1/\pi_0 $$
by setting $\theta_P(\alpha)=[\psi]$ if 
$\alpha=\psi|_{\beta(P)}\circ\beta$ for some $\beta\in\homf(P,S)$ and some 
$\psi{\cdot}\autf(S)\in\pi$.  By definition of $\pi_0$, this is well 
defined, and $\Ker(\theta_P)=\autf(P)$.  Since $\pi_1/\pi_0$ has order 
prime to $p$, this shows that $\autf(P)\nsg\Aut_{\calf.\pi}(P)$ with index 
prime to $p$.

We next claim that $\calf$ and $\calf.\pi$ have the same fully 
centralized, fully normalized and centric subgroups (compare with the 
proof of Proposition \ref{F:p-solv.quot.}(c)).  By definition, $\calf$ and 
$\calf.\pi$ are fusion systems over the same $p$-group $S$. Since each 
$\calf$-conjugacy class is contained in some $\calf.\pi$-conjugacy class, 
any subgroup $P\le S$ which is fully centralized (fully normalized) in 
$\calf.\pi$ is fully centralized (fully normalized) in $\calf$. By the 
same argument, any $\calf.\pi$-centric subgroup is also $\calf$-centric.

Conversely, assume that $P$ is not fully centralized in
$\calf.\pi$, and let $P'\le S$ be a subgroup $\calf.\pi$-conjugate
to $P$ such that $|C_S(P')|>|C_S(P)|$. Let $\psi\colon P'\rTo P$
be an $\calf.\pi$-isomorphism between them. Then, by the argument
above $\psi = \psi'\circ\varphi$, where $\psi'$ is a morphism in
$\calf$ and $\varphi\in\widetilde{\pi}$ is the restriction to $P'$
of an automorphism of $S$. Hence
$|C_S(\varphi(P'))|=|C_S(P')|>|C_S(P)|$, and since $\varphi(P')$ is
$\calf$-conjugate to $P$, this shows that $P$ is not fully
centralized in $\calf$. A similar argument shows that if $P$ is
not fully normalized (or not centric) in $\calf.\pi$, then it is
not fully normalized (or not centric) in $\calf$.

We prove that $\calf.\pi$ is saturated using Theorem
\ref{centr->sat}(b). Thus, we must show that conditions (I) and
(II) of Definition \ref{sat.Frob.} are satisfied for all
$\calf.\pi$-centric subgroups, and that $\calf.\pi$ is generated
by restrictions of morphisms between its centric subgroups. Let
$P\le S$ be a subgroup which is fully normalized in $\calf.\pi$.
Then it is fully normalized in $\calf$, and since $\calf$ is
saturated, it is fully centralized there and
$\Aut_S(P)\in\sylp{\autf(P)} = \sylp{\Aut_{\calf.\pi}}$. Hence,
condition (I) holds for any $P$ in $(\calf.\pi)^c$.

Let $\psi\colon P\rTo Q$ be a morphism in $\calf.\pi$, and write
$\psi = \psi'\circ\varphi$, as before. Set
	\[N_\psi = \{g\in N_S(P)\,|\,\psi\circ c_g\circ\psi^{-1} \in
	\Aut_S(\psi(P))\}; \] 
then $\varphi(N_\psi) = N_{\psi'}$ since 
$\varphi{}c_g\varphi^{-1}=c_{\varphi(g)}$ for all $g\in{}S$. Since 
condition (II) holds for $\calf$, the morphism $\psi'$  can be extended to 
$N_{\psi'}$. Hence $\psi=\psi'\circ\varphi$ can be extended to $N_\psi = 
\varphi(N_{\psi'})$, and so condition (II) holds for $\calf.\pi$.

That all morphisms in $\calf.\pi$ are composites of restrictions of
$\calf.\pi$-morphisms between $\calf.\pi$-centric subgroups holds by
construction. Thus Theorem \ref{centr->sat}(b)applies, and
$\calf.\pi$ is saturated.

\smallskip

\noindent\textbf{(b) } Let $\calz_\calf\:\orb^c(\calf)\rTo\Ab$ be
the functor $\calz_\calf(P)=Z(P)$, and similarly for
$\calz_{\calf.\pi}$ (see Definition \ref{D:orbit}).  By Lemma \ref{H<|G},
restriction of categories induces a monomorphism
    \beq \higherlim{\orb^c(\calf.\pi)}*(\calz_{\calf.\pi})
    \Right5{} \higherlim{\orb^c(\calf)}*(\calz_{\calf}) \tag{2} \eeq
whose image is the subgroup of $\pi$-invariant elements.

By \cite[Proposition 3.1]{BLO2}, the obstruction $\eta(\calf.\pi)$ to the 
existence of a centric linking system associated to $\calf.\pi$ lies in 
$\higherlim{}3(\calz_{\calf.\pi})$.  From the construction in 
\cite{BLO2} of these obstructions (and the fact that $\calf.\pi$ and 
$\calf$ have the same centric subgroups), it is clear that the restriction 
map (2) sends $\eta(\calf.\pi)$ to $\eta(\calf)$.  So if there is a 
linking system $\call$ associated to $\calf$, then $\eta(\calf)=0$, so 
$\eta(\calf.\pi)=0$ by Lemma \ref{H<|G}, and there is a linking system 
$\call.\pi$ associated to $\calf.\pi$.

\smallskip

\noindent\textbf{(c) } Let $\call.\pi$ be a centric linking system 
associated to $\calf.\pi$, as constructed in (b), and let $\call'$ be its 
restriction to $\calf$ (i.e., the inverse image of $\calf$ under the the 
projection $\call.\pi\rTo(\calf.\pi)^c$).  By \cite[Proposition 3.1]{BLO2} 
again, the group $\higherlim{}2(\calz_\calf)$ acts freely and transitively 
on the set of all centric linking systems associated to $\calf$.  If the 
action on $\calf$ of each $\alpha\in\pi$ lifts to an action on $\call$, 
then the element of $\higherlim{}2(\calz_\calf)$ which measures the 
difference between $\call$ and $\call'$ is $\pi$-invariant, and hence (by 
Lemma \ref{H<|G} again) is the restriction of an element of 
$\higherlim{}2(\calz_{\calf.\pi})$.  Upon modifying the equivalence class 
of $\call.\pi$ by this element, if necessary, we get a centric linking 
system associated to $\calf.\pi$ whose restriction to $\calf$ is $\call$.
\end{proof}

As was done in the last section for extensions with $p$-group quotient, we 
now translate this last result to a theorem stated in terms of fibration 
sequences.

\begin{Thm} \label{|L|->X->Bpi'}
Fix a $p$-local finite group $\SFL$, a finite group $\Gamma$ of order 
prime to $p$, and a fibration $E\Right2{v}B\Gamma$ with fiber 
$X\simeq|\call|\pcom$.  Then there is a $p$-local finite group 
$(S,\calf',\call')$ such that $\calf\subseteq\calf'$ is normal of index 
prime to $p$, $\Aut_{\calf'}(S)/\Aut_{\calf}(S)$ is a quotient group of 
$\Gamma$, and $E\pcom\simeq|\call'|\pcom$.
\end{Thm}

\begin{proof} For any space $Y$, let $\Aut(Y)$ denote the topological 
monoid of homotopy equivalences $Y\Right2{\simeq}Y$.  Fibrations with 
fiber $Y$ and base $B$ are classified by homotopy classes of maps $B\rTo 
B\Aut(Y)$.  This follows, for example, as a special case of the main 
theorem in \cite{DKS}. 

We are thus interested  in the classifying space $B\Aut(|\call|\pcom)$, 
whose homotopy groups were determined in \cite[\S8]{BLO2}. To describe 
these, let $\Aut\isotyp(\call)$ be the monoid of \emph{isotypical} self 
equivalences of the category $\call$; i.e., the monoid of all equivalences 
of categories $\psi\in\Aut(\call)$ such that for all $P\in\Ob(\call)$, 
$\psi_{P,P}(\Im(\delta_P))=\Im(\delta_{\psi(P)})$.  Let 
$\Out\isotyp(\call)$ be the group of all isotypical self equivalences 
modulo natural isomorphisms of functors.  By \cite[Theorem 8.1]{BLO2}, 
$\pi_i(B\Aut(|\call|\pcom))$ is a finite $p$-group for $i=2$ and vanishes 
for $i>2$, and 
        $$ \pi_1(B\Aut(|\call|\pcom)) \cong \Out\isotyp(\call). $$

Each isotypical self equivalence of $\call$ is naturally isomorphic to one 
which sends inclusions to inclusions (this was shown in \cite[Lemma 
5.1]{BLO1} for linking systems of a group, and the general case follows by 
the same argument).  Thus each element of $\Out\isotyp(\call)$ is 
represented by some $\beta$ which sends inclusions to inclusions.  This in 
turn implies that $\beta_S$ --- the restriction of $\beta_{S,S}$ to 
$\delta_S(S)\le\Aut_\call(S)$ --- lies in $\Aut\fus(S,\calf)$, and that 
for every $\calf$-centric subgroup $P$ of $S$, the functor $\beta$ sends 
$P$ to the subgroup $\beta_S(P)\le{}S$.  In particular, $\beta$ is an 
automorphism of $\call$, since it induces a bijection on the set of 
objects.  

Now let $B\Gamma\rTo^{f_v}B\Aut(|\call|\pcom)$ be the map which classifies 
the fibration $E\Right2{v}B\Gamma$.  Since $\Gamma$ is a $p'$-group, the 
fibration $E\Right2{v}B\Gamma$ is uniquely determined by the action of 
$\Gamma$ on $\nv{\call}\pcom$; more precisely, by the map induced by $f_v$ 
on fundamental groups: 
	$$\pi_1(f_v)\:\Gamma  \Right4{}  \Out\isotyp(\call)\,.$$

By an argument identical to that used to prove \cite[Theorem 6.2]{BLO1}, 
there is an exact sequence
        $$ 0 \Right3{} \higherlim{\orb^c(\calf)}1(\calz_\calf) \Right4{} 
        \Out\isotyp(\call) \Right4{\mu_\call} \Out\fus(S,\calf), $$
where $\mu_\call$ is defined by restricting a functor $\call\to\call$ to 
$\delta_S(S)\le\Aut_\call(S)$.  Also, $\orb^c(\calf)$ and $\calz_\calf$ are 
as in Definition \ref{D:orbit}, but all that we need to know 
here is that $\Ker(\mu_\call)$ is a $p$-group.  Set
        $$ \psi_v=\mu_\call\circ\pi_1(f_v)\: \Gamma \Right4{} 
        \Out\fus(S,\calf). $$

Set $\pi=\Im(\psi_v)\le\Out\fus(S,\calf)$. By its definition, $\psi_v$ 
comes equipped with a lift to $\Out\isotyp(\call)$.  So by the above 
remarks, every element of $\pi$ lifts to an automorphism of $\call$.  Let 
$(S,\calf.\pi,\call.\pi)$ be the $p$-local finite group constructed in 
Theorem \ref{p'-extn}(a,c) as an extension of $\SFL$.  This induces a 
fibration
        $$ |\call| \Right4{} |\call.\pi| \Right4{} B\pi. $$
By \cite[Corollary I.8.3]{BK}, there is a fiberwise completion $E_w$ of 
$|\call.\pi|$ which sits in a fibration sequence 
        $$ |\call|\pcom \Right4{} E_w \Right4{w} B\pi, $$
and this fibration is classified by a map 
        $$ f_w\: B\pi \Right4{} B\Aut(|\call|\pcom). $$
As was the case for $f_v$, the induced map between fundamental groups 
$\pi_1(f_w)$ determines the fibration. 

Consider the diagram
	\beq \begin{diagram}[w=40pt] 
	\Gamma & \rTo^{\rho} & \pi \\
	\dTo<{\pi_1(f_v)} & \ldTo<{\pi_1(f_w)} & \dTo>{\incl} \\
	\Out\isotyp(\call) & \rTo^{\mu_\call} & \Out\fus(S,\calf) \rlap{\,,}
	\end{diagram} \tag{1} \eeq
where $\rho\:\Gamma\to\pi$ is the restriction of $\psi_v$ to its image, 
and where the square commutes since each composite is equal to $\psi_v$.   
The composite $\mu_\call\circ\pi_1(f_w) \in\Hom(\pi,\Out\fus(S,\calf))$ is 
the homomorphism induced by the homotopy action of $\pi=\pi_1(B\pi)$ on 
the fiber $|\call|\pcom$.  By the construction in Theorem \ref{p'-extn}, 
this is just the inclusion map.  Hence the lower right triangle in the 
above diagram commutes. Since $\Gamma$ has order prime to $p$ and 
$\Ker(\mu_\call)\cong\higherlim{}1(\calz_\calf)$ is a $p$-group, any 
homomorphism from $\Gamma$ to $\Out\fus(S,\calf)$ has a unique conjugacy 
class of lifting to $\Out\isotyp(\call)$ by the Schur-Zassenhaus theorem 
(cf. \cite[Theorem 6.2.1]{Gorenstein}).  In particular, 
$\pi_1(f_w)\circ\rho$ and $\pi_1(f_v)$ are conjugate as homomorphisms from 
$\Gamma$ to $\Out\isotyp(\call)$.  Thus, the upper left triangle in (1) 
commutes up to conjugacy. 

Since $\pi_2(B\Aut(|\call|\pcom))$ is a finite $p$-group and the higher 
homotopy groups all vanish, we have shown that $f_w\circ{}B\rho$ is 
homotopic to $f_v$. In other words, we have a map of fibrations
	\begin{diagram}[w=30pt] 
        |\call|\pcom  &  \rTo  &  E  & \rTo{v} &  B\Gamma  \\
        \dIgual && \dTo{\lambda} && \dTo>{B\rho} \\
        |\call|\pcom  &  \rTo  &  E_w & \rTo{w} &  B\pi, 
	\end{diagram}
A comparison of the spectral sequences for these two fibrations shows that 
$\lambda$ is an $\F_p$-homology equivalence.  Since $|\call|$ is $p$-good 
\cite[Proposition 1.12]{BLO2}, the natural map from $|\call.\pi|$ to its 
fiberwise completion $E_w$ is also an $\F_p$-homology equivalence.  
Hence these maps induce homotopy equivalences 
$E\pcom\simeq(E_w)\pcom\simeq|\call.\pi|\pcom$.  This finishes the proof 
of the theorem, with $\calf'=\calf.\pi$ and $\call'=\call.\pi$.
\end{proof}

Note that we are not assuming that $\Gamma$ injects into 
$\Out\fus(S,\calf)$ in the hypotheses of Theorem \ref{|L|->X->Bpi'}.  
Thus, for example, if we start with a product fibration 
$E=|\call|\pcom\times{}B\Gamma$, then we end up with 
$(S,\calf',\call')=\SFL$ (and $E\pcom\simeq|\call|\pcom$).  
\mynote{C:  do we need this remark? I put the last diagram to make it 
clear what is the situation in which $\Gamma$ in not isomorphic to $\pi$, 
for instance with $\pi=1$, the fibration in the top row is trivial.}  
\mynote{B:  I don't have a strong opinion, but I think it's good to 
emphasize this.}

We do not know whether the index prime to $p$ analog of Corollary 
\ref{exotic-ext} holds. The key point in the proof of Corollary 
\ref{exotic-ext} is the observation that if $G$ is a finite group, whose 
fusion system at $p$ contains a normal subsystem of index $p^m$, then $G$ 
contains a normal subgroup of index $p^m$. This is not true for subsystems 
of index prime to $p$. For example, the fusion system at the prime $3$ of 
the simple groups $J_4$ and $Ru$ has a subsystem of index 2.


\newsect{Central extensions of fusion systems and linking systems} 

Let $\calf$ be a fusion system over a $p$-group $S$. We say that a 
subgroup $A\leq S$ is \emph{central} in $\calf$ if $C_{\calf}(A)=\calf$, 
where $C_\calf(A)$ is the centralizer fusion system defined in 
\cite[Definition A.3]{BLO2}.  Thus $A$ is central in $\calf$ if 
$A\le{}Z(S)$, and each morphism $\varphi\in\homf(P,Q)$ in $\calf$ extends 
to a morphism $\widebar{\varphi}\in\homf(PA,QA)$ such that 
$\widebar{\varphi}|_A=\Id_A$.  

In this section, we first study quotients of fusion systems, and of 
$p$-local finite groups, by central subgroups.  Afterwards, we will invert 
this procedure, and study central extensions of fusion systems and 
$p$-local finite groups.  

\bigskip

\newsub{Central quotients of fusion and linking systems} 
We first note that every saturated fusion system $\calf$ contains a unique 
maximal central subgroup, which we regard as the \emph{center} of $\calf$.

\begin{Prop} \label{ZF(S)}
For any saturated fusion system $\calf$ over a $p$-group $S$, define
	$$ Z_\calf(S) = \bigl\{ x\in{}Z(S) \,\big|\, \varphi(x)=x,
	\text{ all } \varphi\in\Mor(\calfc) \bigr\}
	= \higherlim{\calfc}{}Z(-): $$
the inverse limit of the centers of all $\calf$-centric subgroups of $S$.  
Then $Z_\calf(S)$ is the \emph{center} of $\calf$:  it is central in 
$\calf$, and contains all other central subgroups. 
\end{Prop}

\begin{proof}  By definition, if $A$ is central in $\calf$, then 
$A\le{}Z(S)$, and any morphism in $\calf$ between subgroups containing $A$ 
restricts to the identity on $A$.  Since all $\calf$-centric subgroups 
contain $Z(S)$, this shows that $A\le{}Z_\calf(S)$.

By Alperin's fusion theorem (Theorem \ref{Alp.fusion}(a)), each morphism in 
$\calf$ is a composite of restrictions of morphisms between $\calf$-centric 
subgroups.  In particular, each morphism is a restriction of a morphism 
between subgroups containing $Z_\calf(S)$ which is the identity on 
$Z_\calf(S)$, and thus $Z_\calf(S)$ is central in $\calf$.
\end{proof}

The center of a fusion system $\calf$ has already appeared when studying 
mapping spaces of classifying spaces associated to $\calf$.  By 
\cite[Theorem 8.1]{BLO2}, for any $p$-local finite group $\SFL$, 
	$$ \map(|\call|\pcom,|\call|\pcom)_{\Id} \simeq BZ_\calf(S). $$

We next define the quotient of a fusion system by a central subgroup.

\begin{Defi} \label{D:F/A}
Let $\calf$ be a fusion system over a $p$-group $S$, and let $A$ be a 
central subgroup. Define $\calf/A$ to be the fusion system over $S/A$ with 
morphism sets
	$$ \Hom_{\calf/A}(P/A,Q/A) =
	\Im\bigl[\Hom_\calf(P,Q) \Right4{} \Hom(P/A,Q/A)\bigr]. $$
\end{Defi}

By \cite[Lemma 5.6]{BLO2}, if $\calf$ is a saturated fusion system over a 
$p$-group $S$, and $A\le{}Z_\calf(S)$, then $\calf/A$ is also saturated as 
a fusion system over $S/A$.  We now want to study the opposite question:  
if $\calf/A$ is saturated, is $\calf$ also saturated?  The following 
example shows that it is very easy to construct counterexamples to this.  
In fact, under a mild hypothesis on $S$, if a fusion system over $S/A$ has 
the form $\calf/A$ for any fusion system $\calf$ over $S$, then it has 
that form for a fusion system $\calf$ which is not saturated.

\begin{Ex}
Fix a $p$-group $S$ and a central subgroup $A\le{}Z(S)$.  Assume there is 
$\alpha\in\Aut(S){\sminus}\Inn(S)$ such that $\alpha|_A=\Id_A$ and 
$\alpha$ induces the identity on $S/A$.  Let $\widebar{\calf}$ be a 
saturated fusion system over $S/A$, and assume there is some fusion system 
$\calf_0$ over $S$ such that $\calf_0/A=\widebar{\calf}$.  Let 
$\calf\supseteq\calf_0$ be the fusion system over $S$ defined by setting 
	\begin{multline*}  
	\Hom_{\calf}(P,Q) = \bigl\{ 
	\varphi\in\Hom(P,Q) \,\big|\, 
	\exists\,\widebar{\varphi}\in\Hom(PA,QA),\ \\
	\widebar{\varphi}|_P=\varphi,\ \varphi|_A=\Id_A,\ 
	\widebar\varphi/A\in\Hom_{\widebar{\calf}}(PA/A,QA/A) \bigr\}. 
	\end{multline*}
Then $A$ is a central subgroup of $\calf$, and $\calf/A=\widebar\calf$.  
But $\calf$ is not saturated, since $\Aut_{\calf}(S)$ contains as normal 
subgroup the group of automorphisms of $S$ which are the identity on $A$ 
and on $S/A$ (a $p$-group by Lemma \ref{gorenstein}), and this subgroup is 
by assumption not contained in $\Inn(S)$.  Hence $\Aut_{S}(S)$ is not a 
Sylow subgroup of $\autf(S)$.
\end{Ex}

The hypotheses of the above example are satisfied, for example, by any 
pair of $p$-groups $1\ne{}A\lneqq{}S$ with $A\le{}Z(S)$, such that $S$ is 
abelian, or more generally such that $A\cap[S,S]=1$.


We will now describe conditions which allow us to say when $\calf$ is 
saturated.  Before doing so, in the next two lemmas, we first compare 
properties of subgroups in $\calf$ with those of subgroups of $\calf/A$, 
when $\calf$ is a fusion system with central subgroup $A$.  This is done 
under varying assumptions as to whether $\calf$ or $\calf/A$ is saturated.

\begin{Lem}\label{centric-radical}
The following hold for any fusion system $\calf$ over a $p$-group $S$, and 
any subgroup $A\le{}Z(S)$ central in $\calf$.
\begin{enumerate}
\item If $A\le{}P\le{}S$ and $P/A$ is $\calf/A$-centric, then $P$ is 
$\calf$-centric.

\item If $\calf/A$ is saturated, and $P\le{}S$ is $\calf$-quasicentric, 
then $PA/A$ is $\calf/A$-quasicentric.  

\item If $\calf$ and $\calf/A$ are both saturated and $P/A\le{}S/A$ is 
$\calf/A$-quasicentric, then $P$ is $\calf$-quasicentric.
\end{enumerate}
\end{Lem}

\begin{proof} For each $P\le{}S$ containing $A$, let 
$\widetilde{C}_S(P)\le{}S$ be the subgroup such that 
$\widetilde{C}_S(P)/A=C_{S/A}(P/A)$.  Let 
	$$ \eta_P\: \widetilde{C}_S(P) \Right4{} \Hom(P,A) $$
be the homomorphism $\eta_P(x)(g)=[x,g]$ for $x\in\widetilde{C}_S(P)$ and 
$g\in{}P$.  Thus $\Ker(\eta_P)=C_S(P)$.  

For each $P\le{}S$ containing $A$, set
	$$ \Gamma_P = \Ker\bigl[\autf(P) \Right3{} 
	\Aut_{\calf/A}(P/A)\bigr]. $$
Since every $\alpha\in\autf(P)$ is the identity on $A$, $\Gamma_P$ is a 
$p$-group by Lemma \ref{gorenstein}.  

\smallskip

\noindent\textbf{(a) } If $g,g'$ commute in $S$, then their images commute 
in $S/A$.  Thus for all $Q\le{}S$, $C_S(Q)/A \le{}C_{S/A}(Q/A)$.  

Assume $P/A$ is $\calf/A$-centric.  Then for each $P'$ which is 
$\calf$-conjugate to $P$, $P'/A$ is $\calf/A$-conjugate to $P/A$, and 
hence $C_S(P')/A\le{}C_{S/A}(P'/A)\le{}P'/A$.  Thus $P$ is $\calf$-centric.

\smallskip

\noindent\textbf{(b) }  Fix $P\le{}S$ such that $P$ is 
$\calf$-quasicentric.  Since $C_S(P'A)=C_S(P')$ for all $P'$ which is 
$\calf$-conjugate to $P$, $P'A$ is fully centralized in $\calf$ if and 
only if $P'$ is.  Also, $C_\calf(P')=C_\calf(P'A)$ for such $P'$, and this 
shows that $PA$ is also $\calf$-quasicentric.  So after replacing $P$ by 
$PA$, if necessary, we can assume $P\ge{}A$.  

Choose $P'/A$ which is fully centralized in $\calf/A$ and 
$\calf/A$-conjugate to $P/A$.  Then $P'$ is $\calf$-conjugate to $P$, 
hence still $\calf$-quasicentric.  So upon replacing $P$ by $P'$, we are 
reduced to showing that $P/A$ is $\calf/A$-quasicentric when $P\ge{}A$, 
$P$ is $\calf$-quasicentric, and $P/A$ is fully centralized.

For any $P'$ which is $\calf$-conjugate to $P$, there is a morphism
	$$ \varphi/A\in\Hom_{\calf/A}\bigl((P'{\cdot}\widetilde{C}_S(P'))/A,
	(P{\cdot}\widetilde{C}_S(P))/A\bigr) $$
which sends $P'/A$ to $P/A$ (by axiom (II) for the saturated fusion system 
$\calf/A$).  Then $\varphi$ restricts to a morphism from $C_S(P')$ to 
$C_S(P)$.  Thus $|C_S(P')|\le|C_S(P)|$ for all $P'$ which is 
$\calf$-conjugate to $P$, and this proves that $P$ is fully centralized in 
$\calf$.

If $P/A$ is not $\calf/A$-quasicentric, then by Lemma \ref{L:qcentric} 
(and since $\calf/A$ is saturated), there is some $Q/A\le{} 
P/A{\cdot}C_{S/A}(P/A)$ containing $P/A$, and some $\alpha\in\autf(Q)$, 
such that $\Id\ne\alpha/A\in\Aut_{\calf/A}(Q/A)$ has order prime to $p$ 
and $\alpha/A$ is the identity on $P/A$.  Since $\alpha$ is also the 
identity on $A$, Lemma \ref{gorenstein} implies that $\alpha|_P=\Id_P$.  
Set $Q'=Q\cap\widetilde{C}_S(P)$.  Then $\alpha(Q')=Q'$, and 
$\eta_P\circ\alpha=\eta_P$ since $\alpha$ is the identity on $P\ge{}A$.  
Thus $\alpha$ induces the identity on $Q'/C_{Q'}(P)$ since 
$\Ker(\eta_P)=C_S(P)$.  Since $\alpha$ has order prime to $p$, Lemma 
\ref{gorenstein} now implies that $\alpha|_{C_{Q'}(P)}\ne\Id$.  Thus 
$C_{Q'}(P)\le{}C_S(P)$ and $\alpha|_{P{\cdot}C_{Q'}(P)}$ is a nontrivial 
automorphism in $C_\calf(P)$ of order prime to $p$.  Since $P$ is fully 
centralized in $\calf$, this implies (by definition) that $C_\calf(P)$ is 
not the fusion system of $C_S(P)$, and hence that $P$ is not
$\calf$-quasicentric.

\smallskip

\noindent\textbf{(c) }  Now assume that $\calf$ and $\calf/A$ are both 
saturated, and that $P/A\le{}S/A$ is $\calf/A$-quasicentric.  If $P$ is 
not $\calf$-quasicentric, and $P'$ is $\calf$-conjugate to $P$ and fully 
centralized in $\calf$, then by Lemma \ref{L:qcentric}(b), there is some 
$P'\le{}Q\le{}P'{\cdot}C_S(P')$ and some $\Id\ne\alpha\in\autf(Q)$ such 
that $\alpha|_{P'}=\Id_{P'}$ and $\alpha$ has order prime to $p$.  Then 
$Q/A\le(P'/A){\cdot}C_{S/A}(P'/A)$, $\alpha/A\in\Aut_{\calf/A}(Q/A)$ also 
has order prime to $p$, and so $\alpha/A\ne\Id$ by Lemma \ref{gorenstein} 
again.  But by Lemma \ref{L:qcentric}(a), this contradicts the assumption 
that $P/A$ is $\calf/A$-quasicentric.
\end{proof}

In the next lemma, we compare conditions for being fully normalized in 
$\calf$ and in $\calf/A$.

\begin{Lem}\label{centric-radical1}
The following hold for any fusion system $\calf$ over a $p$-group $S$, and 
any subgroup $A\le{}Z(S)$ central in $\calf$.
\begin{enumerate}
\item Assume $\calf$ is saturated.  Then for all $P,Q\le{}S$ containing 
$A$ such that $P$ is fully normalized in $\calf$, if 
$\varphi,\varphi'\in\homf(P,Q)$ are such that $\varphi/A=\varphi'/A$, then 
$\varphi'=\varphi\circ{}c_x$ for some $x\in{}N_S(P)$ such that 
$xA\in{}C_{S/A}(P/A)$.  

\item Assume $\calf/A$ is saturated, and let $P,P'\le{}S$ be 
$\calf$-conjugate subgroups which contain $A$.  Then $P$ is fully 
normalized in $\calf$ if and only if $P/A$ is fully normalized in 
$\calf/A$.  Moreover, if $P$ is fully normalized in $\calf$, then there is 
$\psi\in\homf(N_S(P'),N_S(P))$ such that $\psi(P')=P$.
\end{enumerate}
\end{Lem}

\begin{proof}
\noindent\textbf{(a) } Assume $\calf$ is saturated, and fix 
$P,Q\le{}S$ containing $A$ such that $P$ is fully normalized in $\calf$.  
Let $\varphi,\varphi'\in\homf(P,Q)$ be such that $\varphi/A=\varphi'/A$.  
Then $\Im(\varphi)=\Im(\varphi')$.  Set $\alpha=\varphi^{-1}\circ 
\varphi'\in\autf(P)$; then $\varphi'=\varphi\circ\alpha$, and 
$\alpha/A=\Id_{P/A}$.  Thus 
	$$ \alpha\in\Ker[\autf(P)\Right2{}\Aut_{\calf/A}(P/A)], $$
which is a normal $p$-subgroup by Lemma \ref{gorenstein}. Since $P$ is fully 
normalized, $\Aut_S(P)\in\sylp{\autf(P)}$  and so $\alpha\in\Aut_S(P)$.  Thus 
$\alpha\in{}c_x$ for some $x\in{}N_S(P)$, 
and $xA\in{}C_{S/A}(P/A)$ since $\alpha/A=\Id_{P/A}$.  

\smallskip

\noindent\textbf{(b) }  We are now assuming that $\calf/A$ is saturated.  
Assume first that $P/A$ is fully normalized in $\calf/A$.  Then by Lemma 
\ref{N->N}, there is a morphism 
	\beq \varphi_0\in\Hom_{\calf/A}(N_{S/A}(P'/A),N_{S/A}(P/A)) 
	\tag{1} \eeq
such that $\varphi_0(P'/A)=P/A$, and this lifts to 
$\varphi\in\homf(N_S(P'),N_S(P))$ such that $\varphi(P')=P$.  Therefore 
$|N_S(P')|\le{}|N_S(P)|$ for any $P'$ which is $\calf$-conjugate to $P$ 
and $P$ is fully normalized in $\calf$.  This also proves that the last 
statement in (b) holds, once we know that $P/A$ is fully normalized.

Assume now that $P$ is fully normalized in $\calf$; we want to show that 
$P/A$ is fully normalized in $\calf/A$.  Fix $P'$ which is 
$\calf$-conjugate to $P$ and such that $P'/A$ is fully normalized in 
$\calf/A$.  By Lemma \ref{N->N} again, there exists
	$$ \varphi_0\in\Hom_{\calf/A}(N_{S/A}(P/A),N_{S/A}(P'/A)) $$
such that $\varphi_0(P/A)=P'/A$.  Then $\varphi_0=\varphi/A$ 
for some $\varphi\in\homf(N_S(P),N_S(P'))$, and $\varphi$ is an 
isomorphism since $P$ is fully normalized in $\calf$.  Thus $\varphi_0$ is 
an isomorphism, and hence $P/A$ is fully normalized in $\calf/A$. 

Thus if $P$ is fully normalized in $\calf$, then $P/A$ is also fully 
normalized, and we have already seen that the last statement in (b) holds 
in this case.
\end{proof}

We are now ready to give conditions under which we show that $\calf$ is 
saturated if $\calf/A$ is saturated.  As in \cite[\S2]{bcglo1}, for any 
fusion system $\calf$ over a $p$-group $S$, and any set $\calh$ of 
subgroups of $S$, we say that $\calf$ is $\calh$-generated if each 
morphism in $\calf$ is a composite of restrictions of morphisms between 
subgroups in $\calh$.

\begin{Prop}\label{prop-saturation-central}
Let $A$ be a central subgroup of a fusion system $\calf$ over a $p$-group 
$S$, such that $\calf/A$ is a saturated fusion system.  Let 
$\calh$ be any set of subgroups of $S$, closed under $\calf$-conjugacy and 
overgroups, which contains all $\calf$-centric subgroups of $S$.  Assume
\begin{enumerate} 
\item $\Ker\bigl[\autf(P)\Right2{}\Aut_{\calf/A}(P/A)\bigr]
\le\Aut_S(P)$ for each $P\in\calh$ which is fully normalized in $\calf$; and 
\item $\calf$ is $\calh$-generated.
\end{enumerate}
Then $\calf$ is saturated.
\end{Prop}

\begin{proof}  Let $\calh_0$ be the set of all $P\in\calh$ such that 
$P\ge{}A$.  Since $A$ is central, each morphism $\varphi\in\homf(P,Q)$ 
extends to some $\widebar{\varphi}\in\homf(PA,QA)$ such that 
$\widebar{\varphi}|_A=\Id_A$.  Thus $\calf$ is $\calh_0$-generated if it 
is $\calh$-generated.  So upon replacing $\calh$ by $\calh_0$, we can 
assume all subgroups in $\calh$ contain $A$.

By assumption, $\calh$ 
contains all $\calf$-centric subgroups of $S$.  So by Theorem 
\ref{centr->sat}(b), it suffices to check that axioms (I) and (II) hold 
for all $P\in\calh$.  For all $P\le{}S$ containing $A$, we write
	$$ \Gamma_P\defeq
	\Ker\bigl[\autf(P)\Right2{}\Aut_{\calf/A}(P/A)\bigr]. $$

\smallskip

\noindent\textbf{(I) }  Assume $P\in\calh$ is fully normalized in $\calf$. 
Then $P/A$ is fully normalized in $\calf/A$ by Lemma 
\ref{centric-radical1}(b).  By assumption, $\Gamma_P\le\Aut_S(P)$.  Hence 
	$$ [\autf(P):\Aut_S(P)] = [\Aut_{\calf/A}(P/A):\Aut_{S/A}(P/A)], $$
and $\Aut_S(P)\in\sylp{\autf(P)}$ since 
$\Aut_{S/A}(P/A)\in\sylp{\Aut_{\calf/A}(P/A)}$.  

Assume $P'\le{}S$ is $\calf$-conjugate to $P$ and fully centralized in 
$\calf$.  By Lemma \ref{centric-radical1}(b) again, there is 
$\psi\in\homf(N_S(P'),N_S(P))$ such that $\psi(P')=P$.  Then 
$\psi(C_S(P'))\le{}C_S(P)$, so $P$ is also fully centralized.

\smallskip

\noindent\textbf{(II) }  Assume $\varphi\in\homf(P,S)$ is such that 
$P\in\calh$ and $\varphi(P)$ is fully centralized.  Set 
	$$ N_\varphi = \{g\in{}N_S(P) \,|\, 
	\varphi{}c_g\varphi^{-1}\in\Aut_S(\varphi(P)) \} , $$
as usual.  Then $N_\varphi/A \le N_{\varphi/A}$.  

Assume first that $\varphi(P)$ is fully normalized in $\calf$.  Then 
$\varphi(P)/A$ is fully normalized in $\calf/A$ by Lemma 
\ref{centric-radical1}(b), and hence also fully centralized.  So by (II), 
applied to the saturated fusion system $\calf/A$, there is 
$\widehat{\varphi}\in\homf(N_{\varphi/A},S/A)$ such that 
$\widehat{\varphi}|_{(P/A)}=\varphi/A$. Let  
$\widebar{\varphi}\in\homf(N_{\varphi},S)$ be a lift of 
$\widehat{\varphi}$, then $(\widebar{\varphi}|_P)/A=\varphi/A$ and 
$\widebar{\varphi}(P)=\varphi(P)$. So there is $\alpha=\varphi \circ 
(\widebar{\varphi}|_P)^{-1}\in\Gamma_{\varphi(P)}\le\autf(\varphi(P))$ such 
that $\alpha\circ\widebar{\varphi}|_P=\varphi$.  By (a), there is 
$x\in{}N_S(\varphi(P))$ such that $\alpha=c_x$; then 
$c_x\circ\widebar{\varphi}$ lies in $\homf(N_\varphi,S)$ and extends 
$\varphi$.  

It remains to prove the general case.  Choose $P'$ which is fully 
normalized in $\calf$ and $\calf$-conjugate to $P$.  By Lemma 
\ref{centric-radical1}(b), there is 
$\psi\in\homf(N_S(\varphi(P)),N_S(P'))$ such that $\psi(\varphi(P))=P'$.  
Then $N_\varphi\le{}N_{\psi\varphi}$. Since $\psi\varphi(P)=P'$ is fully 
normalized, $\psi\varphi$ extends to some 
$\widebar{\psi}\in\homf(N_\varphi,N_S(P'))$.  We will show 
that $\Im(\widebar{\psi})\le\Im(\psi)$, and thus there is 
$\widebar{\varphi}\in\homf(N_\varphi,N_S(\varphi(P)))$ such that 
$\widebar{\varphi}|_P=\varphi$.

For each $g\in{}N_\varphi$, choose $x\in{}N_S(\varphi(P))$ such that 
$\varphi{}c_g\varphi^{-1}=c_x$; then we obtain 
$c_{\psi(x)}=\psi{}c_x\psi^{-1}=c_{\widebar{\psi}(g)}$, and so 
$\psi(x)\widebar{\psi}(g)^{-1}\in C_S(P')$.  Since $\varphi(P)$ is fully 
centralized in $\calf$, $\psi(C_S(\varphi(P)))=C_S(P')$, and thus 
$\widebar{\psi}(g)\in\Im(\psi)$.  
\end{proof}

For example, one can take as the set $\calh$ in Proposition 
\ref{prop-saturation-central} either the set of subgroups $P\le{}S$ 
containing $A$ such that $P/A$ is $\calf/A$-quasicentric (by Lemma 
\ref{centric-radical}(b)), or the set of $\calf$-centric subgroups of $S$. 

Now let $\SFL$ be a $p$-local finite group. One can also define a 
centralizer linking system $C_{\call}(A)$ when $A\le{}S$ is fully 
centralized \cite[Definition 2.4]{BLO2}.  Since this is always a linking 
system associated to $C_\calf(A)$, $A$ is central in $\call$ 
(i.e., $C_\call(A)=\call$) if and only if $A$ is central in $\calf$.  So 
from now on, by a central subgroup of $\SFL$, we just mean a 
subgroup $A\le{}Z(S)$ which is central in $\calf$.  

\begin{Defi} \label{D:L/A}
Let $\SFL$ be a $p$-local finite group with a central subgroup $A$. Define 
$\call/A$ to be the category with objects the subgroups $P/A$ for 
$P\in\Ob(\call)$ (i.e., such that $P$ is $\calf$-centric), 
and with morphism sets 
	$$ \Mor_{\call/A}(P/A,Q/A) = \Mor_{\call}(P,Q)/\delta_P(A). $$
Let $(\call/A)^c\subseteq\call/A$ be the full subcategory with object set 
the $\calf/A$-centric subgroups in $S/A$.  

Similarly, if $\callq$ is the 
associated quasicentric linking system, then define $\callq/A$ to be the 
category with objects the $\calf/A$-quasicentric subgroups of $S/A$ --- 
equivalently, the subgroups $PA/A$ for $P\in\Ob(\callq)$ --- and with 
morphisms 
	$$ \Mor_{\call^q/A}(P/A,Q/A) = \Mor_{\call^q}(P,Q)/\delta_P(A). $$
\end{Defi}

Note that by Lemma \ref{centric-radical}(a,c), for any $P/A\le{}S/A$, if 
$P/A$ is $\calf/A$-centric or $\calf/A$-quasicentric, then $P$ is 
$\calf$-centric or $\calf$-quasicentric, respectively.  Thus the 
categories $(\call/A)^c\subseteq\call/A$ and $\callq/A$ are well defined.

We are now ready to prove our main theorem about quotient fusion and 
linking systems.  

\begin{Thm} \label{P:L/A-centext}
Let $A$ be a central subgroup of a $p$-local finite group $\SFL$ with 
associated quasicentric linking system $\callq$.  Let 
$\callq_{\ge{}A}\subseteq\callq$ be the full subcategory with objects 
those $\calf$-quasicentric subgroups of $S$ which contain $A$.  Then the 
following hold.
\begin{enumerate}  
\item $(S/A,\calf/A,(\call/A)^c)$ is again a $p$-local finite group, and 
$\callq/A$ is a quasicentric linking system associated to $\calf/A$. 

\item The inclusions of categories induce homotopy 
equivalences $|\call|\simeq|\callq_{\ge{}A}|\simeq|\callq|$ and 
$|(\call/A)^c|\simeq|\call/A|\simeq|\callq/A|$.

\item The functor $\tau\:\callq\rTo\callq/A$, defined by $\tau(P)=PA/A$ 
and with the obvious maps on morphisms, induces principal fibration 
sequences 
	$$ BA \Right4{} |\call| \Right4{|\tau|} |\call/A| 
	\qquad\textup{and}\qquad
	BA \Right4{} |\callq_{\ge{}A}|
	\Right4{|\tau|} |\callq/A| $$
which remain principal fibration sequences after $p$-completion.
\end{enumerate}
\end{Thm}

\begin{proof}  \textbf{(a) } The first statement is shown in \cite[Lemma 
5.6]{BLO2} when $|A|=p$, and the general case follows by iteration.  So we 
need only prove that $\callq/A$ is a quasicentric linking system 
associated to $\calf/A$.  Axioms (B)$_q$, (C)$_q$, and (D)$_q$ for 
$\callq/A$ follow immediately from those axioms applied to $\callq$, so it 
remains only to prove (A)$_q$.  

Let $\calh$ be the set of subgroups $P\le{}S$ such that $A\le{}P$ and 
$P/A$ is $\calf/A$-quasicentric and fix $P,Q\in\calh$.  By construction, 
$C_{S/A}(P/A)$ acts freely on $\Mor_{\callq/A}(P/A,Q/A)$ and induces a 
surjection
	$$ \Mor_{\callq/A}(P/A,Q/A)\big/C_{S/A}(P/A) \Onto4{} 
	\Hom_{\calf/A}(P/A,Q/A). $$
We must show that this is a bijection whenever $P/A$ is fully centralized 
in $\calf/A$.  Since any other fully centralized subgroup in the same 
$\calf$-conjugacy class has centralizer of the same order, it suffices to 
show this when $P/A$ is fully normalized in $\calf/A$; or equivalently (by 
Lemma \ref{centric-radical1}(b)) when $P$ is fully normalized in $\calf$.  

Let $\widetilde{C}_S(P)\le{}S$ be the subgroup such that 
$\widetilde{C}_S(P)/A=C_{S/A}(P/A)$.  Fix any $\calf/A$-quasicentric 
subgroup $Q/A$, and consider the following sequence of maps
	$$ \Mor_{\callq}(P,Q) \Right4{} \homf(P,Q) \Right4{}
	\Hom_{\calf/A}(P/A,Q/A). $$
Since $\callq$ is the quasicentric linking system of $\calf$, by property 
(A)$_q$ the first map is the orbit map of the action of $C_S(P)$.  The 
second is the orbit map for the action of $\Aut_{\widetilde{C}_S(P)}(P)$ 
by Lemma \ref{centric-radical1}(a). Thus the composite is the orbit map 
for the free action of $\widetilde{C}_S(P)$.  It now follows that 
$\widetilde{C}_S(P)/A\cong{}C_{S/A}(P/A)$ acts freely on 
$\Mor_{\callq/A}(P/A,Q/A)=\Mor_{\callq}(P,Q)/A$ with orbit set 
$\Hom_{\calf/A}(P/A,Q/A)$. 

\smallskip
	
\noindent\textbf{(b) }  These homotopy equivalences are all special cases 
of Proposition \ref{L-props}(a).

\smallskip

\noindent\textbf{(c) } 
For any category $\calc$ and any $n\ge0$, let $\calc_n$ denote the set of 
$n$-simplices in the nerve of $\calc$; i.e., the set of composable 
$n$-tuples of morphisms.  For each $n\ge0$, the group $\calb(A)_n$ acts 
freely on $(\callq_{\ge{}A})_n$:  an element $(a_1,\dots,a_n)$ acts by 
composing the $i$-th component with $\delta_P(a_i)$ for appropriate $P$.  
This action commutes with the face and degeneracy maps, and its orbit set 
is $(\callq/A)_n$.  It follows that the projection of $|\callq_{\ge{}A}|$ 
onto $|\callq/A|$ (and of $|\call|$ onto $|\call/A|$) is a principal 
fibration with fiber the topological group $BA=|\calb(A)|$ (see, e.g., 
\cite[\S\S18--20]{May} or \cite[Corollary V.2.7]{GJ}). 

By the principal fibration lemma of Bousfield and Kan \cite[II.2.2]{BK}, 
these sequences are still principal fibration sequences after 
$p$-completion.  
\end{proof}

\newsub{Central extensions of $p$-local finite groups}  
We first make it more precise what we mean by this.  

\begin{Defi} \label{D:centext}
A \emph{central extension} of a (saturated) fusion system $\calf$ over a 
$p$-group $S$ by an abelian $p$-group $A$ consists of a (saturated) fusion 
system $\widetilde\calf$ over a $p$-group $\widetilde{S}$, together with 
an inclusion $A \le Z(\widetilde{S})$, such that $A$ is a central subgroup 
of $\widetilde{\calf}$, 
$\widetilde{S}/A\cong{}S$, and $\widetilde\calf/A \cong \calf$ as fusion 
systems over $S$. 

Similarly, a central extension of a $p$-local finite group $\SFL$ by an 
abelian group $A$ consists of a $p$-local finite group 
$(\widetilde{S},\widetilde{\calf},\widetilde{\call})$, together with an 
inclusion $A \le Z_{\widetilde\calf}(\widetilde{S})$, such that 
$(\widetilde S/A,\widetilde\calf/A, (\widetilde\call/A)^c)\cong \SFL$.
\end{Defi}

Extensions of categories were defined and studied by Georges Hoff in 
\cite{Hoff}, where he proves that they are classified by certain 
$\Ext$-groups. We will deal with one particular case of this. Hoff's 
theorem implies that an extension of categories of the type 
$\calb(A)\Right2{}\widetilde{\call}^q_{\ge{}A}\Right2{\tau}\callq$ is 
classified by an element in $\higherlim{\callq}2(A)$.  What this extension 
really means is that $\callq$ is a quotient category of 
$\widetilde{\call}^q_{\ge{}A}$, where each morphism set in 
$\widetilde{\call}^q_{\ge{}A}$ admits a free action of $A$ with orbit set 
the corresponding morphism set in $\callq$. Also, $\higherlim{}2(A)$ means 
the second derived functor of the limit of the constant functor which 
sends each object of $\callq$ to $A$ and each morphism to $\Id_A$. 

We regard $A$ as an \emph{additive} group.  Fix an element
$[\omega]\in\higherlim{\callq}2(A)$, where $\omega$ is a (reduced) 
$2$-cocyle. Thus $\omega$ is a function from pairs of composable 
morphisms in $\callq$ to $A$ such that $\omega(f,g)=0$ if $f$ or $g$ is an 
identity morphism, and such that for any triple $f,g,h$ of composable 
morphisms, the cocycle condition is satisfied
	$$ d\omega(f,g,h) \defeq
	\omega(g,h)-\omega(gf,h)+\omega(f,hg)-\omega(f,g)=0 .$$ 
	
Consider the composite $i\circ\delta_S:\mathcal B S\rightarrow \call^q$, 
where $\delta_S$ is induced by the distinguished morphism $S\rightarrow 
\Aut_{\call^q}(S)$ and $i$ is induced by the inclusion of 
$\Aut_{\call^q}(S)$ in $\call^q$ as the full subcategory with one object 
$S$. Then $\omega_S=(i\circ\delta_S)^*(\omega)$ is a $2$-cocycle defined 
on $S$ which classifies a central extension of $p$-groups
	$$ 1 \Right2{} A \Right4{} \widetilde{S} \Right4{\tau} S 
	\Right2{} 1\,. $$
Thus $\widetilde{S}=S\times{}A$, with group multiplication defined by 
	$$ (h,b){\cdot}(g,a)=(hg,b+a+\omega_S(g,h)), $$
and $\tau(g,a)=g$.  For each $\calf$-quasicentric subgroup $P\le{}S$, set 
$\widetilde{P}=\tau^{-1}(P)\le\widetilde{S}$.
	
Using this cocycle $\omega$, we can define a new category 
$\widetilde{\call}_0$ as follows.  There is one object 
$\widetilde{P}=\tau^{-1}(P)$ of $\widetilde{\call}_0$ for each object $P$ 
of $\callq$.  Morphism sets in $\widetilde{\call}_0$ are defined by 
setting
	$$ \Mor_{\widetilde{\call}_0}(\widetilde P, \widetilde Q)
	=\Mor_{\callq}(P,Q)\times A. $$
Composition in $\widetilde{\call}_0$ is defined by 
	$$ (g,b)\circ (f,a)=(gf,a+b+\omega(f,g)). $$ 
The associativity of this composition law follows since $d\omega=0$.  
Furthermore, if we chose another representative $\omega + d\mu$ 
where $\mu$ is a $1$-cochain, the categories obtained are isomorphic.
This construction comes together with a projection functor 
$\tau\:\widetilde\call_0\rTo\callq$.  

Finally, we define
	$$ \delta_{\widetilde{S}}\: S \Right5{} 
	\Aut_{\widetilde{\call}_0}(S) $$
by setting $\delta_{\widetilde{S}}(g,a)=(\delta_S(g),a)$.  This is 
a group homomorphism by definition of $\omega_S$.  

\begin{Prop}\label{lemconj}
Let $\SFL$ be a $p$-local finite group, and let $\callq$ be the associated 
quasicentric linking system.  Fix an abelian group $A$ and a reduced 
2-cocycle $\omega$ on $\callq$ with coefficients in $A$.  Let 
	$$ \widetilde{\call}_0 \Right5{\tau} \callq 
	\qquad\textup{and}\qquad
	\widetilde{S} \Right5{\tau} S $$
be the induced extensions of categories and of groups, with distinguished 
monomorphism $\delta_{\widetilde{S}}$ as defined above.  

Then there is a unique saturated fusion system $\widetilde{\calf}$ over 
$\widetilde{S}$ and a functor 
$\widetilde{\pi}\:\widetilde{\call}_0\rTo\widetilde{\calf}$, together with 
distinguished monomorphisms $\delta_{\widetilde{P}}\colon 
\widetilde{P}\rightarrow \Aut_{\widetilde\call_0}(P)$.  \mynote{B:  Fix 
this!!}

If $\widetilde{\call}\subseteq\widetilde{\call}_0$ denotes the 
full subcategory whose objects are the $\widetilde{\calf}$-centric 
subgroups of $\widetilde{S}$ then  
$(\widetilde{S},\widetilde{\calf},\widetilde{\call})$ is a 
$p$-local finite group with central subgroup $A\le{}Z(\widetilde{S})$, 
such that 
	\beq 
	(\widetilde{S}/A,\widetilde{\calf}/A,\widetilde{\call}_0/A) 
	\cong (S,\calf,\callq). \tag{1} \eeq
 Also, $\widetilde{\call}_0$ extends 
to a quasicentric linking system $\widetilde{\call}^q$ associated to 
$\widetilde{\calf}$.
\end{Prop}

\begin{proof} Assume that inclusion morphisms $\iota_P$ have been chosen 
in  the quasicentric linking system  $\callq$ associated to 
$(S,\calf,\call)$.  For each $\calf$-quasicentric $P\lneqq{}S$, define
	$$ \iota_{\widetilde{P}} = (\iota_P,0) \in 
	\Mor_{\widetilde{\call}_0}(\widetilde{P},\widetilde{S}). $$
Next, define the distinguished monomorphism
	$$ \delta_{\widetilde{P}}\: 
	\widetilde{P}{\cdot}C_{\widetilde{S}}(\widetilde{P}) \Right5{} 
	\Aut_{\widetilde\call_0}(\widetilde{P}) $$
to be the unique monomorphism such that the following square commutes for 
each $\widetilde{P}$ and each 
$q\in{}\widetilde{P}{\cdot}C_{\widetilde{S}}(\widetilde{P})$:
	\beq \begin{diagram}[w=40pt]
	\widetilde{P} & \rTo^{\iota_{\widetilde{P}}} & \widetilde{S} \\
	\dTo>{\delta_{\widetilde{P}}(q,a)} && 
	\dTo>{\delta_{\widetilde{S}}(q,a)} \\
	\widetilde{P} & \rTo^{\iota_{\widetilde{P}}} & \widetilde{S} 
	\rlap{\,.}
	\end{diagram} \tag{2} \eeq
More precisely, since $\delta_{\widetilde{S}}(q,a)=(\delta_S(q),a)$ by 
definition, this means that
	$$ \delta_{\widetilde{P}}(q,a)= 
	(\delta_P(q),a+\omega(\iota_P,\delta_S(q))-
	\omega(\delta_P(q),\iota_P)). $$

For each morphism 
$(f,a)\in\Mor_{\widetilde\call_0}(\widetilde{P},\widetilde{Q})$,
and each element $(q,b)\in\widetilde{P}$, there is a unique element 
$c\in{}A$ such that the following square commutes:
	\beq
	\begin{diagram}[w=40pt]
	\widetilde{P} & \rTo^{(f,a)} & \widetilde{Q} \\
	\dTo>{\delta_{\widetilde{P}}(q,b)} && 
	\dTo>{\delta_{\widetilde{Q}}(\pi(f)(q),c)} \\
	\widetilde{P} & \rTo^{(f,a)} & \widetilde{Q} 
	\rlap{\,.}
	\end{diagram} \tag{3} \eeq
In this situation, we set 
	$$ \widetilde{\pi}(f,a)(q,b) = (\pi(f)(q),c)\in\widetilde{Q}. $$
By juxtaposing squares of the form (3), we see that 
$\widetilde{\pi}(f,a)\in\Hom(\widetilde{P},\widetilde{Q})$, and that this 
defines a functor $\widetilde{\pi}$ from $\widetilde\call_0$ to the 
category of subgroups of $\widetilde{S}$ with monomorphisms.  

Define $\widetilde{\calf}$ to be the fusion system over $\widetilde{S}$ 
generated by the image of $\widetilde{\pi}$ and restrictions.  By 
construction, the surjection $\tau\:\widetilde{S}\Right2{}S$ induces a 
functor $\tau_*\:\widetilde{\calf}\Right2{}\calf$ between the fusion 
systems, which is surjective since $\calf$ is generated by restrictions of 
morphisms between $\calf$-quasicentric subgroups (Theorem 
\ref{Alp.fusion}(a)).  So we can identify $\calf$ with 
$\widetilde{\calf}/A$.  By Lemma \ref{centric-radical}(b), for each 
$\widetilde{\calf}$-quasicentric subgroup $P\le\widetilde{S}$, $PA/A$ is 
$\calf$-quasicentric.  So we can extend $\widetilde{\call}_0$ to a 
category $\widetilde{\call}^q$ defined on all 
$\widetilde{\calf}$-quasicentric subgroups, by setting
	$$ \Mor_{\callq}(P,Q) = \bigl\{f\in\Mor_{\widetilde\call_0}(PA,QA) 
	\,\big|\, \widetilde{\pi}(f)(P)\le Q \bigr\}
	\quad\textup{and}\quad \delta_P=\delta_{PA} $$
for each pair $P,Q$ of $\widetilde{\calf}$-quasicentric subgroups.  
We extend $\widetilde{\pi}$ to $\widetilde{\call}^q$ in the obvious way.

It remains to prove that $\widetilde{\calf}$ is saturated, and 
that $\widetilde{\call}^q$ is a quasicentric linking system associated to 
$\widetilde{\calf}$.  In this process of doing this, we will also 
prove the isomorphism (1).  

Let $\calh$ be the set of subgroups 
$\widetilde{P}=\tau^{-1}(P)\le\widetilde{S}$ for all $\calf$-quasicentric 
subgroups $P\le{}S$.  

\smallskip

\noindent\textbf{ $\widetilde{\calf}$ is saturated:} 
We want to apply Proposition \ref{prop-saturation-central} to prove that 
$\widetilde{\calf}$ is saturated. By Lemma \ref{centric-radical}(b), 
$\calh$ contains all $\widetilde{\calf}$-quasicentric subgroups of 
$\widetilde{S}$ which contain $A$, and in particular, all 
$\widetilde{\calf}$-centric subgroups (since every 
$\widetilde{\calf}$-centric subgroup of $\widetilde{S}$ must contain $A$).  
Since $\widetilde{\calf}$ is $\calh$-generated by construction, it remains 
only to prove condition (a) in Proposition \ref{prop-saturation-central}.

Fix some $\widetilde{P}=\tau^{-1}(P)$ in $\calh$, and let 
$\varphi\in\Aut_{\widetilde{\calf}}(\widetilde{P})$ be such that 
$\tau_*(\varphi)=\Id_P$.  Choose $(f,a)\in\widetilde{\pi}^{-1}(\varphi)$; 
thus 
	$$ \varphi=\widetilde{\pi}(f,a) \in 
	\Ker\bigl[\Aut_{\widetilde{\calf}}(\widetilde{P}) 
	\Right3{} \autf(P) \bigr]. $$
Then $\pi(f)=\Id_P$, so $f=\delta_P(q)$ for some $q\in{}C_S(P)$, and
$(f,a)=\delta_{\widetilde{P}}(q,c)$ for some $c\in A$.  Since 
$\delta_{\widetilde{P}}$ is a 
homomorphism, the definition of $\varphi=\widetilde{\pi}(f,a)$ via (3) 
shows that $\varphi=\widetilde\pi(f,a)=\widetilde\delta_{\widetilde 
P}(q,c)=c_{(q,c)}$, and thus that 
$\varphi\in\Aut_{\widetilde{S}}(\widetilde{P})$.  Thus condition (a) in 
Proposition \ref{prop-saturation-central} holds, and this finishes the 
proof that $\widetilde{\calf}$ is saturated.

\smallskip

\noindent\textbf{$\widetilde{\call}^q$ is a quasicentric linking system 
associated to $\widetilde{\calf}$:}  The distinguished monomorphisms 
$\delta_{\widetilde P}$, for $\widetilde P\in\Ob(\widetilde{\call}^q)$, 
were chosen so as to satisfy (D)$_q$, and this was independent of the 
choice of inclusion morphisms which lift the chosen inclusion morphisms in 
$\callq$.  Once the $\delta_P$ were determined, then $\widetilde{\pi}$ was 
defined to satisfy (C)$_q$, and $\widetilde{\calf}$ was defined as the 
category generated by $\Im(\widetilde{\pi})$ and restrictions.  Thus all 
of these structures were uniquely determined by the starting data.
Axiom (B)$_q$ follows from (C)$_q$ by Lemma \ref{C=>B}.  

It remains only to prove (A)$_q$.  We have already seen that the functor 
$\widetilde{\pi}$ is surjective on all morphism sets.  Also, since 
$C_{\widetilde{S}}(P)=C_{\widetilde{S}}(PA)$ for all $P\le\widetilde{S}$, 
it suffices to prove (A)$_q$ for morphisms between subgroups 
$\widetilde{P}=\tau^{-1}(P)$ and $\widetilde{Q}=\tau^{-1}(Q)$ containing 
$A$.  By construction, $C_{\widetilde{S}}(\widetilde{P})$ acts freely on 
each morphism set 
$\Mor_{\widetilde{\call}^q}(\widetilde{P},\widetilde{Q})$, and it remains 
to show that if $\widetilde{P}$ is fully centralized, then 
$\widetilde{\pi}_{\widetilde{P},\widetilde{Q}}$ is the orbit map of this 
action.  As in the proof of Theorem \ref{P:L/A-centext}(a), it suffices to 
do this when $\widetilde{P}$ and $P$ are fully normalized.  

Fix two morphisms $(f,a)$ and $(g,b)$ from $\widetilde{P}$ to 
$\widetilde{Q}$ such that $\widetilde{\pi}(f,a)=\widetilde{\pi}(g,b)$.  
Then $\pi(f)=\pi(g)$, so $g=f\circ\delta_P(x)$ for some $x\in{}C_S(P)$, and
	$$ (g,b) = (f,a) \circ (\delta_P(x),b-a-\omega(\delta_P(x),f)) 
	= (f,a) \circ \delta_{\widetilde{P}}(x,c), $$
where $c=b-a-\omega(\delta_P(x),f)- 
(\omega(\iota,\delta_S(x))-\omega(\delta_S(x),\iota_P))$.  Also, 
$\widetilde{\pi}(\delta_{\widetilde{P}}(x,c))=\Id_{\widetilde{P}}$ implies 
(via (3)) that $(x,c)$ commutes with all elements of $\widetilde{P}$, so 
$(x,c)\in{}C_{\widetilde{S}}(\widetilde{P})$, and this finishes the 
argument.  
\end{proof}

We now want to relate the obstruction theory for central extensions of 
linking systems with kernel $A$ to those for central extensions of 
$p$-groups, and to those for principal fibrations with fiber $BA$.  As a 
consequence of this, we will show (in Theorem \ref{centthm'}) that when 
appropriate restrictions are added, these three types of extensions are 
equivalent.

Given a central extension $\widetilde{\calf}$ of $\calf$ by the central 
subgroup $A$, there is an induced central extension $1\to 
A\to\widetilde{S}\to S\to1$ of Sylow subgroups. Restriction to subgroups 
$P\le{}S$ produces corresponding central extensions $1\to 
A\to\widetilde{P}\to P\to1$. The homology classes of these central 
extensions are all compatible with morphisms from the fusion system, and 
hence define an element in $\higherlim{\calf}{}H^2(-;A)$.  This, together 
with notation already used in \cite{BLO2}, motivates the following 
definition:

\begin{Defi} \label{D:H*F}
For any saturated fusion system over a $p$-group $S$, and any finite 
abelian $p$-group $A$, define
	$$ H^*(\calf;A) = \higherlim{\calf}{}H^*(-;A) \cong
	\higherlim{\calfc}{}H^*(-;A). $$
\end{Defi}

The following lemma describes the relation between the cohomology of 
$\calf$ and the cohomology of the geometric realization of any linking 
system associated to $\calf$.

\begin{Lem} \label{H(L)=H(F)}
For any $p$-local finite group $\SFL$, and any finite abelian $p$-group 
$A$, the natural homomorphism
	\beq H^*(|\call|;A) \Right5{} H^*(\calf;A), \tag{1} \eeq
induced by the inclusion of $BS$ into $|\call|$, is an isomorphism.  
Furthermore, there are natural isomorphisms
	\beq \higherlim{\callq}2(A) \cong H^2(|\callq|;A) \cong 
	H^2(|\call|;A). \tag{2} \eeq
\end{Lem}

\begin{proof}  The second isomorphism in (2) follows from Proposition 
\ref{L-props}(a).  The first isomorphism holds for higher limits of any 
constant functor over any small, discrete category $\calc$, since both 
groups are cohomology groups of the same cochain complex
	$$ 0 \Right2{} \prod_{c}A \Right3{} \prod_{c_0\to c_1}A 
	\Right3{} \prod_{c_0\to c_1\to c_2}A \Right3{} \cdots\,. $$
This cochain complex for higher limits is shown in \cite[Appendix II, 
Proposition 3.3]{GZ} (applied with $\mathcal{M}=\Ab\op$).  

To prove the isomorphism (1), it suffices to consider the case where 
$A=\Z/p^n$ for some $n$.  This was shown in \cite[Theorem 5.8]{BLO2} when 
$A=\Z/p$, so we can assume that $n\ge2$, and that the lemma holds when 
$A=\Z/p^{n-1}$.  Consider the following diagram of Bockstein exact 
sequences
	\beq \begin{diagram}
	\rTo^{\ } & H^{i+1}(|\call|;\Z/p) & \rTo & H^i(|\call|;\Z/p^{n-1}) & 
	\rTo & H^{i}(|\call|;\Z/p^n) & \rTo & H^i(|\call|;\Z/p) & \rTo \\
	& \dTo && \dTo && \dTo && \dTo \\
	\rTo^{\ } & H^{i+1}(BS;\Z/p) & \rTo & H^i(BS;\Z/p^{n-1}) & 
	\rTo & H^{i}(BS;\Z/p^n) & \rTo & H^i(BS;\Z/p) & \rTo \rlap{\,.}
	\end{diagram} \tag{1} \eeq
We claim that the bottom row restricts to an exact sequence of groups 
$H^*(\calf;-)$.  Once this is shown, the result follows by the 5-lemma.

By \cite[Proposition 5.5]{BLO2}, there is a certain $(S,S)$-biset $\Omega$ 
which induces, via a sum of composites of transfer maps and maps induced 
by homomorphisms, an idempotent endomorphism of $H^*(BS;\Z/p)$ whose image 
is $H^*(\calf;\Z/p)$.  This biset $\Omega$ also induces endomorphisms 
$[\Omega]$ of $H^*(BS;\Z/p^n)$ and $H^*(BS;\Z/p^{n-1})$ which commute with 
the bottom row in (1), since any exact sequence induced by a short exact 
sequence of coefficient groups will commute with transfer maps and maps 
induced by homomorphisms.  The same argument as that used in the proof of 
\cite[Proposition 5.5]{BLO2} shows that in all of these cases, 
$\Im([\Omega])=H^*(\calf;-)$, and the restriction of $[\Omega]$ to its 
image is multiplication by $|\Omega|/|S|\in1+p\Z$.  Thus the sequence of 
the $H^*(\calf;-)$ splits as a direct summand of the bottom row in (1), 
and hence is exact.
\end{proof}

We can now collect the results about central extensions of $\SFL$ in the 
following theorem, which is the analog for $p$-local finite groups of the 
classical classification of central extensions of groups. 

\begin{Thm}\label{centthm'}
Let $\SFL$ be a $p$-local finite group. For each finite abelian $p$-group 
$A$, the following three sets are in one-to-one correspondence:
\begin{enumerate} 
\item equivalence classes of central extensions of $\SFL$ by $A$;
\item equivalence classes of principal fibrations $ BA \rTo X \rTo 
|\call|\pcom $; and
\item isomorphism classes of central extensions $1\rTo A \rTo \widetilde{S} 
\rTo^{\tau} S \rTo 1$ for which each morphism $\varphi\in\homf(P,Q)$ lifts 
to some $\widetilde{\varphi}\in\Hom(\tau^{-1}(P),\tau^{-1}(Q))$.
\end{enumerate}
The equivalence between the first two is induced by taking classifying 
spaces, and the equivalence between (a) and (c) is induced by restriction 
to the underlying $p$-group.  These sets are all in natural one-to-one 
correspondence with 
	\beq \higherlim{\call}2(A) \cong H^2(|\call|;A) \cong 
	H^2(\calf;A). \tag{1} \eeq
\end{Thm}

\begin{proof} The three groups in (1) are isomorphic by Lemma 
\ref{H(L)=H(F)}.  By Proposition \ref{lemconj}, central extensions of 
$\call$ are classified by $\higherlim{\call}2(A)$.  Principal fibrations 
over $|\call|\pcom$ with fiber $BA$ are classified by 
	$$ \bigl[|\call|\pcom,B(BA)\bigr] \cong H^2(|\call|\pcom;A)\cong 
	H^2(|\call|;A), $$
where the second isomorphism holds since $|\call|$ is $p$-good 
(\cite[Proposition 1.12]{BLO2}) and $A$ is 
an abelian $p$-group.  Central extensions of $S$ by $A$ are classified by 
$H^2(S;A)$, and the central extension satisfies the condition in (c) if 
and only if the corresponding element of $H^2(S;A)$ extends to an element 
in the inverse limit $H^2(\calf;A)$.  

A central extension of $p$-local finite groups induces a principal 
fibration of classifying spaces by Theorem \ref{P:L/A-centext}(c), and 
this principal fibration restricts to the principal fibration over $BS$ of 
classifying spaces of $p$-groups.  Thus, if the fibration over $|\call|$ 
is classified by $\chi\in H^2(|\call|;A)$, then the fibration over $BS$ is 
classified by the restriction of $\chi$ to $H^2(BS;A)$.  It is well known
that the invariant in $H^2(S;A)=H^2(BS;A)$ for a central extension of $S$ 
by $A$ is the same as that which describes the principal fibration over 
$BS$ with fiber $BA$ (see \cite[Lemma IV.1.12]{AM}).  Since 
$H^2(|\call|;A)$ injects into $H^2(\calf;A)$ (Lemma \ref{H(L)=H(F)}), this 
shows that $\chi$ is also the class of the 2-cocycle which defines the 
extension of categories.  So the map between the sets in (a) and (b) 
defined by taking geometric realization is equal to the bijection defined 
by the obstruction theory.  

Since the isomorphism $\higherlim{\call}2(A)\cong{}H^2(\calf;A)$ is 
defined by restriction to $S$ (as a group of automorphisms in $\call$), 
the bijection between (a) and (c) induced by the bijection of obstruction 
groups is the same as that induced by restriction to $S$.  
\end{proof}

The following corollary shows that all minimal examples of ``exotic'' 
fusion systems have trivial center.


\begin{Cor} \label{F.exo=>F/A.exo}
Let $\calf$ be a saturated fusion system over a $p$-group $S$, and assume 
there is a nontrivial subgroup $1\ne{}A\le{}Z(S)$ which is central in 
$\calf$.  Then $\calf$ is the fusion system of some finite group if and 
only if $\calf/A$ is.
\end{Cor}

\begin{proof}  Assume $\calf$ is the fusion system of the finite group 
$G$, with $S\in\sylp{G}$.  By assumption ($A$ is central in $\calf$), each 
morphism in $\calf$ extends to a morphism between subgroups containing $A$ 
which is the identity on $A$.  Hence $\calf$ is also the fusion system of 
$C_G(A)$ over $S$, and so $\calf/A$ is the fusion system of $C_G(A)/A$.

It remains to prove the converse.  Assume $\calf/A$ is isomorphic to the 
fusion system of the finite group $\widebar{G}$, and identify $S/A$ with a 
Sylow $p$-subgroup of $\widebar{G}$.  Since by Lemma \ref{H(L)=H(F)}, 
$H^2(\calf/A;A)\cong H^2(B\widebar{G};A)$ and A is $\calf$-central, the 
cocycle classifying the extension is in $H^2(B\widebar{G};A)$, and hence 
there is an extension of finite groups
	$$ 1 \Right3{} A \Right4{} G \Right4{\tau} 
	\widebar{G} \Right3{} 1 $$
with the same obstruction invariant as the extension $\calf\rTo\calf/A$.  
In particular, we can identify $S=\tau^{-1}(S/A)\in\sylp{G}$, and 
$\widebar{G}=G/A$.  

We will prove the following two statements:
\begin{enumerate} 
\item $\calf$ has an associated centric linking system $\call$; and
\item $\call_{S/A}^c(G/A)$ is the unique centric linking system associated 
to $\calf/A=\calf_{S/A}(G/A)$.
\end{enumerate}
Once these have been shown, then they imply that
	$$ (S/A,\calf/A,(\call/A)^c)\cong
	(S/A,\calf_{S/A}(G/A),\call_{S/A}^c(G/A)) $$
as $p$-local finite groups.  Hence $\SFL\cong(S,\calf_S(G),\call_S^c(G))$ 
as $p$-local finite groups by Theorem \ref{centthm'}, and thus $\calf$ is 
the fusion system of $G$. 

It remains to prove (a) and (b).  Let 
	$$ \calz_\calf\:\orb^c(\calf) \Right3{} \zploc\mod
	\qquad\textup{and}\qquad
	\calz_G\:\orb^c_S(G) \Right3{} \zploc\mod $$
be the categories and functors of Definition \ref{D:orbit}.  By 
\cite[Lemma 2.1]{Oliver-odd}, $\calz_G$ can also be regarded as a functor 
on $\orb^c(\calf_S(G))$, and 
\begin{enumerate}\setcounter{enumi}{2}
\item $\higherlim{\orb_S^c(G)}*(\calz_G)\cong
\higherlim{\orb^c(\calf_S(G))}*(\calz_G)$.  
\end{enumerate}

By \cite[Proposition 3.1]{BLO2}, the existence and uniqueness of a centric 
linking system depends on the vanishing of certain obstruction classes:  
the obstruction to existence lies in 
$\higherlim{\orb^c_{S}(G)}3(\calz_\calf)$ and the obstruction to uniqueness 
in $\higherlim{\orb^c_{S}(G)}2(\calz_\calf)$.  Thus 
(b) follows from \cite[Proposition 3.1]{BLO2} and (c), once we know that 
$\higherlim{}2(\calz_{G/A})=0$; and this is shown in \cite[Theorem 
A]{Oliver-odd} (if $p$ is odd) or \cite[Theorem A]{Oliver-2} (if $p=2$). 

It remains to prove point (a), and we will do this by showing that 
	\beq \higherlim{\orb^c_{}(\calf)}3(\calz_\calf)\cong 
	\higherlim{\orb^c_{S}(G)}3(\calz_G)=0. \tag{1} \eeq
The last equality follows from \cite[Theorem A]{Oliver-odd} or \cite[Theorem 
A]{Oliver-2} again, so it remains only to prove the isomorphism. 

Let $\calh$ be the set of subgroups $P\le{}S$ containing $A$ such that 
$P/A$ is $\calf/A$-centric; or equivalently, $p$-centric in $G/A$.  Let 
$\orb^\calh(\calf)\subseteq\orb^c(\calf)$ and 
$\orb^\calh(\calf_S(G))\subseteq\orb^c(\calf_S(G))$ be the full 
subcategories with object sets $\calh$.

We claim that 
	\beq \higherlim{\orb^\calh(\calf)}*(\calz_\calf)\cong
	\higherlim{\orb^c(\calf)}*(\calz_\calf)
	\qquad\textup{and}\qquad
	\higherlim{\orb^\calh(\calf_S(G))}*(\calz_G)\cong
	\higherlim{\orb^c(\calf_S(G))}*(\calz_G) \,. \tag{2} \eeq
If $P\le{}S$ is $\calf$-centric but not in $\calh$, then there is 
$x\in{}N_S(P){\sminus}P$ such that $c_x\in\autf(P)$ is the identity on 
$P/A$ and on $A$, but is not in $\Inn(P)$.  Thus 
$1\ne[c_x]\in{}O^p(\outf(P))$, and so $P$ is not $\calf$-radical.  
By \cite[Proposition 3.2]{BLO2}, if $F\:\orb^c(\calf)\op\rTo\zploc\mod$ is 
any functor which vanishes except on the conjugacy class of $P$, then 
$\higherlim{}*(F)\cong\Lambda^*(\outf(P);F(P))$, where $\Lambda^*$ is a 
certain functor defined in \cite[\S5]{JMO}.  By \cite[Proposition 
6.1(ii)]{JMO}, $\Lambda^*(\outf(P);F(P))=0$ for any $F$, since 
$O_p(\outf(P))\ne1$.  From the long exact sequences of higher limits 
induced by extension of functors, it now follows that 
$\higherlim{\orb^c(\calf)}*(F)=0$ for any functor $F$ on $\orb^c(\calf)$ 
which vanishes on $\orb^\calh(\calf)$; and this proves the first 
isomorphism in (2).  The second isomorphism follows by a similar argument.

The natural surjection of $\calf$ onto $\calf/A$ induces isomorphisms of 
categories
	\beq \orb^\calh(\calf) \cong \orb^c(\calf/A) 
	\cong \orb^\calh(\calf_S(G))\,. \tag{3} \eeq
By the definition of $\calh$, we clearly have bijections between the sets 
of objects in these categories, so it remains only to compare morphism 
sets.  The result follows from Lemma \ref{centric-radical1}(a) for sets of 
morphisms between pairs of fully normalized subgroups, and the general 
case follows since every object is isomorphic to one which is fully 
normalized. 

We next claim that for all $i>0$, 
	\beq \higherlim{\orb^{\calh}(\calf)}i(\calz_\calf) 
	\cong \higherlim{\orb^c(\calf/A)}i(\calz_\calf/A)
	\qquad\textup{and}\qquad
	\higherlim{\orb^{\calh}(\calf_S(G))}i(\calz_G) \cong
	\higherlim{\orb^c(\calf/A)}i(\calz_G/A) \,. \tag{4} \eeq
To show this, by (3), together with the long exact sequence of higher 
limits induced by the short exact sequence of functors 
	$$ 1\rTo A \Right4{} \calz_\calf \Right4{} \calz_\calf/A\rTo 1, $$ 
we need only show that $\higherlim{\orb^c(\calf/A)}i(A)=0$ for all 
$i>0$.  Here, $A$ denotes the constant functor on $\orb^c(\calf/A)$ which 
sends all objects to $A$ and all morphisms to $\Id_A$.  For each 
$P\in\calh$, let $F_{A,P}$ be the functor on $\orb^c(\calf/A)$ where 
$F_{A,P}(P'/A)=A$ if $P'$ is $\calf$-conjugate to $P$ and 
$F_{A,P}(P'/A)=0$ otherwise; and which sends isomorphisms between 
subgroups conjugate to $P$ to $\Id_A$.  By \cite[Proposition 3.2]{BLO2}, 
$\higherlim{}*(F_{A,P})\cong\Lambda^*(\Out_{\calf/A}(P/A);A)$.  Also, 
$\Lambda^i(\Out_{\calf/A}(P/A);A)=0$ for $i>0$ by \cite[Proposition 
6.1(i,ii)]{JMO}, since the action of $\Out_{\calf/A}(P/A)$ on $A$ is 
trivial.  From the long exact sequences of higher limits induced by 
extension of functors, we now see that $\higherlim{}i(A)=0$ for all $i>0$. 

Finally, we claim that 
	\beq \calz_\calf/A\cong\calz_G/A \tag{5} \eeq
as functors on $\orb^c(\calf/A)$ under the identifications in (3).  To
see this, note that for each $P$, $(\calz_\calf/A)(P)=Z(P)/A= 
(\calz_G/A)(P)$; and since these are both identified as subgroups of 
$P/A$, any morphism in $\orb^c(\calf/A)$ from $P/A$ to $Q/A$ induces the 
same map (under the two functors) from $Z(Q)/A$ to $Z(P)/A$.  (This 
argument does \emph{not} apply to prove that $\calz_\calf\cong\calz_G$.  
These two functors send each object to the same subgroup of $S$, but we do 
not know that they send each morphism to the same homomorphism between the 
subgroups.) 

Thus by (2), (4), and (5), for all $i>0$, 
	\begin{align*} 
	\higherlim{\orb^{c}(\calf)}i(\calz_\calf) \cong 
	\higherlim{\orb^{\calh}(\calf)}i(\calz_\calf) 
	&\cong \higherlim{\orb^c(\calf/A)}i(\calz_\calf/A) \\
	&\cong \higherlim{\orb^c(\calf/A)}i(\calz_G/A) \cong 
	\higherlim{\orb^{\calh}(\calf_S(G))}i(\calz_G) \cong
	\higherlim{\orb^c_S(G)}i(\calz_G). 
	\end{align*}
This finishes the proof of (1), and hence finishes the proof of the 
corollary.
\end{proof}



\end{document}